\documentclass[aos,preprint]{imsart}

\RequirePackage[OT1]{fontenc}
\RequirePackage{amsthm,amsmath,natbib}
\RequirePackage[colorlinks,citecolor=blue,urlcolor=blue]{hyperref}

\usepackage{graphicx}
\usepackage{float}

\setattribute{journal}{name}{}

\startlocaldefs
\numberwithin{equation}{section}
\theoremstyle{plain}

\endlocaldefs

\usepackage{color}   
\definecolor{AliceBlue}{rgb}{0.94,0.97,1.00}
\definecolor{AntiqueWhite1}{rgb}{1.00,0.94,0.86}
\definecolor{AntiqueWhite2}{rgb}{0.93,0.87,0.80}
\definecolor{AntiqueWhite3}{rgb}{0.80,0.75,0.69}
\definecolor{AntiqueWhite4}{rgb}{0.55,0.51,0.47}
\definecolor{AntiqueWhite}{rgb}{0.98,0.92,0.84}
\definecolor{BlanchedAlmond}{rgb}{1.00,0.92,0.80}
\definecolor{BlueViolet}{rgb}{0.54,0.17,0.89}
\definecolor{CadetBlue1}{rgb}{0.60,0.96,1.00}
\definecolor{CadetBlue2}{rgb}{0.56,0.90,0.93}
\definecolor{CadetBlue3}{rgb}{0.48,0.77,0.80}
\definecolor{CadetBlue4}{rgb}{0.33,0.53,0.55}
\definecolor{CadetBlue}{rgb}{0.37,0.62,0.63}
\definecolor{CornflowerBlue}{rgb}{0.39,0.58,0.93}
\definecolor{DarkBlue}{rgb}{0.00,0.00,0.55}
\definecolor{DarkCyan}{rgb}{0.00,0.55,0.55}
\definecolor{DarkGoldenrod1}{rgb}{1.00,0.73,0.06}
\definecolor{DarkGoldenrod2}{rgb}{0.93,0.68,0.05}
\definecolor{DarkGoldenrod3}{rgb}{0.80,0.58,0.05}
\definecolor{DarkGoldenrod4}{rgb}{0.55,0.40,0.03}
\definecolor{DarkGoldenrod}{rgb}{0.72,0.53,0.04}
\definecolor{DarkGray}{rgb}{0.66,0.66,0.66}
\definecolor{DarkGreen}{rgb}{0.00,0.39,0.00}
\definecolor{DarkGrey}{rgb}{0.66,0.66,0.66}
\definecolor{DarkKhaki}{rgb}{0.74,0.72,0.42}
\definecolor{DarkMagenta}{rgb}{0.55,0.00,0.55}
\definecolor{DarkOliveGreen1}{rgb}{0.79,1.00,0.44}
\definecolor{DarkOliveGreen2}{rgb}{0.74,0.93,0.41}
\definecolor{DarkOliveGreen3}{rgb}{0.64,0.80,0.35}
\definecolor{DarkOliveGreen4}{rgb}{0.43,0.55,0.24}
\definecolor{DarkOliveGreen}{rgb}{0.33,0.42,0.18}
\definecolor{DarkOrange1}{rgb}{1.00,0.50,0.00}
\definecolor{DarkOrange2}{rgb}{0.93,0.46,0.00}
\definecolor{DarkOrange3}{rgb}{0.80,0.40,0.00}
\definecolor{DarkOrange4}{rgb}{0.55,0.27,0.00}
\definecolor{DarkOrange}{rgb}{1.00,0.55,0.00}
\definecolor{DarkOrchid1}{rgb}{0.75,0.24,1.00}
\definecolor{DarkOrchid2}{rgb}{0.70,0.23,0.93}
\definecolor{DarkOrchid3}{rgb}{0.60,0.20,0.80}
\definecolor{DarkOrchid4}{rgb}{0.41,0.13,0.55}
\definecolor{DarkOrchid}{rgb}{0.60,0.20,0.80}
\definecolor{DarkRed}{rgb}{0.55,0.00,0.00}
\definecolor{DarkSalmon}{rgb}{0.91,0.59,0.48}
\definecolor{DarkSeaGreen1}{rgb}{0.76,1.00,0.76}
\definecolor{DarkSeaGreen2}{rgb}{0.71,0.93,0.71}
\definecolor{DarkSeaGreen3}{rgb}{0.61,0.80,0.61}
\definecolor{DarkSeaGreen4}{rgb}{0.41,0.55,0.41}
\definecolor{DarkSeaGreen}{rgb}{0.56,0.74,0.56}
\definecolor{DarkSlateBlue}{rgb}{0.28,0.24,0.55}
\definecolor{DarkSlateGray1}{rgb}{0.59,1.00,1.00}
\definecolor{DarkSlateGray2}{rgb}{0.55,0.93,0.93}
\definecolor{DarkSlateGray3}{rgb}{0.47,0.80,0.80}
\definecolor{DarkSlateGray4}{rgb}{0.32,0.55,0.55}
\definecolor{DarkSlateGray}{rgb}{0.18,0.31,0.31}
\definecolor{DarkSlateGrey}{rgb}{0.18,0.31,0.31}
\definecolor{DarkTurquoise}{rgb}{0.00,0.81,0.82}
\definecolor{DarkViolet}{rgb}{0.58,0.00,0.83}
\definecolor{DeepPink1}{rgb}{1.00,0.08,0.58}
\definecolor{DeepPink2}{rgb}{0.93,0.07,0.54}
\definecolor{DeepPink3}{rgb}{0.80,0.06,0.46}
\definecolor{DeepPink4}{rgb}{0.55,0.04,0.31}
\definecolor{DeepPink}{rgb}{1.00,0.08,0.58}
\definecolor{DeepSkyBlue1}{rgb}{0.00,0.75,1.00}
\definecolor{DeepSkyBlue2}{rgb}{0.00,0.70,0.93}
\definecolor{DeepSkyBlue3}{rgb}{0.00,0.60,0.80}
\definecolor{DeepSkyBlue4}{rgb}{0.00,0.41,0.55}
\definecolor{DeepSkyBlue}{rgb}{0.00,0.75,1.00}
\definecolor{DimGray}{rgb}{0.41,0.41,0.41}
\definecolor{DimGrey}{rgb}{0.41,0.41,0.41}
\definecolor{DodgerBlue1}{rgb}{0.12,0.56,1.00}
\definecolor{DodgerBlue2}{rgb}{0.11,0.53,0.93}
\definecolor{DodgerBlue3}{rgb}{0.09,0.45,0.80}
\definecolor{DodgerBlue4}{rgb}{0.06,0.31,0.55}
\definecolor{DodgerBlue}{rgb}{0.12,0.56,1.00}
\definecolor{FloralWhite}{rgb}{1.00,0.98,0.94}
\definecolor{ForestGreen}{rgb}{0.13,0.55,0.13}
\definecolor{GhostWhite}{rgb}{0.97,0.97,1.00}
\definecolor{GreenYellow}{rgb}{0.68,1.00,0.18}
\definecolor{HotPink1}{rgb}{1.00,0.43,0.71}
\definecolor{HotPink2}{rgb}{0.93,0.42,0.65}
\definecolor{HotPink3}{rgb}{0.80,0.38,0.56}
\definecolor{HotPink4}{rgb}{0.55,0.23,0.38}
\definecolor{HotPink}{rgb}{1.00,0.41,0.71}
\definecolor{IndianRed1}{rgb}{1.00,0.42,0.42}
\definecolor{IndianRed2}{rgb}{0.93,0.39,0.39}
\definecolor{IndianRed3}{rgb}{0.80,0.33,0.33}
\definecolor{IndianRed4}{rgb}{0.55,0.23,0.23}
\definecolor{IndianRed}{rgb}{0.80,0.36,0.36}
\definecolor{LavenderBlush1}{rgb}{1.00,0.94,0.96}
\definecolor{LavenderBlush2}{rgb}{0.93,0.88,0.90}
\definecolor{LavenderBlush3}{rgb}{0.80,0.76,0.77}
\definecolor{LavenderBlush4}{rgb}{0.55,0.51,0.53}
\definecolor{LavenderBlush}{rgb}{1.00,0.94,0.96}
\definecolor{LawnGreen}{rgb}{0.49,0.99,0.00}
\definecolor{LemonChiffon1}{rgb}{1.00,0.98,0.80}
\definecolor{LemonChiffon2}{rgb}{0.93,0.91,0.75}
\definecolor{LemonChiffon3}{rgb}{0.80,0.79,0.65}
\definecolor{LemonChiffon4}{rgb}{0.55,0.54,0.44}
\definecolor{LemonChiffon}{rgb}{1.00,0.98,0.80}
\definecolor{LightBlue1}{rgb}{0.75,0.94,1.00}
\definecolor{LightBlue2}{rgb}{0.70,0.87,0.93}
\definecolor{LightBlue3}{rgb}{0.60,0.75,0.80}
\definecolor{LightBlue4}{rgb}{0.41,0.51,0.55}
\definecolor{LightBlue}{rgb}{0.68,0.85,0.90}
\definecolor{LightCoral}{rgb}{0.94,0.50,0.50}
\definecolor{LightCyan1}{rgb}{0.88,1.00,1.00}
\definecolor{LightCyan2}{rgb}{0.82,0.93,0.93}
\definecolor{LightCyan3}{rgb}{0.71,0.80,0.80}
\definecolor{LightCyan4}{rgb}{0.48,0.55,0.55}
\definecolor{LightCyan}{rgb}{0.88,1.00,1.00}
\definecolor{LightGoldenrod1}{rgb}{1.00,0.93,0.55}
\definecolor{LightGoldenrod2}{rgb}{0.93,0.86,0.51}
\definecolor{LightGoldenrod3}{rgb}{0.80,0.75,0.44}
\definecolor{LightGoldenrod4}{rgb}{0.55,0.51,0.30}
\definecolor{LightGoldenrodYellow}{rgb}{0.98,0.98,0.82}
\definecolor{LightGoldenrod}{rgb}{0.93,0.87,0.51}
\definecolor{LightGray}{rgb}{0.83,0.83,0.83}
\definecolor{LightGreen}{rgb}{0.56,0.93,0.56}
\definecolor{LightGrey}{rgb}{0.83,0.83,0.83}
\definecolor{LightPink1}{rgb}{1.00,0.68,0.73}
\definecolor{LightPink2}{rgb}{0.93,0.64,0.68}
\definecolor{LightPink3}{rgb}{0.80,0.55,0.58}
\definecolor{LightPink4}{rgb}{0.55,0.37,0.40}
\definecolor{LightPink}{rgb}{1.00,0.71,0.76}
\definecolor{LightSalmon1}{rgb}{1.00,0.63,0.48}
\definecolor{LightSalmon2}{rgb}{0.93,0.58,0.45}
\definecolor{LightSalmon3}{rgb}{0.80,0.51,0.38}
\definecolor{LightSalmon4}{rgb}{0.55,0.34,0.26}
\definecolor{LightSalmon}{rgb}{1.00,0.63,0.48}
\definecolor{LightSeaGreen}{rgb}{0.13,0.70,0.67}
\definecolor{LightSkyBlue1}{rgb}{0.69,0.89,1.00}
\definecolor{LightSkyBlue2}{rgb}{0.64,0.83,0.93}
\definecolor{LightSkyBlue3}{rgb}{0.55,0.71,0.80}
\definecolor{LightSkyBlue4}{rgb}{0.38,0.48,0.55}
\definecolor{LightSkyBlue}{rgb}{0.53,0.81,0.98}
\definecolor{LightSlateBlue}{rgb}{0.52,0.44,1.00}
\definecolor{LightSlateGray}{rgb}{0.47,0.53,0.60}
\definecolor{LightSlateGrey}{rgb}{0.47,0.53,0.60}
\definecolor{LightSteelBlue1}{rgb}{0.79,0.88,1.00}
\definecolor{LightSteelBlue2}{rgb}{0.74,0.82,0.93}
\definecolor{LightSteelBlue3}{rgb}{0.64,0.71,0.80}
\definecolor{LightSteelBlue4}{rgb}{0.43,0.48,0.55}
\definecolor{LightSteelBlue}{rgb}{0.69,0.77,0.87}
\definecolor{LightYellow1}{rgb}{1.00,1.00,0.88}
\definecolor{LightYellow2}{rgb}{0.93,0.93,0.82}
\definecolor{LightYellow3}{rgb}{0.80,0.80,0.71}
\definecolor{LightYellow4}{rgb}{0.55,0.55,0.48}
\definecolor{LightYellow}{rgb}{1.00,1.00,0.88}
\definecolor{LimeGreen}{rgb}{0.20,0.80,0.20}
\definecolor{MediumAquamarine}{rgb}{0.40,0.80,0.67}
\definecolor{MediumBlue}{rgb}{0.00,0.00,0.80}
\definecolor{MediumOrchid1}{rgb}{0.88,0.40,1.00}
\definecolor{MediumOrchid2}{rgb}{0.82,0.37,0.93}
\definecolor{MediumOrchid3}{rgb}{0.71,0.32,0.80}
\definecolor{MediumOrchid4}{rgb}{0.48,0.22,0.55}
\definecolor{MediumOrchid}{rgb}{0.73,0.33,0.83}
\definecolor{MediumPurple1}{rgb}{0.67,0.51,1.00}
\definecolor{MediumPurple2}{rgb}{0.62,0.47,0.93}
\definecolor{MediumPurple3}{rgb}{0.54,0.41,0.80}
\definecolor{MediumPurple4}{rgb}{0.36,0.28,0.55}
\definecolor{MediumPurple}{rgb}{0.58,0.44,0.86}
\definecolor{MediumSeaGreen}{rgb}{0.24,0.70,0.44}
\definecolor{MediumSlateBlue}{rgb}{0.48,0.41,0.93}
\definecolor{MediumSpringGreen}{rgb}{0.00,0.98,0.60}
\definecolor{MediumTurquoise}{rgb}{0.28,0.82,0.80}
\definecolor{MediumVioletRed}{rgb}{0.78,0.08,0.52}
\definecolor{MidnightBlue}{rgb}{0.10,0.10,0.44}
\definecolor{MintCream}{rgb}{0.96,1.00,0.98}
\definecolor{MistyRose1}{rgb}{1.00,0.89,0.88}
\definecolor{MistyRose2}{rgb}{0.93,0.84,0.82}
\definecolor{MistyRose3}{rgb}{0.80,0.72,0.71}
\definecolor{MistyRose4}{rgb}{0.55,0.49,0.48}
\definecolor{MistyRose}{rgb}{1.00,0.89,0.88}
\definecolor{NavajoWhite1}{rgb}{1.00,0.87,0.68}
\definecolor{NavajoWhite2}{rgb}{0.93,0.81,0.63}
\definecolor{NavajoWhite3}{rgb}{0.80,0.70,0.55}
\definecolor{NavajoWhite4}{rgb}{0.55,0.47,0.37}
\definecolor{NavajoWhite}{rgb}{1.00,0.87,0.68}
\definecolor{NavyBlue}{rgb}{0.00,0.00,0.50}
\definecolor{OldLace}{rgb}{0.99,0.96,0.90}
\definecolor{OliveDrab1}{rgb}{0.75,1.00,0.24}
\definecolor{OliveDrab2}{rgb}{0.70,0.93,0.23}
\definecolor{OliveDrab3}{rgb}{0.60,0.80,0.20}
\definecolor{OliveDrab4}{rgb}{0.41,0.55,0.13}
\definecolor{OliveDrab}{rgb}{0.42,0.56,0.14}
\definecolor{OrangeRed1}{rgb}{1.00,0.27,0.00}
\definecolor{OrangeRed2}{rgb}{0.93,0.25,0.00}
\definecolor{OrangeRed3}{rgb}{0.80,0.22,0.00}
\definecolor{OrangeRed4}{rgb}{0.55,0.15,0.00}
\definecolor{OrangeRed}{rgb}{1.00,0.27,0.00}
\definecolor{PaleGoldenrod}{rgb}{0.93,0.91,0.67}
\definecolor{PaleGreen1}{rgb}{0.60,1.00,0.60}
\definecolor{PaleGreen2}{rgb}{0.56,0.93,0.56}
\definecolor{PaleGreen3}{rgb}{0.49,0.80,0.49}
\definecolor{PaleGreen4}{rgb}{0.33,0.55,0.33}
\definecolor{PaleGreen}{rgb}{0.60,0.98,0.60}
\definecolor{PaleTurquoise1}{rgb}{0.73,1.00,1.00}
\definecolor{PaleTurquoise2}{rgb}{0.68,0.93,0.93}
\definecolor{PaleTurquoise3}{rgb}{0.59,0.80,0.80}
\definecolor{PaleTurquoise4}{rgb}{0.40,0.55,0.55}
\definecolor{PaleTurquoise}{rgb}{0.69,0.93,0.93}
\definecolor{PaleVioletRed1}{rgb}{1.00,0.51,0.67}
\definecolor{PaleVioletRed2}{rgb}{0.93,0.47,0.62}
\definecolor{PaleVioletRed3}{rgb}{0.80,0.41,0.54}
\definecolor{PaleVioletRed4}{rgb}{0.55,0.28,0.36}
\definecolor{PaleVioletRed}{rgb}{0.86,0.44,0.58}
\definecolor{PapayaWhip}{rgb}{1.00,0.94,0.84}
\definecolor{PeachPuff1}{rgb}{1.00,0.85,0.73}
\definecolor{PeachPuff2}{rgb}{0.93,0.80,0.68}
\definecolor{PeachPuff3}{rgb}{0.80,0.69,0.58}
\definecolor{PeachPuff4}{rgb}{0.55,0.47,0.40}
\definecolor{PeachPuff}{rgb}{1.00,0.85,0.73}
\definecolor{PowderBlue}{rgb}{0.69,0.88,0.90}
\definecolor{RosyBrown1}{rgb}{1.00,0.76,0.76}
\definecolor{RosyBrown2}{rgb}{0.93,0.71,0.71}
\definecolor{RosyBrown3}{rgb}{0.80,0.61,0.61}
\definecolor{RosyBrown4}{rgb}{0.55,0.41,0.41}
\definecolor{RosyBrown}{rgb}{0.74,0.56,0.56}
\definecolor{RoyalBlue1}{rgb}{0.28,0.46,1.00}
\definecolor{RoyalBlue2}{rgb}{0.26,0.43,0.93}
\definecolor{RoyalBlue3}{rgb}{0.23,0.37,0.80}
\definecolor{RoyalBlue4}{rgb}{0.15,0.25,0.55}
\definecolor{RoyalBlue}{rgb}{0.25,0.41,0.88}
\definecolor{SaddleBrown}{rgb}{0.55,0.27,0.07}
\definecolor{SandyBrown}{rgb}{0.96,0.64,0.38}
\definecolor{SeaGreen1}{rgb}{0.33,1.00,0.62}
\definecolor{SeaGreen2}{rgb}{0.31,0.93,0.58}
\definecolor{SeaGreen3}{rgb}{0.26,0.80,0.50}
\definecolor{SeaGreen4}{rgb}{0.18,0.55,0.34}
\definecolor{SeaGreen}{rgb}{0.18,0.55,0.34}
\definecolor{SkyBlue1}{rgb}{0.53,0.81,1.00}
\definecolor{SkyBlue2}{rgb}{0.49,0.75,0.93}
\definecolor{SkyBlue3}{rgb}{0.42,0.65,0.80}
\definecolor{SkyBlue4}{rgb}{0.29,0.44,0.55}
\definecolor{SkyBlue}{rgb}{0.53,0.81,0.92}
\definecolor{SlateBlue1}{rgb}{0.51,0.44,1.00}
\definecolor{SlateBlue2}{rgb}{0.48,0.40,0.93}
\definecolor{SlateBlue3}{rgb}{0.41,0.35,0.80}
\definecolor{SlateBlue4}{rgb}{0.28,0.24,0.55}
\definecolor{SlateBlue}{rgb}{0.42,0.35,0.80}
\definecolor{SlateGray1}{rgb}{0.78,0.89,1.00}
\definecolor{SlateGray2}{rgb}{0.73,0.83,0.93}
\definecolor{SlateGray3}{rgb}{0.62,0.71,0.80}
\definecolor{SlateGray4}{rgb}{0.42,0.48,0.55}
\definecolor{SlateGray}{rgb}{0.44,0.50,0.56}
\definecolor{SlateGrey}{rgb}{0.44,0.50,0.56}
\definecolor{SpringGreen1}{rgb}{0.00,1.00,0.50}
\definecolor{SpringGreen2}{rgb}{0.00,0.93,0.46}
\definecolor{SpringGreen3}{rgb}{0.00,0.80,0.40}
\definecolor{SpringGreen4}{rgb}{0.00,0.55,0.27}
\definecolor{SpringGreen}{rgb}{0.00,1.00,0.50}
\definecolor{SteelBlue1}{rgb}{0.39,0.72,1.00}
\definecolor{SteelBlue2}{rgb}{0.36,0.67,0.93}
\definecolor{SteelBlue3}{rgb}{0.31,0.58,0.80}
\definecolor{SteelBlue4}{rgb}{0.21,0.39,0.55}
\definecolor{SteelBlue}{rgb}{0.27,0.51,0.71}
\definecolor{VioletRed1}{rgb}{1.00,0.24,0.59}
\definecolor{VioletRed2}{rgb}{0.93,0.23,0.55}
\definecolor{VioletRed3}{rgb}{0.80,0.20,0.47}
\definecolor{VioletRed4}{rgb}{0.55,0.13,0.32}
\definecolor{VioletRed}{rgb}{0.82,0.13,0.56}
\definecolor{WhiteSmoke}{rgb}{0.96,0.96,0.96}
\definecolor{YellowGreen}{rgb}{0.60,0.80,0.20}
\definecolor{aliceblue}{rgb}{0.94,0.97,1.00}
\definecolor{antiquewhite}{rgb}{0.98,0.92,0.84}
\definecolor{aquamarine1}{rgb}{0.50,1.00,0.83}
\definecolor{aquamarine2}{rgb}{0.46,0.93,0.78}
\definecolor{aquamarine3}{rgb}{0.40,0.80,0.67}
\definecolor{aquamarine4}{rgb}{0.27,0.55,0.45}
\definecolor{aquamarine}{rgb}{0.50,1.00,0.83}
\definecolor{azure1}{rgb}{0.94,1.00,1.00}
\definecolor{azure2}{rgb}{0.88,0.93,0.93}
\definecolor{azure3}{rgb}{0.76,0.80,0.80}
\definecolor{azure4}{rgb}{0.51,0.55,0.55}
\definecolor{azure}{rgb}{0.94,1.00,1.00}
\definecolor{beige}{rgb}{0.96,0.96,0.86}
\definecolor{bisque1}{rgb}{1.00,0.89,0.77}
\definecolor{bisque2}{rgb}{0.93,0.84,0.72}
\definecolor{bisque3}{rgb}{0.80,0.72,0.62}
\definecolor{bisque4}{rgb}{0.55,0.49,0.42}
\definecolor{bisque}{rgb}{1.00,0.89,0.77}
\definecolor{black}{rgb}{0.00,0.00,0.00}
\definecolor{blanchedalmond}{rgb}{1.00,0.92,0.80}
\definecolor{blue1}{rgb}{0.00,0.00,1.00}
\definecolor{blue2}{rgb}{0.00,0.00,0.93}
\definecolor{blue3}{rgb}{0.00,0.00,0.80}
\definecolor{blue4}{rgb}{0.00,0.00,0.55}
\definecolor{blueviolet}{rgb}{0.54,0.17,0.89}
\definecolor{blue}{rgb}{0.00,0.00,1.00}
\definecolor{brown1}{rgb}{1.00,0.25,0.25}
\definecolor{brown2}{rgb}{0.93,0.23,0.23}
\definecolor{brown3}{rgb}{0.80,0.20,0.20}
\definecolor{brown4}{rgb}{0.55,0.14,0.14}
\definecolor{brown}{rgb}{0.65,0.16,0.16}
\definecolor{burlywood1}{rgb}{1.00,0.83,0.61}
\definecolor{burlywood2}{rgb}{0.93,0.77,0.57}
\definecolor{burlywood3}{rgb}{0.80,0.67,0.49}
\definecolor{burlywood4}{rgb}{0.55,0.45,0.33}
\definecolor{burlywood}{rgb}{0.87,0.72,0.53}
\definecolor{cadetblue}{rgb}{0.37,0.62,0.63}
\definecolor{chartreuse1}{rgb}{0.50,1.00,0.00}
\definecolor{chartreuse2}{rgb}{0.46,0.93,0.00}
\definecolor{chartreuse3}{rgb}{0.40,0.80,0.00}
\definecolor{chartreuse4}{rgb}{0.27,0.55,0.00}
\definecolor{chartreuse}{rgb}{0.50,1.00,0.00}
\definecolor{chocolate1}{rgb}{1.00,0.50,0.14}
\definecolor{chocolate2}{rgb}{0.93,0.46,0.13}
\definecolor{chocolate3}{rgb}{0.80,0.40,0.11}
\definecolor{chocolate4}{rgb}{0.55,0.27,0.07}
\definecolor{chocolate}{rgb}{0.82,0.41,0.12}
\definecolor{coral1}{rgb}{1.00,0.45,0.34}
\definecolor{coral2}{rgb}{0.93,0.42,0.31}
\definecolor{coral3}{rgb}{0.80,0.36,0.27}
\definecolor{coral4}{rgb}{0.55,0.24,0.18}
\definecolor{coral}{rgb}{1.00,0.50,0.31}
\definecolor{cornflowerblue}{rgb}{0.39,0.58,0.93}
\definecolor{cornsilk1}{rgb}{1.00,0.97,0.86}
\definecolor{cornsilk2}{rgb}{0.93,0.91,0.80}
\definecolor{cornsilk3}{rgb}{0.80,0.78,0.69}
\definecolor{cornsilk4}{rgb}{0.55,0.53,0.47}
\definecolor{cornsilk}{rgb}{1.00,0.97,0.86}
\definecolor{cyan1}{rgb}{0.00,1.00,1.00}
\definecolor{cyan2}{rgb}{0.00,0.93,0.93}
\definecolor{cyan3}{rgb}{0.00,0.80,0.80}
\definecolor{cyan4}{rgb}{0.00,0.55,0.55}
\definecolor{cyan}{rgb}{0.00,1.00,1.00}
\definecolor{darkblue}{rgb}{0.00,0.00,0.55}
\definecolor{darkcyan}{rgb}{0.00,0.55,0.55}
\definecolor{darkgoldenrod}{rgb}{0.72,0.53,0.04}
\definecolor{darkgray}{rgb}{0.66,0.66,0.66}
\definecolor{darkgreen}{rgb}{0.00,0.39,0.00}
\definecolor{darkgrey}{rgb}{0.66,0.66,0.66}
\definecolor{darkkhaki}{rgb}{0.74,0.72,0.42}
\definecolor{darkmagenta}{rgb}{0.55,0.00,0.55}
\definecolor{darkolive}{rgb}{0.33,0.42,0.18}
\definecolor{darkorange}{rgb}{1.00,0.55,0.00}
\definecolor{darkorchid}{rgb}{0.60,0.20,0.80}
\definecolor{darkred}{rgb}{0.55,0.00,0.00}
\definecolor{darksalmon}{rgb}{0.91,0.59,0.48}
\definecolor{darksea}{rgb}{0.56,0.74,0.56}
\definecolor{darkslate}{rgb}{0.18,0.31,0.31}
\definecolor{darkslate}{rgb}{0.18,0.31,0.31}
\definecolor{darkslate}{rgb}{0.28,0.24,0.55}
\definecolor{darkturquoise}{rgb}{0.00,0.81,0.82}
\definecolor{darkviolet}{rgb}{0.58,0.00,0.83}
\definecolor{deeppink}{rgb}{1.00,0.08,0.58}
\definecolor{deepsky}{rgb}{0.00,0.75,1.00}
\definecolor{dimgray}{rgb}{0.41,0.41,0.41}
\definecolor{dimgrey}{rgb}{0.41,0.41,0.41}
\definecolor{dodgerblue}{rgb}{0.12,0.56,1.00}
\definecolor{firebrick1}{rgb}{1.00,0.19,0.19}
\definecolor{firebrick2}{rgb}{0.93,0.17,0.17}
\definecolor{firebrick3}{rgb}{0.80,0.15,0.15}
\definecolor{firebrick4}{rgb}{0.55,0.10,0.10}
\definecolor{firebrick}{rgb}{0.70,0.13,0.13}
\definecolor{floralwhite}{rgb}{1.00,0.98,0.94}
\definecolor{forestgreen}{rgb}{0.13,0.55,0.13}
\definecolor{gainsboro}{rgb}{0.86,0.86,0.86}
\definecolor{ghostwhite}{rgb}{0.97,0.97,1.00}
\definecolor{gold1}{rgb}{1.00,0.84,0.00}
\definecolor{gold2}{rgb}{0.93,0.79,0.00}
\definecolor{gold3}{rgb}{0.80,0.68,0.00}
\definecolor{gold4}{rgb}{0.55,0.46,0.00}
\definecolor{goldenrod1}{rgb}{1.00,0.76,0.15}
\definecolor{goldenrod2}{rgb}{0.93,0.71,0.13}
\definecolor{goldenrod3}{rgb}{0.80,0.61,0.11}
\definecolor{goldenrod4}{rgb}{0.55,0.41,0.08}
\definecolor{goldenrod}{rgb}{0.85,0.65,0.13}
\definecolor{gold}{rgb}{1.00,0.84,0.00}
\definecolor{gray0}{rgb}{0.00,0.00,0.00}
\definecolor{gray100}{rgb}{1.00,1.00,1.00}
\definecolor{gray10}{rgb}{0.10,0.10,0.10}
\definecolor{gray11}{rgb}{0.11,0.11,0.11}
\definecolor{gray12}{rgb}{0.12,0.12,0.12}
\definecolor{gray13}{rgb}{0.13,0.13,0.13}
\definecolor{gray14}{rgb}{0.14,0.14,0.14}
\definecolor{gray15}{rgb}{0.15,0.15,0.15}
\definecolor{gray16}{rgb}{0.16,0.16,0.16}
\definecolor{gray17}{rgb}{0.17,0.17,0.17}
\definecolor{gray18}{rgb}{0.18,0.18,0.18}
\definecolor{gray19}{rgb}{0.19,0.19,0.19}
\definecolor{gray1}{rgb}{0.01,0.01,0.01}
\definecolor{gray20}{rgb}{0.20,0.20,0.20}
\definecolor{gray21}{rgb}{0.21,0.21,0.21}
\definecolor{gray22}{rgb}{0.22,0.22,0.22}
\definecolor{gray23}{rgb}{0.23,0.23,0.23}
\definecolor{gray24}{rgb}{0.24,0.24,0.24}
\definecolor{gray25}{rgb}{0.25,0.25,0.25}
\definecolor{gray26}{rgb}{0.26,0.26,0.26}
\definecolor{gray27}{rgb}{0.27,0.27,0.27}
\definecolor{gray28}{rgb}{0.28,0.28,0.28}
\definecolor{gray29}{rgb}{0.29,0.29,0.29}
\definecolor{gray2}{rgb}{0.02,0.02,0.02}
\definecolor{gray30}{rgb}{0.30,0.30,0.30}
\definecolor{gray31}{rgb}{0.31,0.31,0.31}
\definecolor{gray32}{rgb}{0.32,0.32,0.32}
\definecolor{gray33}{rgb}{0.33,0.33,0.33}
\definecolor{gray34}{rgb}{0.34,0.34,0.34}
\definecolor{gray35}{rgb}{0.35,0.35,0.35}
\definecolor{gray36}{rgb}{0.36,0.36,0.36}
\definecolor{gray37}{rgb}{0.37,0.37,0.37}
\definecolor{gray38}{rgb}{0.38,0.38,0.38}
\definecolor{gray39}{rgb}{0.39,0.39,0.39}
\definecolor{gray3}{rgb}{0.03,0.03,0.03}
\definecolor{gray40}{rgb}{0.40,0.40,0.40}
\definecolor{gray41}{rgb}{0.41,0.41,0.41}
\definecolor{gray42}{rgb}{0.42,0.42,0.42}
\definecolor{gray43}{rgb}{0.43,0.43,0.43}
\definecolor{gray44}{rgb}{0.44,0.44,0.44}
\definecolor{gray45}{rgb}{0.45,0.45,0.45}
\definecolor{gray46}{rgb}{0.46,0.46,0.46}
\definecolor{gray47}{rgb}{0.47,0.47,0.47}
\definecolor{gray48}{rgb}{0.48,0.48,0.48}
\definecolor{gray49}{rgb}{0.49,0.49,0.49}
\definecolor{gray4}{rgb}{0.04,0.04,0.04}
\definecolor{gray50}{rgb}{0.50,0.50,0.50}
\definecolor{gray51}{rgb}{0.51,0.51,0.51}
\definecolor{gray52}{rgb}{0.52,0.52,0.52}
\definecolor{gray53}{rgb}{0.53,0.53,0.53}
\definecolor{gray54}{rgb}{0.54,0.54,0.54}
\definecolor{gray55}{rgb}{0.55,0.55,0.55}
\definecolor{gray56}{rgb}{0.56,0.56,0.56}
\definecolor{gray57}{rgb}{0.57,0.57,0.57}
\definecolor{gray58}{rgb}{0.58,0.58,0.58}
\definecolor{gray59}{rgb}{0.59,0.59,0.59}
\definecolor{gray5}{rgb}{0.05,0.05,0.05}
\definecolor{gray60}{rgb}{0.60,0.60,0.60}
\definecolor{gray61}{rgb}{0.61,0.61,0.61}
\definecolor{gray62}{rgb}{0.62,0.62,0.62}
\definecolor{gray63}{rgb}{0.63,0.63,0.63}
\definecolor{gray64}{rgb}{0.64,0.64,0.64}
\definecolor{gray65}{rgb}{0.65,0.65,0.65}
\definecolor{gray66}{rgb}{0.66,0.66,0.66}
\definecolor{gray67}{rgb}{0.67,0.67,0.67}
\definecolor{gray68}{rgb}{0.68,0.68,0.68}
\definecolor{gray69}{rgb}{0.69,0.69,0.69}
\definecolor{gray6}{rgb}{0.06,0.06,0.06}
\definecolor{gray70}{rgb}{0.70,0.70,0.70}
\definecolor{gray71}{rgb}{0.71,0.71,0.71}
\definecolor{gray72}{rgb}{0.72,0.72,0.72}
\definecolor{gray73}{rgb}{0.73,0.73,0.73}
\definecolor{gray74}{rgb}{0.74,0.74,0.74}
\definecolor{gray75}{rgb}{0.75,0.75,0.75}
\definecolor{gray76}{rgb}{0.76,0.76,0.76}
\definecolor{gray77}{rgb}{0.77,0.77,0.77}
\definecolor{gray78}{rgb}{0.78,0.78,0.78}
\definecolor{gray79}{rgb}{0.79,0.79,0.79}
\definecolor{gray7}{rgb}{0.07,0.07,0.07}
\definecolor{gray80}{rgb}{0.80,0.80,0.80}
\definecolor{gray81}{rgb}{0.81,0.81,0.81}
\definecolor{gray82}{rgb}{0.82,0.82,0.82}
\definecolor{gray83}{rgb}{0.83,0.83,0.83}
\definecolor{gray84}{rgb}{0.84,0.84,0.84}
\definecolor{gray85}{rgb}{0.85,0.85,0.85}
\definecolor{gray86}{rgb}{0.86,0.86,0.86}
\definecolor{gray87}{rgb}{0.87,0.87,0.87}
\definecolor{gray88}{rgb}{0.88,0.88,0.88}
\definecolor{gray89}{rgb}{0.89,0.89,0.89}
\definecolor{gray8}{rgb}{0.08,0.08,0.08}
\definecolor{gray90}{rgb}{0.90,0.90,0.90}
\definecolor{gray91}{rgb}{0.91,0.91,0.91}
\definecolor{gray92}{rgb}{0.92,0.92,0.92}
\definecolor{gray93}{rgb}{0.93,0.93,0.93}
\definecolor{gray94}{rgb}{0.94,0.94,0.94}
\definecolor{gray95}{rgb}{0.95,0.95,0.95}
\definecolor{gray96}{rgb}{0.96,0.96,0.96}
\definecolor{gray97}{rgb}{0.97,0.97,0.97}
\definecolor{gray98}{rgb}{0.98,0.98,0.98}
\definecolor{gray99}{rgb}{0.99,0.99,0.99}
\definecolor{gray9}{rgb}{0.09,0.09,0.09}
\definecolor{gray}{rgb}{0.75,0.75,0.75}
\definecolor{green1}{rgb}{0.00,1.00,0.00}
\definecolor{green2}{rgb}{0.00,0.93,0.00}
\definecolor{green3}{rgb}{0.00,0.80,0.00}
\definecolor{green4}{rgb}{0.00,0.55,0.00}
\definecolor{greenyellow}{rgb}{0.68,1.00,0.18}
\definecolor{green}{rgb}{0.00,1.00,0.00}
\definecolor{grey0}{rgb}{0.00,0.00,0.00}
\definecolor{grey100}{rgb}{1.00,1.00,1.00}
\definecolor{grey10}{rgb}{0.10,0.10,0.10}
\definecolor{grey11}{rgb}{0.11,0.11,0.11}
\definecolor{grey12}{rgb}{0.12,0.12,0.12}
\definecolor{grey13}{rgb}{0.13,0.13,0.13}
\definecolor{grey14}{rgb}{0.14,0.14,0.14}
\definecolor{grey15}{rgb}{0.15,0.15,0.15}
\definecolor{grey16}{rgb}{0.16,0.16,0.16}
\definecolor{grey17}{rgb}{0.17,0.17,0.17}
\definecolor{grey18}{rgb}{0.18,0.18,0.18}
\definecolor{grey19}{rgb}{0.19,0.19,0.19}
\definecolor{grey1}{rgb}{0.01,0.01,0.01}
\definecolor{grey20}{rgb}{0.20,0.20,0.20}
\definecolor{grey21}{rgb}{0.21,0.21,0.21}
\definecolor{grey22}{rgb}{0.22,0.22,0.22}
\definecolor{grey23}{rgb}{0.23,0.23,0.23}
\definecolor{grey24}{rgb}{0.24,0.24,0.24}
\definecolor{grey25}{rgb}{0.25,0.25,0.25}
\definecolor{grey26}{rgb}{0.26,0.26,0.26}
\definecolor{grey27}{rgb}{0.27,0.27,0.27}
\definecolor{grey28}{rgb}{0.28,0.28,0.28}
\definecolor{grey29}{rgb}{0.29,0.29,0.29}
\definecolor{grey2}{rgb}{0.02,0.02,0.02}
\definecolor{grey30}{rgb}{0.30,0.30,0.30}
\definecolor{grey31}{rgb}{0.31,0.31,0.31}
\definecolor{grey32}{rgb}{0.32,0.32,0.32}
\definecolor{grey33}{rgb}{0.33,0.33,0.33}
\definecolor{grey34}{rgb}{0.34,0.34,0.34}
\definecolor{grey35}{rgb}{0.35,0.35,0.35}
\definecolor{grey36}{rgb}{0.36,0.36,0.36}
\definecolor{grey37}{rgb}{0.37,0.37,0.37}
\definecolor{grey38}{rgb}{0.38,0.38,0.38}
\definecolor{grey39}{rgb}{0.39,0.39,0.39}
\definecolor{grey3}{rgb}{0.03,0.03,0.03}
\definecolor{grey40}{rgb}{0.40,0.40,0.40}
\definecolor{grey41}{rgb}{0.41,0.41,0.41}
\definecolor{grey42}{rgb}{0.42,0.42,0.42}
\definecolor{grey43}{rgb}{0.43,0.43,0.43}
\definecolor{grey44}{rgb}{0.44,0.44,0.44}
\definecolor{grey45}{rgb}{0.45,0.45,0.45}
\definecolor{grey46}{rgb}{0.46,0.46,0.46}
\definecolor{grey47}{rgb}{0.47,0.47,0.47}
\definecolor{grey48}{rgb}{0.48,0.48,0.48}
\definecolor{grey49}{rgb}{0.49,0.49,0.49}
\definecolor{grey4}{rgb}{0.04,0.04,0.04}
\definecolor{grey50}{rgb}{0.50,0.50,0.50}
\definecolor{grey51}{rgb}{0.51,0.51,0.51}
\definecolor{grey52}{rgb}{0.52,0.52,0.52}
\definecolor{grey53}{rgb}{0.53,0.53,0.53}
\definecolor{grey54}{rgb}{0.54,0.54,0.54}
\definecolor{grey55}{rgb}{0.55,0.55,0.55}
\definecolor{grey56}{rgb}{0.56,0.56,0.56}
\definecolor{grey57}{rgb}{0.57,0.57,0.57}
\definecolor{grey58}{rgb}{0.58,0.58,0.58}
\definecolor{grey59}{rgb}{0.59,0.59,0.59}
\definecolor{grey5}{rgb}{0.05,0.05,0.05}
\definecolor{grey60}{rgb}{0.60,0.60,0.60}
\definecolor{grey61}{rgb}{0.61,0.61,0.61}
\definecolor{grey62}{rgb}{0.62,0.62,0.62}
\definecolor{grey63}{rgb}{0.63,0.63,0.63}
\definecolor{grey64}{rgb}{0.64,0.64,0.64}
\definecolor{grey65}{rgb}{0.65,0.65,0.65}
\definecolor{grey66}{rgb}{0.66,0.66,0.66}
\definecolor{grey67}{rgb}{0.67,0.67,0.67}
\definecolor{grey68}{rgb}{0.68,0.68,0.68}
\definecolor{grey69}{rgb}{0.69,0.69,0.69}
\definecolor{grey6}{rgb}{0.06,0.06,0.06}
\definecolor{grey70}{rgb}{0.70,0.70,0.70}
\definecolor{grey71}{rgb}{0.71,0.71,0.71}
\definecolor{grey72}{rgb}{0.72,0.72,0.72}
\definecolor{grey73}{rgb}{0.73,0.73,0.73}
\definecolor{grey74}{rgb}{0.74,0.74,0.74}
\definecolor{grey75}{rgb}{0.75,0.75,0.75}
\definecolor{grey76}{rgb}{0.76,0.76,0.76}
\definecolor{grey77}{rgb}{0.77,0.77,0.77}
\definecolor{grey78}{rgb}{0.78,0.78,0.78}
\definecolor{grey79}{rgb}{0.79,0.79,0.79}
\definecolor{grey7}{rgb}{0.07,0.07,0.07}
\definecolor{grey80}{rgb}{0.80,0.80,0.80}
\definecolor{grey81}{rgb}{0.81,0.81,0.81}
\definecolor{grey82}{rgb}{0.82,0.82,0.82}
\definecolor{grey83}{rgb}{0.83,0.83,0.83}
\definecolor{grey84}{rgb}{0.84,0.84,0.84}
\definecolor{grey85}{rgb}{0.85,0.85,0.85}
\definecolor{grey86}{rgb}{0.86,0.86,0.86}
\definecolor{grey87}{rgb}{0.87,0.87,0.87}
\definecolor{grey88}{rgb}{0.88,0.88,0.88}
\definecolor{grey89}{rgb}{0.89,0.89,0.89}
\definecolor{grey8}{rgb}{0.08,0.08,0.08}
\definecolor{grey90}{rgb}{0.90,0.90,0.90}
\definecolor{grey91}{rgb}{0.91,0.91,0.91}
\definecolor{grey92}{rgb}{0.92,0.92,0.92}
\definecolor{grey93}{rgb}{0.93,0.93,0.93}
\definecolor{grey94}{rgb}{0.94,0.94,0.94}
\definecolor{grey95}{rgb}{0.95,0.95,0.95}
\definecolor{grey96}{rgb}{0.96,0.96,0.96}
\definecolor{grey97}{rgb}{0.97,0.97,0.97}
\definecolor{grey98}{rgb}{0.98,0.98,0.98}
\definecolor{grey99}{rgb}{0.99,0.99,0.99}
\definecolor{grey9}{rgb}{0.09,0.09,0.09}
\definecolor{grey}{rgb}{0.75,0.75,0.75}
\definecolor{honeydew1}{rgb}{0.94,1.00,0.94}
\definecolor{honeydew2}{rgb}{0.88,0.93,0.88}
\definecolor{honeydew3}{rgb}{0.76,0.80,0.76}
\definecolor{honeydew4}{rgb}{0.51,0.55,0.51}
\definecolor{honeydew}{rgb}{0.94,1.00,0.94}
\definecolor{hotpink}{rgb}{1.00,0.41,0.71}
\definecolor{indianred}{rgb}{0.80,0.36,0.36}
\definecolor{ivory1}{rgb}{1.00,1.00,0.94}
\definecolor{ivory2}{rgb}{0.93,0.93,0.88}
\definecolor{ivory3}{rgb}{0.80,0.80,0.76}
\definecolor{ivory4}{rgb}{0.55,0.55,0.51}
\definecolor{ivory}{rgb}{1.00,1.00,0.94}
\definecolor{khaki1}{rgb}{1.00,0.96,0.56}
\definecolor{khaki2}{rgb}{0.93,0.90,0.52}
\definecolor{khaki3}{rgb}{0.80,0.78,0.45}
\definecolor{khaki4}{rgb}{0.55,0.53,0.31}
\definecolor{khaki}{rgb}{0.94,0.90,0.55}
\definecolor{lavenderblush}{rgb}{1.00,0.94,0.96}
\definecolor{lavender}{rgb}{0.90,0.90,0.98}
\definecolor{lawngreen}{rgb}{0.49,0.99,0.00}
\definecolor{lemonchiffon}{rgb}{1.00,0.98,0.80}
\definecolor{lightblue}{rgb}{0.68,0.85,0.90}
\definecolor{lightcoral}{rgb}{0.94,0.50,0.50}
\definecolor{lightcyan}{rgb}{0.88,1.00,1.00}
\definecolor{lightgoldenrod}{rgb}{0.93,0.87,0.51}
\definecolor{lightgoldenrod}{rgb}{0.98,0.98,0.82}
\definecolor{lightgray}{rgb}{0.83,0.83,0.83}
\definecolor{lightgreen}{rgb}{0.56,0.93,0.56}
\definecolor{lightgrey}{rgb}{0.83,0.83,0.83}
\definecolor{lightpink}{rgb}{1.00,0.71,0.76}
\definecolor{lightsalmon}{rgb}{1.00,0.63,0.48}
\definecolor{lightsea}{rgb}{0.13,0.70,0.67}
\definecolor{lightsky}{rgb}{0.53,0.81,0.98}
\definecolor{lightslate}{rgb}{0.47,0.53,0.60}
\definecolor{lightslate}{rgb}{0.47,0.53,0.60}
\definecolor{lightslate}{rgb}{0.52,0.44,1.00}
\definecolor{lightsteel}{rgb}{0.69,0.77,0.87}
\definecolor{lightyellow}{rgb}{1.00,1.00,0.88}
\definecolor{limegreen}{rgb}{0.20,0.80,0.20}
\definecolor{linen}{rgb}{0.98,0.94,0.90}
\definecolor{magenta1}{rgb}{1.00,0.00,1.00}
\definecolor{magenta2}{rgb}{0.93,0.00,0.93}
\definecolor{magenta3}{rgb}{0.80,0.00,0.80}
\definecolor{magenta4}{rgb}{0.55,0.00,0.55}
\definecolor{magenta}{rgb}{1.00,0.00,1.00}
\definecolor{maroon1}{rgb}{1.00,0.20,0.70}
\definecolor{maroon2}{rgb}{0.93,0.19,0.65}
\definecolor{maroon3}{rgb}{0.80,0.16,0.56}
\definecolor{maroon4}{rgb}{0.55,0.11,0.38}
\definecolor{maroon}{rgb}{0.69,0.19,0.38}
\definecolor{mediumaquamarine}{rgb}{0.40,0.80,0.67}
\definecolor{mediumblue}{rgb}{0.00,0.00,0.80}
\definecolor{mediumorchid}{rgb}{0.73,0.33,0.83}
\definecolor{mediumpurple}{rgb}{0.58,0.44,0.86}
\definecolor{mediumsea}{rgb}{0.24,0.70,0.44}
\definecolor{mediumslate}{rgb}{0.48,0.41,0.93}
\definecolor{mediumspring}{rgb}{0.00,0.98,0.60}
\definecolor{mediumturquoise}{rgb}{0.28,0.82,0.80}
\definecolor{mediumviolet}{rgb}{0.78,0.08,0.52}
\definecolor{midnightblue}{rgb}{0.10,0.10,0.44}
\definecolor{mintcream}{rgb}{0.96,1.00,0.98}
\definecolor{mistyrose}{rgb}{1.00,0.89,0.88}
\definecolor{moccasin}{rgb}{1.00,0.89,0.71}
\definecolor{navajowhite}{rgb}{1.00,0.87,0.68}
\definecolor{navyblue}{rgb}{0.00,0.00,0.50}
\definecolor{navy}{rgb}{0.00,0.00,0.50}
\definecolor{oldlace}{rgb}{0.99,0.96,0.90}
\definecolor{olivedrab}{rgb}{0.42,0.56,0.14}
\definecolor{orange1}{rgb}{1.00,0.65,0.00}
\definecolor{orange2}{rgb}{0.93,0.60,0.00}
\definecolor{orange3}{rgb}{0.80,0.52,0.00}
\definecolor{orange4}{rgb}{0.55,0.35,0.00}
\definecolor{orangered}{rgb}{1.00,0.27,0.00}
\definecolor{orange}{rgb}{1.00,0.65,0.00}
\definecolor{orchid1}{rgb}{1.00,0.51,0.98}
\definecolor{orchid2}{rgb}{0.93,0.48,0.91}
\definecolor{orchid3}{rgb}{0.80,0.41,0.79}
\definecolor{orchid4}{rgb}{0.55,0.28,0.54}
\definecolor{orchid}{rgb}{0.85,0.44,0.84}
\definecolor{palegoldenrod}{rgb}{0.93,0.91,0.67}
\definecolor{palegreen}{rgb}{0.60,0.98,0.60}
\definecolor{paleturquoise}{rgb}{0.69,0.93,0.93}
\definecolor{paleviolet}{rgb}{0.86,0.44,0.58}
\definecolor{papayawhip}{rgb}{1.00,0.94,0.84}
\definecolor{peachpuff}{rgb}{1.00,0.85,0.73}
\definecolor{peru}{rgb}{0.80,0.52,0.25}
\definecolor{pink1}{rgb}{1.00,0.71,0.77}
\definecolor{pink2}{rgb}{0.93,0.66,0.72}
\definecolor{pink3}{rgb}{0.80,0.57,0.62}
\definecolor{pink4}{rgb}{0.55,0.39,0.42}
\definecolor{pink}{rgb}{1.00,0.75,0.80}
\definecolor{plum1}{rgb}{1.00,0.73,1.00}
\definecolor{plum2}{rgb}{0.93,0.68,0.93}
\definecolor{plum3}{rgb}{0.80,0.59,0.80}
\definecolor{plum4}{rgb}{0.55,0.40,0.55}
\definecolor{plum}{rgb}{0.87,0.63,0.87}
\definecolor{powderblue}{rgb}{0.69,0.88,0.90}
\definecolor{purple1}{rgb}{0.61,0.19,1.00}
\definecolor{purple2}{rgb}{0.57,0.17,0.93}
\definecolor{purple3}{rgb}{0.49,0.15,0.80}
\definecolor{purple4}{rgb}{0.33,0.10,0.55}
\definecolor{purple}{rgb}{0.63,0.13,0.94}
\definecolor{red1}{rgb}{1.00,0.00,0.00}
\definecolor{red2}{rgb}{0.93,0.00,0.00}
\definecolor{red3}{rgb}{0.80,0.00,0.00}
\definecolor{red4}{rgb}{0.55,0.00,0.00}
\definecolor{red}{rgb}{1.00,0.00,0.00}
\definecolor{rosybrown}{rgb}{0.74,0.56,0.56}
\definecolor{royalblue}{rgb}{0.25,0.41,0.88}
\definecolor{saddlebrown}{rgb}{0.55,0.27,0.07}
\definecolor{salmon1}{rgb}{1.00,0.55,0.41}
\definecolor{salmon2}{rgb}{0.93,0.51,0.38}
\definecolor{salmon3}{rgb}{0.80,0.44,0.33}
\definecolor{salmon4}{rgb}{0.55,0.30,0.22}
\definecolor{salmon}{rgb}{0.98,0.50,0.45}
\definecolor{sandybrown}{rgb}{0.96,0.64,0.38}
\definecolor{seagreen}{rgb}{0.18,0.55,0.34}
\definecolor{seashell1}{rgb}{1.00,0.96,0.93}
\definecolor{seashell2}{rgb}{0.93,0.90,0.87}
\definecolor{seashell3}{rgb}{0.80,0.77,0.75}
\definecolor{seashell4}{rgb}{0.55,0.53,0.51}
\definecolor{seashell}{rgb}{1.00,0.96,0.93}
\definecolor{sienna1}{rgb}{1.00,0.51,0.28}
\definecolor{sienna2}{rgb}{0.93,0.47,0.26}
\definecolor{sienna3}{rgb}{0.80,0.41,0.22}
\definecolor{sienna4}{rgb}{0.55,0.28,0.15}
\definecolor{sienna}{rgb}{0.63,0.32,0.18}
\definecolor{skyblue}{rgb}{0.53,0.81,0.92}
\definecolor{slateblue}{rgb}{0.42,0.35,0.80}
\definecolor{slategray}{rgb}{0.44,0.50,0.56}
\definecolor{slategrey}{rgb}{0.44,0.50,0.56}
\definecolor{snow1}{rgb}{1.00,0.98,0.98}
\definecolor{snow2}{rgb}{0.93,0.91,0.91}
\definecolor{snow3}{rgb}{0.80,0.79,0.79}
\definecolor{snow4}{rgb}{0.55,0.54,0.54}
\definecolor{snow}{rgb}{1.00,0.98,0.98}
\definecolor{springgreen}{rgb}{0.00,1.00,0.50}
\definecolor{steelblue}{rgb}{0.27,0.51,0.71}
\definecolor{tan1}{rgb}{1.00,0.65,0.31}
\definecolor{tan2}{rgb}{0.93,0.60,0.29}
\definecolor{tan3}{rgb}{0.80,0.52,0.25}
\definecolor{tan4}{rgb}{0.55,0.35,0.17}
\definecolor{tan}{rgb}{0.82,0.71,0.55}
\definecolor{thistle1}{rgb}{1.00,0.88,1.00}
\definecolor{thistle2}{rgb}{0.93,0.82,0.93}
\definecolor{thistle3}{rgb}{0.80,0.71,0.80}
\definecolor{thistle4}{rgb}{0.55,0.48,0.55}
\definecolor{thistle}{rgb}{0.85,0.75,0.85}
\definecolor{tomato1}{rgb}{1.00,0.39,0.28}
\definecolor{tomato2}{rgb}{0.93,0.36,0.26}
\definecolor{tomato3}{rgb}{0.80,0.31,0.22}
\definecolor{tomato4}{rgb}{0.55,0.21,0.15}
\definecolor{tomato}{rgb}{1.00,0.39,0.28}
\definecolor{turquoise1}{rgb}{0.00,0.96,1.00}
\definecolor{turquoise2}{rgb}{0.00,0.90,0.93}
\definecolor{turquoise3}{rgb}{0.00,0.77,0.80}
\definecolor{turquoise4}{rgb}{0.00,0.53,0.55}
\definecolor{turquoise}{rgb}{0.25,0.88,0.82}
\definecolor{violetred}{rgb}{0.82,0.13,0.56}
\definecolor{violet}{rgb}{0.93,0.51,0.93}
\definecolor{wheat1}{rgb}{1.00,0.91,0.73}
\definecolor{wheat2}{rgb}{0.93,0.85,0.68}
\definecolor{wheat3}{rgb}{0.80,0.73,0.59}
\definecolor{wheat4}{rgb}{0.55,0.49,0.40}
\definecolor{wheat}{rgb}{0.96,0.87,0.70}
\definecolor{whitesmoke}{rgb}{0.96,0.96,0.96}
\definecolor{white}{rgb}{1.00,1.00,1.00}
\definecolor{yellow1}{rgb}{1.00,1.00,0.00}
\definecolor{yellow2}{rgb}{0.93,0.93,0.00}
\definecolor{yellow3}{rgb}{0.80,0.80,0.00}
\definecolor{yellow4}{rgb}{0.55,0.55,0.00}
\definecolor{yellowgreen}{rgb}{0.60,0.80,0.20}
\definecolor{yellow}{rgb}{1.00,1.00,0.00}

\usepackage[usenames,dvipsnames]{xcolor}

\RequirePackage[colorlinks]{hyperref}
\RequirePackage{hypernat}

\usepackage[psamsfonts]{amsfonts}
 
\usepackage{pdflscape}

\usepackage{dsfont}

\newtheorem{Lem}{Lemma}[section]

\newtheorem{Theor}{Theorem}[section]

\usepackage{hyperref}
\hypersetup{pdfpagemode=FullScreen,  
colorlinks=true,
citecolor=Bittersweet,
linkcolor=Bittersweet,
urlcolor=Bittersweet 
}

\newcommand{\R}{\mathbb{R}}
\newcommand{\cqfd}{\hfill $\square$}
\newcommand{\VS}{V}

\begin{document}

\begin{frontmatter}
\title{Halfspace Depths for Scatter,\\ Concentration and Shape Matrices}
\runtitle{Scatter, Concentration and Shape Halfspace Depths}
\vspace{1mm}

\begin{aug}
\author{\fnms{Davy} \snm{Paindaveine}\thanksref{t1}\ead[label=e1]{dpaindav@ulb.ac.be}
\ead[label=u1,url]{http://homepages.ulb.ac.be/\textasciitilde dpaindav}}
\and
\author{\fnms{ Germain} \snm{Van Bever}\thanksref{t2}
\ead[label=e2]{gvbever@ulb.ac.be}
\ead[label=u2,url]{https://sites.google.com/site/germainvanbever}}

\thankstext{t1}{Research is supported by\vspace{-1mm} the IAP research network grant 
\vspace{.8mm}
\mbox{nr.} P7/06 of the Belgian government (Belgian Science Policy), the Cr\'{e}dit de Recherche  J.0113.16 of the FNRS (Fonds National pour la Recherche Scientifique), Communaut\'{e} 
\vspace{1mm}
Fran\c{c}aise de Belgique, and a grant from the National Bank of Belgium.}
\thankstext{t2}{Research is supported by\vspace{-1mm} the FC84444 grant 
\vspace{.8mm}
of the FNRS (Fonds National pour la Recherche Scientifique), Communaut\'{e} 
\vspace{1mm}
Fran\c{c}aise de Belgique.}
\runauthor{D. Paindaveine and G. Van Bever}
\vspace{2mm}

\affiliation{Universit\'{e} libre de Bruxelles}

\address{Universit\'{e} libre de Bruxelles\\
ECARES and D\'{e}partement de Math\'{e}matique\\
Avenue Roosevelt, 50, ECARES, CP114/04\\
B-1050, Brussels\\
Belgium\\
\printead{e1}\\
\printead{u1}}

\address{Universit\'{e} libre de Bruxelles\\
ECARES, 
Avenue Roosevelt, 50, CP114/04\\
B-1050, Brussels\\
Belgium\\
E-mail:\ \printead*{e2}\\
\printead{u2}}

\end{aug}

\vspace{4mm}
\begin{abstract}
We propose halfspace depth concepts for scatter, concentration and shape matrices. For scatter matrices, our concept is similar to those from \cite{Chenetal16} and  \cite{Zha2002}. Rather than focusing, as in these earlier works, on deepest scatter matrices, we thoroughly investigate the properties of the proposed depth and of the corresponding depth regions. We do so under minimal assumptions and, in particular, we do not restrict to elliptical distributions nor to absolutely continuous distributions. Interestingly, fully understanding scatter halfspace depth requires considering different geometries/topologies on the space of scatter matrices. We also discuss, in the spirit of \cite{ZuoSer2000A}, the structural properties a scatter depth should satisfy, and investigate whether or not these are met by scatter halfspace depth. Companion concepts of depth for concentration matrices and shape matrices are also proposed and studied. We show the practical relevance of the depth concepts considered in a real-data example from finance. 
\vspace{-1mm}
\end{abstract}



\end{frontmatter}

\section{Introduction}
\label{Secintro}

Statistical depth measures the centrality of a given location in~$\mathbb{R}^k$ with respect to a sample of $k$-variate observations, or, more generally, with respect to a probability measure~$P$ over~$\R^k$. The most famous depths include the halfspace depth (\citealp{Tuk1975}), the simplicial depth (\citealp{Liu1990}), the spatial depth (\citealp{VarZha2000}) and the projection depth (\citealp{Zuo2003}). In the last decade, depth has also known much success in functional data analysis, where it measures the centrality of a function with respect to a sample of functional data. Some instances are the band depth (\citealp{LopRom2009}), the functional halfspace depth (\citealp{Claetal2014}) and the functional spatial depth (\citealp{ChaCha2014b}). The large variety of available depths made it necessary to introduce an axiomatic approach identifying the most desirable properties of a depth function; see \cite{ZuoSer2000A} in the multivariate case and \cite{NieBat2016} in the functional one. 
 
Statistical depth provides a center-outward ordering of the observations that allows to tackle in a robust and nonparametric way a broad range of inference problems; see \cite{Liuetal1999}. For most depths, the deepest point is a robust location functional that extends the univariate median to the multivariate or functional setups; see, in particular, \cite{Caretal2017} for a recent work on the functional spatial median. Beyond the median, depth plays a key role in the classical problem of defining multivariate quantiles; see, e.g., \cite{Haletal2010} or \cite{Ser2010}. In line with this, the  collections of locations in~$\R^k$ whose depth does not exceed a given level are sometimes called \emph{quantile regions}; see, e.g., \cite{HeEin2016} in a multivariate extreme value theory framework. In the functional case, the quantiles in~\cite{Cha1996} may be seen as those associated with functional spatial depth; see \cite{ChaCha2014b}. Both in the multivariate and functional cases, supervised classification and outlier detection are standard applications of depth; we refer, e.g., to \cite{Cueetal2007}, \cite{PaiVanB2012}, \cite{Ser2010z}, \cite{Hubetal2015} and to the references therein.
 
In \cite{Miz2002}, statistical depth was extended to a virtually arbitrary parametric framework. In a generic parametric model indexed by an \mbox{$\ell$-dimensional} parameter~$\vartheta$, the resulting \emph{tangent depth}~$D_{P_n}(\vartheta_0)$ measures how appropriate a parameter value~$\vartheta_0$ is, with respect to the empirical measure~$P_n$ of a sample of $k$-variate observations~$X_1,\ldots,X_n$ at hand, as one could alternatively do based on the likelihood~$L_{P_n}(\vartheta_0)$. Unlike the MLE of~$\vartheta$, the depth-based estimator maximizing~$D_{P_n}(\vartheta)$ is robust under mild conditions; see Section~4 of \cite{Miz2002}. The construction, that for linear regression provides the \cite{RouHub1999} depth, proved useful in various contexts. However, tangent depth requires evaluating the halfspace depth of a given location in~$\mathbb{R}^\ell$, hence can only deal with low-dimensional parameters. In particular, tangent depth cannot cope with covariance or scatter matrix parameters ($\ell=k(k+1)/2$), unless~$k$ is as small as~$2$ or~$3$.

The crucial role played by scatter matrices in multivariate statistics, however, makes it highly desirable to have a satisfactory depth for such parameters, as phrased by \cite{Ser2004b}, that calls for an extension of the \cite{MizMul2004} location-scale depth concept into a location-scatter one. While computational issues prevent from basing this extension on tangent depth, a more ad hoc approach such as the one proposed in \cite{Zha2002} is suitable. Recently, another concept of scatter depth, that is very close in spirit to the one from \cite{Zha2002}, was introduced in \cite{Chenetal16}. Both proposals dominate tangent depth in the sense that, for $k$-variate observations, they rely on projection pursuit in~$\R^k$ rather than in~$\R^{k(k+1)/2}$, which allowed \cite{Chenetal16} to consider their depth even in high dimensions, under, e.g.,  sparsity assumptions. Both works, however, mainly focus on asymptotic, robustness and/or minimax convergence properties of the sample deepest scatter matrix. The properties of these scatter depths thus remain largely unknown, which severely affects the interpretation of the sample concepts. 

In the present work, we consider a concept of halfspace depth for scatter matrices that is close to the \cite{Zha2002} and \cite{Chenetal16} ones. Unlike these previous works, however, we thoroughly study the properties of the scatter depth and of the corresponding depth regions. We do so under minimal assumptions and, in particular, we do not restrict to elliptical distributions nor to absolutely continuous distributions. Interestingly, fully understanding scatter halfspace depth requires considering different geometries/topologies on the space of scatter matrices. Like \cite{DonGas1992} and \cite{RouRut1999} did for location halfspace depth, we study continuity and quasi-concavity properties of scatter halfspace depth, as well as the boundedness, convexity and compacity properties of the corresponding depth regions. Existence of a deepest halfspace scatter matrix, which is not guaranteed a priori, is also investigated. We further discuss, in the spirit of \cite{ZuoSer2000A}, the structural properties a scatter depth should satisfy and we investigate whether or not these are met by scatter halfspace depth. Moreover, companion concepts of depth for concentration matrices and shape matrices are proposed and studied. To the best of our knowledge, our results are the first providing structural and topological properties of depth regions outside the classical location framework. Throughout, numerical results illustrate our theoretical findings. Finally, we show the practical relevance of the depth concepts considered in a real-data example from finance.
 
The outline of the paper is as follows. In Section~\ref{SecScatter}, we define scatter halfspace depth and investigate its affine invariance and uniform consistency. We also obtain explicit expressions of this depth for two distributions we will use as running examples in the paper. In Section~\ref{SecFrobenius}, we derive the properties of scatter halfspace depth and scatter halfspace depth regions when considering the Frobenius topology on the space of scatter matrices, whereas we do the same for the geodesic topology in Section~\ref{SecGeod}. In Section~\ref{SecAxiom}, we identify the desirable properties a generic scatter depth should satisfy and investigate whether or not these are met by scatter halfspace depth. In Sections~\ref{Secconcentration} and~\ref{Secshape}, we extend this depth to concentration and shape matrices, respectively. In Section~\ref{Secrealdata}, we treat a real-data example from finance. Final comments and perspectives for future work are provided in Section~\ref{Secfinal}. Proofs and further numerical results are provided in the appendix. 

Before proceeding, we list here, for the sake of convenience, some notation to be used throughout. The  collection of $k\times k$ matrices, $k\times k$ invertible matrices, and $k\times k$ symmetric matrices will be denoted as~$\mathcal{M}_k$, $GL_k$, and~$\mathcal{S}_k$, respectively (all matrices in this paper are real matrices). The identity matrix in~$\mathcal{M}_k$ will be denoted as~$I_k$. For any~$A\in\mathcal{M}_k$, ${\rm diag}(A)$ will stand for the $k$-vector collecting the diagonal entries of~$A$, whereas, for any $k$-vector~$v$, ${\rm diag}(v)$ will stand for the diagonal matrix such that~${\rm diag}({\rm diag}(v))=v$. For~$p\geq 2$ square matrices~$A_1,\ldots,A_p$, ${\rm diag}(A_1,\ldots,A_p)$ will stand for the block-diagonal matrix with diagonal blocks~$A_1,\ldots,A_p$. Any matrix~$A$ in~$\mathcal{S}_k$ can be diagonalized into~$A=O\, {\rm diag}(\lambda_1(A),\ldots,\lambda_k(A))\, O'$, where $\lambda_1(A)\geq\ldots\geq\lambda_k(A)$ are the eigenvalues of~$A$ and where the columns of the $k\times k$ orthogonal matrix~$O=(v_1(A),\ldots,v_k(A))$ are corresponding unit eigenvectors (as usual, eigenvectors, and possibly eigenvalues, are only partly identified, but this will not play a role in the sequel). The spectral interval of~$A$ is~${\rm Sp}(A):=[\lambda_k(A),\lambda_1(A)]$. For any mapping~$f:\R\to\R$, we let~$f(A)=O\, {\rm diag}(f(\lambda_1(A)),\ldots,f(\lambda_k(A)))\, O'$. If~$\Sigma$ is a \emph{scatter matrix}, in the sense that~$\Sigma$ belongs to the collection~$\mathcal{P}_k$ of symmetric and positive definite $k\times k$ matrices, then this defines $\log(\Sigma)$ and~$\Sigma^{t}$ for any~$t\in\R$. In particular, $\Sigma^{1/2}$ is the unique~$A\in\mathcal{P}_k$ such that~$\Sigma=AA'$, and~$\Sigma^{-1/2}$ is the inverse of this symmetric and positive definite square root.
Throughout,~$T$ will denote a location functional, that is, a function mapping a probability measure~$P$ to a real $k$-vector~$T_P$. A location functional~$T$ is affine-equivariant if~$T_{P_{A,b}}=AT_P+b$ for any~$A\in GL_k$ and~$b\in\R^k$, where the probability measure~$P_{A,b}$ \label{pagedefPAb} is the distribution of~$AX+b$ when~$X$ has distribution~$P$. A much weaker equivariance concept is centro-equivariance, for which~$T_{P_{A,b}}=AT_P+b$ is imposed for~$A=-I_k$ and~$b=0$ only. For a probability measure~$P$ over~$\R^k$ and a location functional~$T$, we will let~$\alpha_{P,T}:=\min(s_{P,T},1-s_{P,T})$, where~$s_{P,T}:=\sup_{u\in \mathcal{S}^{k-1}} P[\{x\in\R^k:u'(x-T_P)=0\}]$ involves the unit sphere~$\mathcal{S}^{k-1}:=\{x\in\R^k:\|x\|^2=x'x=1\}$ of~$\R^k$. We will say that~$P$ is \emph{smooth at~$\theta$}($\in\R^k$) if the $P$-probability of any \vspace{-1.2mm}
hyperplane of~$\R^k$ containing~$\theta$ is zero 
and that it is \emph{smooth} if it is smooth at any~$\theta$. Finally, $\stackrel{\mathcal{D}}{=}$ will denote equality in distribution.

\section{Scatter halfspace depth}
\label{SecScatter}

We start by recalling the classical concept of location halfspace depth. To do so, let~$P$ be a probability measure over~$\R^k$ and~$X$ be a random $k$-vector with distribution~$P$, which allows us  throughout to write~$P[X\in B]$ instead of~$P[B]$ for any $k$-Borel set~$B$. The \emph{location halfspace depth} of~$\theta(\in\R^k)$ with respect to~$P$ is then
$$
H\!D_P^{\rm loc}(\theta)
:=
\inf_{u\in\mathcal{S}^{k-1}} 
P[ u'(X-\theta) \geq 0 ]
.	
$$ 
The corresponding \emph{depth regions}~$R_{P}^{\rm loc}(\alpha):=\{\theta\in\R^k: H\!D_P^{\rm loc}(\theta)\geq \alpha\}$ form a nested family of closed convex subsets of~$\R^k$. The innermost depth region, namely~$M_P^{\rm loc}:=\{\theta\in\R^k: H\!D^{\rm loc}_P(\theta)=\max_{\eta \in\R^k}H\!D^{\rm loc}_P(\eta)\}$ (the maximum always exists; see, e.g., Proposition~7 in \citealp{RouRut1999}), is a set-valued location functional. When a unique representative of~$M_P^{\rm loc}$ is needed, it is customary to consider the \emph{Tukey median}~$\theta_P$ of~$P$, that is defined as the barycenter of~$M_P^{\rm loc}$. The Tukey median has maximal depth (which follows from the convexity of~$M_P^{\rm loc}$) and is an affine-equivariant location functional.

In this paper, for a location functional~$T$, we define the \emph{$T$-scatter halfspace depth} of~$\Sigma(\in\mathcal{P}_k)$ with respect to~$P$ as
\begin{eqnarray}
\lefteqn{
H\!D^{\rm sc}_{P,T}(\Sigma)
:=
\!
\inf_{u\in\mathcal{S}^{k-1}} 
\min
\!
\big( 
P\big[ |u'(X-T_P)| \leq \sqrt{u'\Sigma u}\, \big]
,
}
\nonumber
\\[2mm]
& & 
\hspace{43mm} 
P\big[ |u'(X-T_P)| \geq \sqrt{u'\Sigma u}\, \big]
\big)
.
\label{def}
\end{eqnarray}

This extends to a probability measure with arbitrary location the centered matrix depth concept from \cite{Chenetal16}. If~$P$ is smooth, then the depth in~(\ref{def}) is also equivalent to the (Tukey version of) the dispersion depth introduced in \cite{Zha2002}, but for the fact that the latter, in the spirit of projection depth, involves centering through a univariate location functional (both \cite{Zha2002} and \cite{Chenetal16} also propose bypassing centering through a pairwise difference approach that will be discussed in Section~\ref{Secfinal}). 
While they were not considered in these prior works, it is of interest to introduce the corresponding depth regions
\begin{equation}
\label{defbisbis}
R^{\rm sc}_{P,T}(\alpha)
:=
\big\{\Sigma\in\mathcal{P}_k : H\!D^{\rm sc}_{P,T}(\Sigma)\geq \alpha \big\}
,
\quad
\alpha\geq 0
.
\end{equation}
We will refer to~$R^{\rm sc}_{P,T}(\alpha)$ as the \emph{order-$\alpha$ $(T$-scatter halfspace$)$ depth region of~$P$}. Obviously, one always has~$R^{\rm sc}_{P,T}(0)=\mathcal{P}_k$. Clearly, the concepts in (\ref{def})-(\ref{defbisbis}) give practitioners the flexibility to freely choose the location functional~$T$;  numerical results below, however, will focus on the depth~$H\!D^{{\rm sc}}_P(\Sigma)$ and on the depth regions~$R^{\rm sc}_{P}(\alpha)$ based on the Tukey median~$\theta_P$, that is the natural location functional whenever halfspace depth objects are considered. 

To get a grasp of the scatter depth~$H\!D^{\rm sc}_P(\Sigma)$, it is helpful to start with the univariate case~$k=1$. There, the location halfspace deepest region is the ``median interval"~$
M^{\rm loc}_P=
\arg\max_{\theta\in \R} \min(P[ X \leq \theta ],P[ X \geq \theta ])$ and the Tukey median~$\theta_P$, that is, the midpoint of~$M_P^{\rm loc}$, is the usual representative of the univariate median. The scatter halfspace deepest region is then the median 
\vspace{-.7mm}
 interval~$
M^{\rm sc}_P
:=
\arg\max_{\Sigma\in \R^+_0}
 \min(
P[ (X-\theta_P)^2 \leq \Sigma\, ]
,
P[ (X-\theta_P)^2 \geq \Sigma\, ]
)
$ of~$(X-\theta_P)^2$; call it the \emph{median squared deviation} interval~$\mathcal{I}_{\rm MSD}[X]$ (or~$\mathcal{I}_{\rm MSD}[P]$) of~$X\sim P$. Below, parallel to what is done for the median, ${\rm MSD}[X]$ (or~${\rm MSD}[P]$) will denote the midpoint of this MSD interval. In particular, if~$\mathcal{I}_{\rm MSD}[P]$ is a singleton, then scatter halfspace depth is uniquely maximized at~$\Sigma={\rm MSD}[P]=({\rm MAD}[P])^2$, where~${\rm MAD}[P]$ denotes the median absolute deviation of~$P$. Obviously, the depth regions~$R^{\rm sc}_P(\alpha)$ form a family of nested intervals,~$[\Sigma^-_\alpha,\Sigma_\alpha^+]$ say, included in~$\mathcal{P}_1=\R^+_0$. It is easy to check that, if~$P$ is symmetric about zero
with an invertible cumulative distribution function~$F$ and if~$T$ is centro-equivariant, then 
\begin{eqnarray}
	\lefteqn{
	H\!D^{\rm sc}_P(\Sigma)
=
H\!D^{\rm sc}_{P,T}(\Sigma)
=
2 
\min
\big( 
F(
\sqrt{\Sigma}) - {\textstyle \frac{1}{2}} , 1-F(
\sqrt{\Sigma})
\big)
\quad\textrm{ and}
}
\label{scduniv}
\\[2mm]
& & 
\hspace{3mm} 
R^{\rm sc}_P(\alpha)
=
R^{\rm sc}_{P,T}(\alpha)
=
\big[
(F^{-1}({\textstyle \frac{1}{2}+\frac{\alpha}{2}})
)^2
,
(F^{-1}(1-{\textstyle \frac{\alpha}{2}})
)^2
\big]
.
\label{screguniv}	
\end{eqnarray}
This is compatible with the fact that the maximal value of~$\Sigma\mapsto H\!D^{\rm sc}_P(\Sigma)$ (that is equal to~$1/2$) is achieved at~$\Sigma=({\rm MAD}[P])^2$ only. 

For~$k>1$, elliptical distributions provide an important particular case. We will say that~$P=P^X$ is $k$-variate
\vspace{-.7mm}
  elliptical with location~$\theta(\in\R^k)$ and scatter~$\Sigma(\in\mathcal{P}_k)$ if and only
  \vspace{-.7mm}
 if~$X\stackrel{\mathcal{D}}{=}\theta+\Sigma^{1/2} Z$, where~$Z=(Z_1,\ldots,Z_k)'$ is (i) spherically symmetric about the origin of~$\R^k$ (that is,~$OZ\stackrel{\mathcal{D}}{=}Z$ for any $k\times k$ orthogonal matrix~$O$) and is (ii) standardized in such a way that~${\rm MSD}[Z_1]=1$ (one then has~$T_P=\theta$ for any affine-equivariant location functional~$T$). Denoting by~$\Phi$ the cumulative distribution function of the standard normal, the $k$-variate normal distribution with location zero and scatter~$I_k$ is then the distribution of~$X:=W/b$, where~$b:=\Phi^{-1}(\frac{3}{4})$ and~$W$ is a standard normal random $k$-vector. In this Gaussian case, we obtain
\begin{eqnarray} 
H\!D_{P,T}^{\rm sc}(\Sigma)
&\!\!=\!\!&
\inf_{u\in\mathcal{S}^{k-1}} 
\min
\big( 
P\big[ |u'X| \leq \sqrt{u' \Sigma u}\, \big]
,
P\big[ |u'X| \geq \sqrt{u' \Sigma u}\, \big]
\big)
\nonumber
\\[2mm]
&\!\!=\!\!&
2
\min\Big(
\Phi\big(b \lambda^{1/2}_k(\Sigma)\big) - {\textstyle\frac{1}{2}}
,
1-\Phi\big(b \lambda^{1/2}_1(\Sigma)\big)
\Big)
.
\label{normaltheo}	
\end{eqnarray}
One can check directly that~$H\!D_{P,T}^{\rm sc}(\Sigma)\leq H\!D_{P,T}^{\rm sc}(I_k)=1/2$, with equality if and only if~$\Sigma$ coincides with the ``true" scatter matrix~$I_k$ (we refer to Theorem~\ref{Fishconst} for a more general result). Also,~$\Sigma$ belongs to the depth region~$R^{\rm sc}_{P,T}(\alpha)$ if and only if 
$
{\rm Sp}(\Sigma)
\subset 
[
(
{\textstyle\frac{1}{b}}
\Phi^{-1}({\textstyle\frac{1}{2}+\frac{\alpha}{2}}))^2
,
(
{\textstyle\frac{1}{b}}
\Phi^{-1}(1-{\textstyle\frac{\alpha}{2}}))^2
]
.
$

Provided that the location functional used is affine-equivariant, extension to an arbitrary multinormal is based on the following affine-invariance result, which ensures in particular that scatter halfspace depth will not be affected by possible changes in the marginal measurement units (a similar result is stated in \cite{Zha2002} for the dispersion depth concept considered there).

\begin{Theor}\label{thaffinv}
Let~$T$ be an affine-equivariant location functional.  Then, 
(i) scatter halfspace depth is affine-invariant in the sense that, for any probability measure~$P$ over~$\R^k$, $\Sigma\in\mathcal{P}_k$, $A\in GL_k$ and~$b\in\R^k$, we have
$
H\!D^{\rm sc}_{P_{A,b},T}(A\Sigma A')
=
H\!D_{P,T}^{\rm sc}(\Sigma)
,
$
where~$P_{A,b}$ is as defined on page~\pageref{pagedefPAb}. Consequently, (ii) the regions~$R^{\rm sc}_{P,T}(\alpha)$ are affine-equivariant, in the sense that, for any probability measure~$P$ over~$\R^k$, $\alpha\geq 0$, $A\in GL_k$ and~$b\in\R^k$, we have~$R^{\rm sc}_{P_{A,b},T}(\alpha)=AR^{\rm sc}_{P,T}(\alpha)A'$.
\end{Theor}

This result readily entails that if~$P$ is the $k$-variate normal with location~$\theta_0$ and scatter~$\Sigma_0$, then, provided that~$T$ is affine-equivariant, 
\begin{equation}
\label{Gaussiandepth}
H\!D_{P,T}^{\rm sc}(\Sigma)
=
2
\min\Big(
\Phi\big(b \lambda^{1/2}_k(\Sigma_0^{-1}\Sigma)\big) - {\textstyle\frac{1}{2}}
,
1-\Phi\big(b \lambda^{1/2}_1(\Sigma_0^{-1}\Sigma)\big)
\Big)
\end{equation}
and~$R^{\rm sc}_{P,T}(\alpha)$ is the collection of scatter matrices~$\Sigma$ for which 
$
{\rm Sp}(\Sigma_0^{-1}\Sigma)
\subset 
[(
{\textstyle\frac{1}{b}}
\Phi^{-1}({\textstyle\frac{1}{2}+\frac{\alpha}{2}}))^2
,
(
{\textstyle\frac{1}{b}}
\Phi^{-1}(1-{\textstyle\frac{\alpha}{2}}))^2
]
.
$
For a non-Gaussian elliptical probability measure~$P$ with
\vspace{-.8mm}
  location~$\theta_0$ and scatter~$\Sigma_0$, it is easy to show that $H\!D_{P,T}^{\rm sc}(\Sigma)$ will still depend on~$\Sigma$ only through~$\lambda_1(\Sigma_0^{-1}\Sigma)$ and~$\lambda_k(\Sigma_0^{-1}\Sigma)$. 

\begin{figure}[h!]
\includegraphics[width=\textwidth]{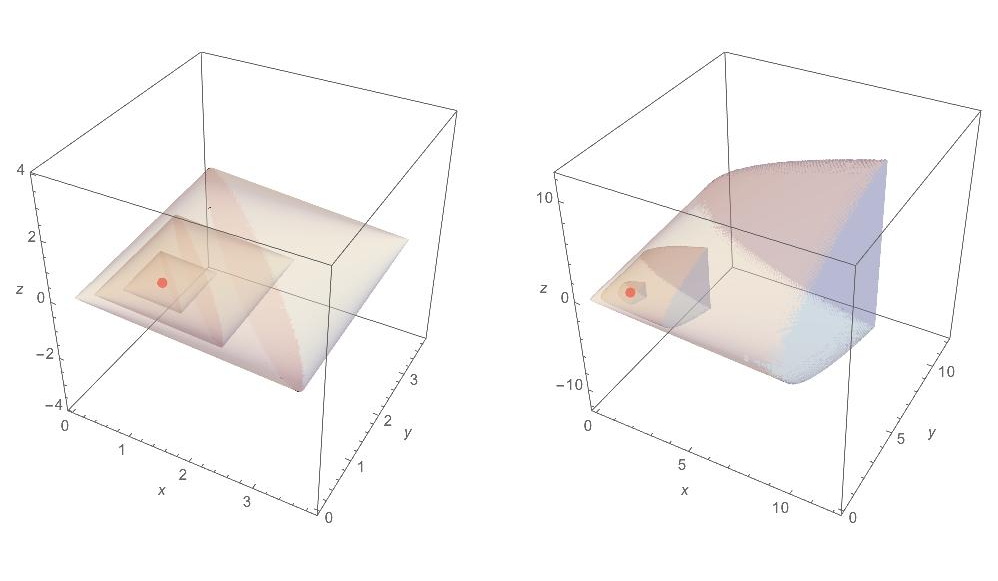}   
\vspace{-6mm}
 \caption{Level sets of order~$\alpha=.2,.3$ and~$.4$, for any centro-symmetric~$T$, of~$(x,y,z)\mapsto H\!D^{\rm sc}_{P,T}(\Sigma_{x,y,z})$, where~$H\!D^{\rm sc}_{P,T}(\Sigma_{x,y,z})$ is the $T$-scatter halfspace depth of~$\Sigma_{x,y,z}={\hspace{-0mm} x\ \hspace{0mm} z\, \choose \,z\ \hspace{0mm} y\,}$ with respect to two probability measures~$P$, namely the bivariate multinormal distribution with location zero and scatter~$I_2$ (left) and the bivariate distribution with independent Cauchy marginals (right). The red points are those associated with~$I_2$ (left) and~$\sqrt{2}I_2$ (right), which are the corresponding deepest scatter matrices (see Sections~\ref{SecGeod} and~\ref{SecAxiom}). 
}
\label{FigSec2a}  
\end{figure}

As already mentioned, we also intend to consider non-elliptical probability measures. A running non-elliptical example will be the one for which~$P$ is the distribution of a random vector~$X=(X_1,\ldots,X_k)'$ with independent Cauchy marginals. If~$T$ is centro-equiva\-riant, then \begin{equation}
\label{Cauchytheo}	
H\!D_{P,T}^{\rm sc}(\Sigma)
=
2
\min
\Big(
\Psi\big(  1/\max_s \sqrt{s'\Sigma^{-1}s}\, \big) - 
{\textstyle\frac{1}{2}}
,
1 - \Psi\big( \sqrt{\max({\rm diag}(\Sigma))}  \big) 
\Big)
,
\end{equation}
where~$\Psi$ is the Cauchy cumulative distribution function and where the maximum in~$s$ is over all sign vectors~$s=(s_1,\ldots,s_k)\in\{-1,1\}^k$; we establish the explicit expression~(\ref{Cauchytheo}) in Appendix~\ref{AppSecScatter}. For~$k=1$, this 
simplifies to
$
H\!D_{P,T}^{\rm sc}(\Sigma)
=
2
\min
\big(
\Psi\big( \sqrt{\Sigma}\, \big) - 
{\textstyle\frac{1}{2}}
,
1 - \Psi\big( \sqrt{\Sigma}  \big) 
\big)
,
$
which agrees with~(\ref{scduniv}). For~$k=2$, we obtain
$$
H\!D_{P,T}^{\rm sc}(\Sigma)
=
2
\min
\Big(
\Psi\big(  \sqrt{\det(\Sigma)/s_\Sigma}\, \big) - 
{\textstyle\frac{1}{2}}
,
1 - \Psi\big( \sqrt{\max(\Sigma_{11},\Sigma_{22})}  \big) 
\Big)
,
$$
where we let~$s_\Sigma:=\Sigma_{11}+\Sigma_{22}+2|\Sigma_{12}|$. For a general~$k$, a scatter matrix~$\Sigma$ belongs to~$R^{\rm sc}_{P,T}(\alpha)$ if and only if
$
1/(s'\Sigma^{-1}s)
\geq
(\Psi^{-1}( {\textstyle\frac{1}{2}}+{\textstyle\frac{\alpha}{2}}))^2
$
for all $s\in\{-1,1\}^k$
and
$
\Sigma_{\ell\ell}
\leq
(\Psi^{-1}\big( 1-{\textstyle\frac{\alpha}{2}}))^2
$
for all $\ell=1,\ldots,k$. The problem of identifying the scatter matrix achieving maximal depth, if any (existence is not guaranteed), will be considered in Section~\ref{SecGeod}. Figure~\ref{FigSec2a} plots scatter halfspace depth regions in the Gaussian and independent Cauchy cases above. Examples involving distributions that are not absolutely continuous with respect to the Lebesgue measure will be considered in the next sections.

In Appendix~\ref{SecSupMonteCarlo}, we validate through a Monte Carlo exercise the expressions for~$H\!D_{P,T}^{\rm sc}(\Sigma)$ obtained in~(\ref{Gaussiandepth})-(\ref{Cauchytheo}) above. Such a numerical validation is justified by the following uniform consistency result; see~(6.2) and~(6.6) in \cite{DonGas1992} for the corresponding location halfspace depth result, and Proposition~2.2(ii) in \cite{Zha2002} for the dispersion depth concept considered there. 

\begin{Theor} 
\label{Consistency}
Let~$P$ be a smooth probability measure over~$\R^k$ and $T$ be a location functional. Let~$P_n$ denote the empirical probability measure associated with a random sample of size~$n$ from~$P$ and assume that~$T_{P_n}\to T_P$ almost surely as~$n\to\infty$. Then $\sup_{\Sigma\in\mathcal{P}_k}
|H\!D^{\rm sc}_{P_n,T}(\Sigma)-H\!D_{P,T}^{\rm sc}(\Sigma)|
\to 0$ almost surely as~$n\to\infty$. 
\end{Theor}
 
This result applies in particular to the scatter halfspace depth~$H\!D_{P}^{\rm sc}(\Sigma)$, as the Tukey median is strongly consistent without any assumption on~$P$ (for completeness, we show this in Lemma~\ref{LemConvergence}). 
Inspection of the proof of Theorem~\ref{Consistency} reveals that the smoothness assumption is only needed to control the estimation of~$T_P$, hence is superfluous when a constant location functional is used. 
This is relevant when the location is fixed, as in \cite{Chenetal16}.

\section{Frobenius topology} 
\label{SecFrobenius} 

Our investigation of the further structural properties of the scatter halfspace depth~$H\!D_{P,T}^{\rm sc}(\Sigma)$ and of the corresponding depth regions~$R^{\rm sc}_{P,T}(\alpha)$ depends on the topology that is considered on~$\mathcal{P}_k$. In this section, we focus on the topology induced by the \emph{Frobenius metric space}~$(\mathcal{P}_k,d_F)$, where $d_F(\Sigma_a,\Sigma_b)=\|\Sigma_b-\Sigma_a\|_F$ is the distance on~$\mathcal{P}_k$ that is inherited from the Frobenius norm~$\|A\|_F=\sqrt{{\rm tr}[AA']}$ on~$\mathcal{M}_k$. The resulting \emph{Frobenius topology} (or simply \emph{$F$-topology}), generated by the $F$-balls~$B_F(\Sigma_0,r):=\{\Sigma\in\mathcal{P}_k : d_F(\Sigma,\Sigma_0)< r \}$ with center~$\Sigma_0$ and radius~$r$, gives a precise meaning to what we call below \emph{$F$-continuous} functions on~$\mathcal{P}_k$, \emph{$F$-open/$F$-closed} subsets of~$\mathcal{P}_k$, etc. We then have the following result.

\begin{Theor}
\label{Propcontinuity}
Let~$P$ be a probability measure over~$\R^k$ and $T$ be a location functional. Then, (i) $\Sigma\mapsto H\!D_{P,T}^{\rm sc}(\Sigma)$ is upper $F$-semicontinuous on~$\mathcal{P}_k$, so that (ii) the depth region~$R^{\rm sc}_{P,T}(\alpha)$ is $F$-closed for any~$\alpha\geq 0$. (iii) If~$P$ is smooth at~$T_P$, then~$\Sigma\mapsto H\!D_{P,T}^{\rm sc}(\Sigma)$ is $F$-continuous on~$\mathcal{P}_k$.  
\end{Theor}
\vspace{-1mm}

For location halfspace depth, the corresponding result was derived in Lemma~6.1 of \cite{DonGas1992}, where the metric on~$\R^k$ is the Euclidean one. The similarity between the location and scatter halfspace depths also extends to the boundedness of depth regions, in the sense that, like location halfspace depth (Proposition~5 in \citealp{RouRut1999}), the order-$\alpha$ scatter halfspace depth region is bounded if and only if~$\alpha>0$.

\begin{Theor}
\label{boundedness}
Let~$P$ be a probability measure over~$\R^k$ and $T$ be a location functional. Then, for any~$\alpha>0$, $R^{\rm sc}_{P,T}(\alpha)$ is $F$-bounded $($that is, it is included, for some~$r>0$, in the $F$-ball~$B_F(I_k,r))$. 
\end{Theor}

This shows that, for any probability measure~$P$, $H\!D_{P,T}^{\rm sc}(\Sigma)$ goes to zero \mbox{as~$\|\Sigma\|_F\to\infty$.} Since $\|\Sigma\|_F \geq \lambda_1(\Sigma)$, this means that explosion of~$\Sigma$ (that is, $\lambda_1(\Sigma)\to \infty$) leads to arbitrarily small depth, which is confirmed in the multinormal case in~(\ref{normaltheo}). In this Gaussian case, however, implosion of~$\Sigma$ (that is, $\lambda_k(\Sigma)\to 0$) also provides arbitrarily small depth, but this is not captured by the general result in Theorem~\ref{boundedness} (similar comments can be given for the independent Cauchy example in~(\ref{Cauchytheo})). Irrespective of the topology adopted (so that the $F$-topology is not to be blamed for this behavior), it is actually possible to have implosion without depth going to zero. We show this by considering the following example. Let~$P=(1-s)P_1+sP_2$, where~$s\in(\frac{1}{2},1)$, $P_1$ is the bivariate standard normal 
\vspace{-.5mm}
 and $P_2$ is the distribution of~${0 \choose Z}$, where~$Z$ is univariate standard normal. Then, it can be showed that, for~$\Sigma_n:={\,1/n\  0\, \choose \ \ 0\ \ 1\,}$ and any centro-equivariant~$T$, we have $H\!D^{\rm sc}_{P,T}(\Sigma_n)\to 1-s>0$ as~$n\to\infty$. 

In the metric space~$(\mathcal{P}_k,d_F)$, any bounded set is also \emph{totally bounded}, that is, can be covered, for any~$\varepsilon>0$, by finitely many balls of the form~$B_F(\Sigma,\varepsilon)$. Theorems~\ref{Propcontinuity}-\ref{boundedness} thus show that, for any~$\alpha>0$, $R^{\rm sc}_{P,T}(\alpha)$ is both $F$-closed and totally $F$-bounded. However, since~$(\mathcal{P}_k,d_F)$ is not complete, there is no guarantee that these regions are $F$-compact. Actually, these regions may fail to be $F$-compact, as we show through the example from the previous paragraph. For any~$\alpha\in (0,1-s)$, the scatter matrix~$\Sigma_n$ belongs to~$R^{\rm sc}_{P,T}(\alpha)$ for~$n$ large enough. However, the sequence~$(\Sigma_n)$ $F$-converges to~${\,0\ 0\, \choose \,0\ 1\,}$, that does not belong to~$R^{\rm sc}_{P,T}(\alpha)$ (since it does not even belong to~$\mathcal{P}_2$). Since this will also hold for any subsequence of~$(\Sigma_n)$, we conclude that, for~$\alpha\in (0,1-s)$, $R^{\rm sc}_{P,T}(\alpha)$ is not $F$-compact in this example. This provides a first discrepancy between location and scatter halfspace depths, since location halfspace depth regions associated with a positive order~$\alpha$ are always compact. 

The lack of compacity of scatter halfspace depth regions may allow for probability measures for which no halfspace deepest scatter exists. This is actually the case in the bivariate mixture example above. There, letting~$e_1=(1,0)'$ and assuming again that~$T$ is centro-equivariant, any~$\Sigma\in\mathcal{P}_2$ indeed satisfies~$H\!D_{P,T}^{\rm sc}(\Sigma)\leq P[ |e_1'X| \geq \sqrt{e_1'\Sigma e_1}]=P[ |X_1| \geq \sqrt{\Sigma_{11}} ]=(1-s) P[ |Z| \geq \sqrt{\Sigma_{11}}]<1-s=\sup_{\Sigma\in\mathcal{P}_2} H\!D_{P,T}^{\rm sc}(\Sigma)$, where the last equality follows from the fact that we identified a sequence~$(\Sigma_n)$ such that~$H\!D_{P,T}^{\rm sc}(\Sigma_n)\to 1-s$. This is again in sharp contrast with the location case, for which a halfspace deepest location always exists; see, e.g., Propositions~5 and~7 in \cite{RouRut1999}. Identifying sufficient conditions under which a halfspace deepest scatter exists requires considering another topology, namely the \emph{geodesic topology} considered in Section~\ref{SecGeod} below. 

The next result states that scatter halfspace depth is a quasi-concave function, which ensures convexity of the corresponding depth regions; we refer to Proposition~1 (and to its corollary) in \cite{RouRut1999} for the corresponding results on location halfspace depth.  

\begin{Theor}
\label{Propquasiconcavity}
Let~$P$ be a probability measure over~$\R^k$ and $T$ be a location functional. Then, (i) $\Sigma\mapsto H\!D_{P,T}^{\rm sc}(\Sigma)$ is quasi-concave, in the sense that, for any~$\Sigma_a,\Sigma_b\in\mathcal{P}_k$ and~$t\in [0,1]$, $H\!D^{\rm sc}_{P,T}(\Sigma_t)\!\geq\! \min(H\!D^{\rm sc}_{P,T}(\Sigma_a),H\!D^{\rm sc}_{P,T}(\Sigma_b))$, where we let~$\Sigma_t:=(1-t)\Sigma_a+t \Sigma_b$; (ii) for any~$\alpha\geq 0$, $R^{\rm sc}_{P,T}(\alpha)$ is convex.
\end{Theor}

Strictly speaking, Theorem~\ref{Propquasiconcavity} is not directly related to the $F$-topology considered on~$\mathcal{P}_k$. Yet we state the result in this section due to the link between the linear paths~$t\mapsto \Sigma_t=(1-t)\Sigma_a+t \Sigma_b$ it involves and the ``flat" nature of the $F$-topology (this link will become clearer below when we will compare with what occurs for the geodesic topology). Illustration of Theorem~\ref{Propquasiconcavity} will be provided in Figure~\ref{FigSec5} below, as well as in Appendix~\ref{SecSupQuasi}.


\section{Geodesic topology}
\label{SecGeod}


Equipped with the inner product~$<\!A,B\!>\,={\rm tr}[A'B]$, $\mathcal{M}_k$ is a Hilbert space. The resulting norm and distance are the Frobenius ones considered in the previous section. As an open set in~$\mathcal{S}_k$, the parameter space~$\mathcal{P}_k$ of interest is a differentiable manifold of dimension~$k(k+1)/2$. The corresponding tangent space at~$\Sigma$, which is isomorphic (via translation) to~$\mathcal{S}_k$, can be equipped with the inner product~$\textrm{\mbox{$<\!A,B\!>$}}\,={\rm tr}[\Sigma^{-1}A\Sigma^{-1}B]$. This leads to considering~$\mathcal{P}_k$ as a Riemannian manifold, with the metric at~$\Sigma$ given by the differential
$
ds
=
\| \Sigma^{-1/2} d\Sigma\, \Sigma^{-1/2} \|_F
;
$
see, e.g., \cite{Bha07}. The length of a path~$\gamma:[0,1]\to \mathcal{P}_k$ is then given by
$$
L(\gamma)
=
\int_0^1
\Big\|
\gamma^{-1/2}(t) \frac{d\gamma(t)}{dt} \gamma^{-1/2}(t)
\Big\|_F
\,dt
.
$$ 
The resulting geodesic distance between~$\Sigma_a,\Sigma_b\in\mathcal{P}_k$ is defined as
\begin{equation}
\label{defindist}	
d_g(\Sigma_a,\Sigma_b)
:=
\inf\big\{ L(\gamma) : \gamma\in \mathcal{G}(\Sigma_a,\Sigma_b) \big\}
=
\| \log( \Sigma_a^{-1/2} \Sigma_b \Sigma_a^{-1/2} ) \|_F
,
\end{equation}
where~$\mathcal{G}(\Sigma_a,\Sigma_b)$ denotes the collection of paths~$\gamma$ from~$\gamma(0)=\Sigma_a$ to~$\gamma(1)=\Sigma_b$ (the second equality in~(\ref{defindist}) is Theorem~6.1.6 in \citealp{Bha07}). It directly follows from the definition of~$d_g(\Sigma_a,\Sigma_b)$ that the geodesic distance satisfies the triangle inequality. Theorem~6.1.6 in \cite{Bha07} also states that all paths~$\gamma$ achieving the infimum in~(\ref{defindist}) provide the same geodesic~$\{\gamma(t):t\in[0,1]\}$ joining~$\Sigma_a$ and~$\Sigma_b$, and that this geodesic can be parametrized as 
\begin{equation}
\label{geod}
\gamma(t)
=
\tilde{\Sigma}_t
:=
\Sigma_a^{1/2} 
\big( \Sigma_a^{-1/2} \Sigma_b \Sigma_a^{-1/2} \big)^t 
\Sigma_a^{1/2}
,
\qquad
t \in [0,1]
.
\end{equation}
By using the explicit formula in~(\ref{defindist}), it is easy to check that this particular parametrization of this unique geodesic is natural in the sense that $d_g(\Sigma_a,\tilde{\Sigma}_t)=t d_g(\Sigma_a,\Sigma_b)$ for any~$t\in[0,1]$. 

Below, we consider the natural topology associated with the metric space $(\mathcal{P}_k,d_g)$, that is, the topology whose open sets are generated by geodesic balls of the form~$B_g(\Sigma_0,r):=\{ \Sigma\in\mathcal{P}_k : d_g(\Sigma,\Sigma_0)<r\}$. This topology --- call it the \emph{geodesic topology}, or simply \emph{$g$-topology} --- defines subsets of~$\mathcal{P}_k$ that are~$g$-open, $g$-closed, $g$-compact, and functions that are $g$-semicontinuous, $g$-continuous, etc. We will say that a subset~$R$ of~$\mathcal{P}_k$ is $g$-bounded if and only if~$R\subset B_g(I_k,r)$ for some~$r>0$ (we can safely restrict to balls centered at~$I_k$ since the triangle inequality guarantees that $R$ is included in a finite-radius $g$-ball centered at~$I_k$ if and only if it is included in a finite-radius $g$-ball centered at an arbitrary~$\Sigma_0\in\mathcal{P}_k$). A $g$-bounded subset of~$\mathcal{P}_k$ is also totally $g$-bounded, still in the sense that, for any~$\varepsilon>0$, it can be covered by finitely many balls of the form~$B_g(\Sigma,\varepsilon)$; for completeness, we prove this in Lemma~\ref{totalboundlemma}. Since~$(\mathcal{P}_k,d_g)$ is complete (see, e.g., Proposition~10 in \citealp{BhaHol06}), a $g$-bounded and $g$-closed subset of~$\mathcal{P}_k$ is then $g$-compact. 

We omit the proof of the next result as it follows along the exact same lines as the proof of Theorem~\ref{Propcontinuity}, once it is seen that a sequence~$(\Sigma_n)$ converging to~$\Sigma_0$ in~$(\mathcal{P}_k,d_g)$ also converges to~$\Sigma_0$ in~$(\mathcal{P}_k,d_F)$. 

\begin{Theor}
\label{geodcontinuity}
Let~$P$ be a probability measure over~$\R^k$ and $T$ be a location functional. Then, (i) $\Sigma\mapsto H\!D_{P,T}^{\rm sc}(\Sigma)$ is upper $g$-semicontinuous on~$\mathcal{P}_k$, so that (ii) the depth region~$R^{\rm sc}_{P,T}(\alpha)$ is $g$-closed for any~$\alpha\geq 0$. (iii) If~$P$ is smooth at~$T_P$, then~$\Sigma\mapsto H\!D_{P,T}^{\rm sc}(\Sigma)$ is $g$-continuous on~$\mathcal{P}_k$.  
\end{Theor}
\vspace{-1mm}
 
The following result uses the notation~$s_{P,T}:=\sup_{u\in \mathcal{S}^{k-1}} P[u'(X-T_P)=0]$ and~$\alpha_{P,T}:=\min(s_{P,T},1-s_{P,T})$ defined in the introduction.

\begin{Theor}
\label{geodboundedness}
Let~$P$ be a probability measure over~$\R^k$ and $T$ be a location functional. Then, for any~$\alpha>\alpha_{P,T}$, $R^{\rm sc}_{P,T}(\alpha)$ is $g$-bounded, hence $g$-compact $($if~$s_{P,T}\geq 1/2$, then this result is trivial in the sense that $R^{\rm sc}_{P,T}(\alpha)$ is empty for any~$\alpha>\alpha_{P,T})$. In particular,  if~$P$ is smooth at~$T_P$, then~$R^{\rm sc}_{P,T}(\alpha)$ is $g$-compact for any~$\alpha>0$.
\end{Theor}

This result complements Theorem~\ref{boundedness} by showing that implosion always leads to a depth that is smaller than or equal to~$\alpha_{P,T}$. In particular, in the multinormal and independent Cauchy examples in Section~\ref{SecScatter}, this shows that both explosion and implosion lead to arbitrarily small depth, whereas Theorem~\ref{boundedness} was predicting this collapsing for explosion only. Therefore, while the behavior of $H\!D_{P,T}^{\rm sc}(\Sigma)$ under implosion/explosion of~$\Sigma$ is independent of the topology adopted, the use of the $g$-topology  provides a better understanding of this behavior than the $F$-topology.

It is not possible to improve the result in Theorem~\ref{geodboundedness}, in the sense that $R^{\rm sc}_{P,T}(\alpha_{P,T})$ may fail to be $g$-bounded. For instance, consider the probability measure~$P$ over~$\R^2$ putting probability mass~$1/6$ on each of the six points~$(0,\pm 1/2)$ and~$(\pm 2,\pm 2)$, and let~$T$ be a centro-equivariant location functional. Clearly, $\alpha_{P,T}=s_{P,T}= 1/3$. Now, 
\vspace{-.5mm}
 letting~$\Sigma_n:={\,1/n\  0\, \choose \ \ 0\ \ 1\,}$, 
we have
$
P[ |u' \Sigma_n^{-1/2} X| \leq 1]\geq 1/3
$
and
$
P[ |u' \Sigma_n^{-1/2} X| \geq 1]\geq 1/3
$
for any~$u\in\mathcal{S}^{1}$ (here, $X$ is a random vector with  distribution~$P$), which entails that
\begin{eqnarray*}
\lefteqn{
H\!D^{\rm sc}_{P,T}(\Sigma_n)
=
\inf_{u\in\mathcal{S}^{1}} 
\min
\!
\big( 
P\big[ |u'X| \leq \sqrt{u'\Sigma_n u}\, \big]
,
P\big[ |u'X| \geq \sqrt{u'\Sigma_n u}\, \big]
\big)
}
\\[2mm]
& &
\hspace{13mm}
=
\inf_{u\in\mathcal{S}^{1}} 
\min
\!
\big( 
P[ |u' \Sigma_n^{-1/2} X| \leq 1],
P[ |u' \Sigma_n^{-1/2} X| \leq 1]
\big)
\geq 
\frac{1}{3}
=
\alpha_{P,T}
,
\end{eqnarray*}
so that~$\Sigma_n\in R^{\rm sc}_{P,T}(\alpha_{P,T})$ for any~$n$. Since~$d_g(\Sigma_n,I_2)\to\infty$, $R^{\rm sc}_{P,T}(\alpha_{P,T})$ is indeed $g$-unbounded. 

An important benefit of working with the $g$-topology is that, unlike the $F$-topology, it allows to show that, under mild assumptions, a halfspace deepest scatter does exist. More precisely, we have the following result. 

\begin{Theor}
\label{geodmaxdepth}
Let~$P$ be a probability measure over~$\R^k$ and $T$ be a location functional. Assume that~$R^{\rm sc}_{P,T}(\alpha_{P,T})$ is non-empty. Then, $\alpha_{*P,T}:=\sup_{\Sigma\in\mathcal{P}_k} H\!D_{P,T}^{\rm sc}(\Sigma)=H\!D^{\rm sc}_{P,T}(\Sigma_*)$ for some~$\Sigma_*\in\mathcal{P}_k$. 
\end{Theor}

In particular, this result shows that for any probability measure~$P$ that is smooth at~$T_P$, there exists a halfspace deepest scatter~$\Sigma_*$. For the $k$-variate multinormal distribution with location zero and scatter~$I_k$ (and any centro-equivariant~$T$), we already stated in Section~\ref{SecScatter} that~$\Sigma\mapsto H\!D_{P,T}^{\rm sc}(\Sigma)$ is uniquely maximized at~$\Sigma_*=I_k$, with a corresponding maximal depth equal to~$1/2$. The next result identifies the halfspace deepest scatter (and the corresponding maximal depth) in the independent Cauchy case. 

\begin{Theor}
\label{propmaxCauchy}
Let~$P$ be the $k$-variate probability measure with independent Cauchy marginals and let $T$ be a centro-equivariant location functional. Then,~$\Sigma\mapsto H\!D_{P,T}^{\rm sc}(\Sigma)$ 
\vspace{-.5mm}
is uniquely maximized at~$\Sigma_*=\sqrt{k}I_k$, and the corresponding maximal depth is~$H\!D^{\rm sc}_{P,T}(\Sigma_*)=\frac{2}{\pi} \arctan\big( k^{-1/4} \big)$. 
\end{Theor}

For~$k=1$, the Cauchy distribution in this result is symmetric (hence, elliptical) about zero, which is compatible with the maximal depth being equal to~$1/2$ there (Theorem~\ref{Fishconst} below shows that the maximal depth for absolutely continuous elliptical distributions is always equal to~$1/2$). For larger values of~$k$, however, this provides an example where the maximal depth is strictly smaller than~$1/2$. Interestingly, this maximal depth goes (monotonically) to zero as~$k\to\infty$. Note that, for the same distribution, location halfspace depth has, irrespective of~$k$, maximal value~$1/2$ (this follows, e.g., from Lemma~1 and Theorem~1 in \citealp{RouStr2004}).

In general, the halfspace deepest scatter~$\Sigma_*$ is not unique. This is typically the case for empirical probability measures~$P_n$ (note that the existence of a halfspace deepest scatter in the empirical case readily follows from the fact that $H\!D^{\rm sc}_{P_n,T}(\Sigma)$ takes its values in~$\{\ell/n:\ell=0,1,\ldots,n\}$). For several purposes, it is needed to identify a unique representative of the halfspace deepest scatters, that would play a similar role for scatter as the one played by the Tukey median for location. To this end, one may consider here a \emph{center of mass}, that is, a scatter matrix of the form 
\begin{equation}
	\label{centerofmass}
\Sigma_{P,T}
:=
\arg\min_{\Sigma\in \mathcal{P}_k} 
\int_{R^{\rm sc}_{P,T}(\alpha_{*P,T})} d_g^2(m,\Sigma) 
\, 
dm
,
\end{equation}
where~$dm$ is a mass distribution on~$R^{\rm sc}_{P,T}(\alpha_{*P,T})$ with total mass one (the natural choice being the uniform over~$R^{\rm sc}_{P,T}(\alpha_{*P,T})$). This is a suitable solution if~$R^{\rm sc}_{P,T}(\alpha_{*P,T})$ is $g$-bounded (hence, $g$-compact), since~\cite{Car29} showed that, in a simply connected manifold with non-positive curvature (as~$\mathcal{P}_k$), every compact set has a unique center of mass; see also Proposition~60 in~\cite{Ber03}. Convexity of~$R^{\rm sc}_{P,T}(\alpha_{*P,T})$ then ensures that~$\Sigma_{P,T}$ has maximal depth. Like for location, this choice of~$\Sigma_{P,T}$ as a representative of the deepest scatters guarantees affine equivariance (in the sense that~$\Sigma_{P_{A,b},T}=A\Sigma_{P,T}A'$ for any~$A\in GL_k$ and any~$b\in\R^k$), provided that~$T$ itself is affine-equivariant. An alternative approach is to consider the scatter matrix~$\Sigma_{P,T}$ whose vectorized form~${\rm vec}\,\Sigma_{P,T}$ is the barycenter of~${\rm vec}\,R^{\rm sc}_{P,T}(\alpha_{*P,T})$. While this is a more practical solution for scatter matrices, the non-flat nature of some of the parameter spaces in Section~\ref{Secshape} will require the more involved, manifold-type, approach in~(\ref{centerofmass}). 

As a final comment related to Theorem~\ref{geodmaxdepth}, note that if $R^{\rm sc}_{P,T}(\alpha_{P,T})$ is empty, then it may actually be so that no halfspace deepest scatter does exist. An example is provided by the bivariate mixture distribution~$P$ in Section~\ref{SecFrobenius}. There, we saw that, for any centro-equivariant~$T$, no halfspace deepest scatter does exist, which is compatible with the fact that, for any~$\Sigma$, $H\!D_{P,T}^{\rm sc}(\Sigma)<1-s=\alpha_{P,T}$, so that~$R^{\rm sc}_{P,T}(\alpha_{P,T})$ is empty.

\section{An axiomatic approach for scatter depth}
\label{SecAxiom}

Building on the properties derived in \cite{Liu1990} for simplicial depth,  \cite{ZuoSer2000A} introduced an axiomatic approach suggesting that a generic location depth $D^{\rm loc}_P(\,\cdot\,):\R^k\to [0,1]$ should satisfy the following properties: (P1) affine invariance, (P2) maximality at the symmetry center (if any), (P3) monotonicity relative to any deepest point, and (P4) vanishing at infinity. Without entering into details, these properties are to be understood as follows: (P1) means that~$D^{\rm loc}_{P_{A,b}}(A\theta+b)=D^{\rm loc}_{P}(\theta)$ for any~$A\in GL_k$ and~$b\in\R^k$, where~$P_{A,b}$ is as defined on page~\pageref{pagedefPAb}; (P2) states that if~$P$ is symmetric (in some sense), then the symmetry center should maximize~$D^{\rm loc}_P(\,\cdot\,)$; according to~(P3), $D^{\rm loc}_P(\,\cdot\,)$ should be monotone non-increasing along any halfline originating from any $P$-deepest point; finally, (P4) states that as~$\theta$ exits  any compact set in~$\R^k$, its depth should converge to zero. There is now an almost universal agreement in the literature that (P1)-(P4) are the natural desirable properties for location depths. 

In view of this, one may wonder what are the desirable properties for a scatter depth. Inspired by (P1)-(P4), we argue that a generic scatter depth~$D^{\rm sc}_P(\,\cdot\,):\mathcal{P}_k\to [0,1]$ should satisfy the following properties, all involving an (unless otherwise specified) arbitrary probability measure~$P$ over~$\R^k$:

\begin{enumerate}
\item[(Q1)] \emph{Affine invariance}: for any~$A\in GL_k$ and~$b\in\R^k$, 
$
D^{\rm sc}_{P_{A,b}}(A\Sigma A')
=
D^{\rm sc}_P(\Sigma)
,
$ 
where~$P_{A,b}$ is still as defined on page~\pageref{pagedefPAb};
\item[(Q2)] \emph{Fisher consistency under ellipticity}: if $P$ is elliptically symmetric with location~$\theta_0$ and scatter~$\Sigma_0$, then $D^{\rm sc}_P(\Sigma_0)\geq D^{\rm sc}_P(\Sigma)$ for any~$\Sigma\in\mathcal{P}_k$;
\item[(Q3)] \emph{Monotonicity relative to any deepest scatter}: if~$\Sigma_a$ maximizes~$D^{\rm sc}_P(\,\cdot\,)$, then, for any~$\Sigma_b\in\mathcal{P}_k$, $t\mapsto D^{\rm sc}_P((1-t)\Sigma_a+t\Sigma_b)$ is monotone non-increasing over~$[0,1]$;
\item[(Q4)] \emph{Vanishing at the boundary of the parameter space}: if~$(\Sigma_n)$ $F$-converges to the boundary of~$\mathcal{P}_k$ (in the sense that either~$d_F(\Sigma_n,\Sigma)\to 0$ for some~$\Sigma\in \mathcal{S}_k\setminus\mathcal{P}_k$ or~$d_F(\Sigma_n,I_k)\to\infty$), then $D^{\rm sc}_P(\Sigma_n)\to 0$.  
\end{enumerate}

While (Q1) and~(Q3) are the natural scatter counterparts of~(P1) and~(P3), respectively, some comments are in order for~(Q2) and~(Q4). 
We start with~(Q2). In essence, (P2) requires that, whenever an indisputable location center exists (as it is the case for symmetric distributions), this location should be flagged as most central by the location depth at hand. A similar reasoning leads to~(Q2): we argue that, for an elliptical probability measure, the ``true" value of the scatter parameter is indisputable, and (Q2) then imposes that the scatter depth at hand should identify this true scatter value as the (or at least, as a) deepest one. One might actually strengthen~(Q2) by replacing the elliptical model there by a broader model in which the true scatter would still be clearly defined. In such a case, of course, the larger the model for which scatter depth satisfies~(Q2), the better (a possibility, that we do not explore here, is to consider the union of the elliptical model and the \emph{independent component model}; see~\citealp{IlmPai2011} and the references therein). This is parallel to what happens in~(P2): the weaker the symmetry assumption under which~(P2) is satisfied, the better (for instance, having~(P2) satisfied with angular symmetry is better than having it satisfied with central symmetry only); see \cite{ZuoSer2000A}. 

We then turn to~(Q4), whose location counterpart~(P4) is typically read by saying that the depth/centrality~$D^{\rm loc}_P(\theta_n)$ goes to zero when the point~$\theta_n$ goes to the boundary of the \emph{sample space}. In the spirit of parametric depth (\citealp{Miz2002,MizMul2004}), however, it is more appropriate to look at~$\theta_n$ as a candidate location fit and to consider that~(P4) imposes that the appropriateness~$D^{\rm loc}_P(\theta_n)$ of this fit goes to zero as~$\theta_n$ goes to the boundary of the \emph{parameter space}. For location, the confounding between the sample space and parameter space (both are~$\R^k$) allows for both interpretations. For scatter, however, there is no such confounding (the sample space is~$\R^k$ and the parameter space is~$\mathcal{P}_k$), and we argue~(Q4) above is the natural scatter version of~(P4): whenever~$\Sigma_n$ goes to the boundary of the parameter space~$\mathcal{P}_k$, scatter depth should flag it as an arbitrarily poor candidate fit. 

Theorem~\ref{thaffinv} states that scatter halfspace depth satisfies~(Q1) as soon as it is based on an affine-equivariant~$T$. Scatter halfspace depth satisfies~(Q3) as well: if~$\Sigma_a$ maximizes~$H\!D^{\rm sc}_{P,T}(\,\cdot\,)$, then Theorem~\ref{Propquasiconcavity} indeed readily implies that~$H\!D^{\rm sc}_{P,T}((1-t)\Sigma_a+t\Sigma_b)\geq \min(H\!D^{\rm sc}_{P,T}(\Sigma_a),H\!D^{\rm sc}_{P,T}(\Sigma_b))=H\!D^{\rm sc}_{P,T}(\Sigma_b)$ for any~$\Sigma_b\in\mathcal{P}_k$ and~$t\in [0,1]$. The next Fisher consistency result shows that, provided that~$T$ is affine-equivariant,~(Q2) is also met.

\begin{Theor}
\label{Fishconst}
Let~$P$ be an elliptical probability measure over~$\R^k$ with location~$\theta_0$ and scatter~$\Sigma_0$, and let~$T$ be an affine-equivariant location functional. Then,
\vspace{-1.2mm}
  (i)  
$
H\!D^{\rm sc}_{P,T}(\Sigma_0)
\geq
H\!D_{P,T}^{\rm sc}(\Sigma)
$ for any~$\Sigma\in\mathcal{P}_k$, and the equality holds if and only
\vspace{-.4mm}
  if~${\rm Sp}(\Sigma_0^{-1}\Sigma)\subset \mathcal{I}_{\rm MSD}[Z_1]$, where~$Z=(Z_1,\ldots,Z_k)'\stackrel{\mathcal{D}}{=}\Sigma_0^{-1/2}(X-\theta_0)$; 
(ii) in particular, if~$\mathcal{I}_{\rm MSD}[Z_1]$ is a singleton $($equivalently, if~$\mathcal{I}_{\rm MSD}[Z_1]=\{1\})$, then  $
\Sigma\mapsto H\!D_{P,T}^{\rm sc}(\Sigma)$ is uniquely maximized at~$\Sigma_0$. 
\end{Theor}

While~(Q1)-(Q3) are satisfied by scatter halfspace depth without any assumption on~$P$,~(Q4) is not, as the mixture example considered in Section~\ref{SecFrobenius} shows (since the sequence~$(\Sigma_n)$ considered there has limiting depth~$1-s>0$). However, Theorem~\ref{boundedness} reveals that~(Q4) may fail only when~$d_F(\Sigma_n,\Sigma)\to 0$ for some~$\Sigma\in \mathcal{S}_k\setminus\mathcal{P}_k$. More importantly, Theorem~\ref{geodboundedness} implies that $T$-scatter halfspace depth will satisfy~(Q4) at any~$P$ that is smooth at~$T_P$. 

In a generic parametric depth setup,~(Q3) would require that the parameter space is convex. If the parameter space rather is a non-flat Riemannian manifold, then it is natural to replace the ``linear" monotonicity property~(Q3) with a ``geodesic" one. In the context of scatter depth, this would lead to replacing~(Q3) with 

\begin{enumerate}
\item[($\widetilde{\mathrm{Q}3}$)] \emph{Geodesic monotonicity relative to any deepest scatter}: if~$\Sigma_a$ maximizes~$D^{\rm sc}_P(\,\cdot\,)$, then, for any~$\Sigma_b\in\mathcal{P}_k$, $t\mapsto D^{\rm sc}_P(\tilde{\Sigma}_t)$ is monotone non-increasing over~$[0,1]$ along the geodesic path~$\tilde{\Sigma}_t$ from~$\Sigma_a$ to~$\Sigma_b$ in~(\ref{geod}).
\end{enumerate}
 
We refer to
\vspace{-.6mm}
  Section~\ref{Secshape} for a parametric framework where~(Q3) cannot be considered and where~($\widetilde{\mathrm{Q}3}$) needs to be adopted instead. For scatter, however, the hybrid nature of~$\mathcal{P}_k$, which is both flat (as a convex subset of the vector space~$\mathcal{S}_k$) and curved (as a Riemannian manifold with non-positive curvature), allows to consider both~(Q3) and~($\widetilde{\mathrm{Q}3}$). Just like~(Q3) follows from quasi-concavity of the mapping~$\Sigma\mapsto H\!D_{P,T}^{\rm sc}(\Sigma)$,~($\widetilde{\mathrm{Q}3}$) would follow from the same mapping being \emph{geodesic quasi-concave}, in the sense that $H\!D_{P,T}(\tilde{\Sigma}_t)\geq \min(H\!D_{P,T}(\Sigma_a),H\!D_{P,T}(\Sigma_b))$ along the geodesic path~$\tilde{\Sigma}_t$ from~$\Sigma_a$ to~$\Sigma_b$. Geodesic quasi-concavity would actually imply that scatter halfspace depth regions are \emph{geodesic convex}, in the sense that, for any~$\Sigma_a,\Sigma_b\in R^{\rm sc}_{P,T}(\alpha)$, the geodesic from~$\Sigma_a$ to~$\Sigma_b$ is contained in~$R^{\rm sc}_{P,T}(\alpha)$. We refer to \cite{DumTyl2016} for an application of geodesic convex functions to inference on (high-dimensional) scatter matrices. 

Theorem~\ref{Propquasiconcavity} shows that $\Sigma\mapsto H\!D_{P,T}^{\rm sc}(\Sigma)$ is quasi-concave for any~$P$. A natural question is then whether or not this extends to geodesic quasi-concavity. The answer is positive at any $k$-variate elliptical probability measure and at the $k$-variate probability measure with independent Cauchy marginals.

\begin{Theor}
\label{propgeodesicquasiconc}
Let~$P$ be an elliptical probability measure over~$\R^k$ or the $k$-variate probability measure with independent Cauchy marginals, and let $T$ be an affine-equivariant location functional. Then, (i) $\Sigma\mapsto H\!D_{P,T}^{\rm sc}(\Sigma)$ is geodesic quasi-concave, so that (ii) $R^{\rm sc}_{P,T}(\alpha)$ is geodesic convex for any~$\alpha\geq 0$. 
\end{Theor}

\begin{figure}[h!] 
\includegraphics[width=.8\textwidth]{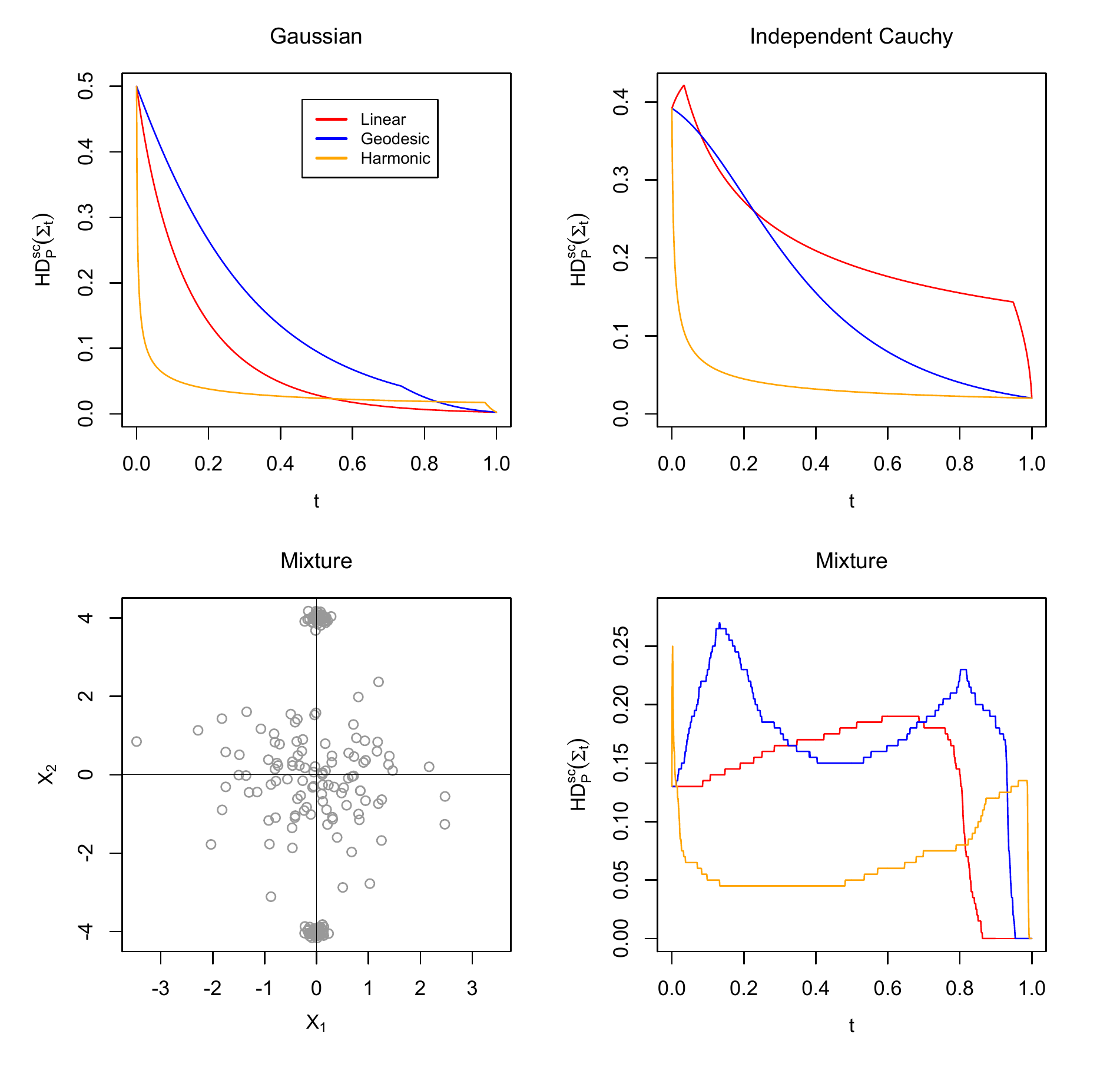}  
\vspace{-5mm}
 \caption{Plots, for various bivariate probability measures~$P$, of the scatter halfspace depth function~$\Sigma\mapsto H\!D_P^{\rm sc}(\Sigma)$ along the linear path~$\Sigma_t=(1-t)\Sigma_a+t\Sigma_b$ (red), the geodesic path~$\tilde{\Sigma}_t=\Sigma_a^{1/2}(\Sigma_a^{-1/2} \Sigma_b \Sigma_a^{-1/2})^t \Sigma_a^{1/2}$ (blue), and the harmonic path~$\Sigma_t^*=((1-t)\Sigma_a^{-1}+t\Sigma_b^{-1})^{-1}$ (orange), from~$\Sigma_a=I_2$ to~$\Sigma_b={\rm diag}(0.001,20)$; harmonic paths are introduced in Section~\ref{Secconcentration}. The probability measures considered are the bivariate normal with location zero and scatter~$I_2$ (top left), the bivariate distribution with independent Cauchy marginals (top right), and the empirical probability measure associated with a random sample of size~$n=200$ from the bivariate mixture distribution described in Section~\ref{SecAxiom} (bottom right). The scatter plot of the sample used in the mixture case is provided in the bottom left panel.}
\label{FigSec5}
\end{figure}

We close this section with a numerical illustration of the quasi-concavity results in Theorems~\ref{Propquasiconcavity} and~\ref{propgeodesicquasiconc} and with an example showing that geodesic quasi-concavity may actually fail to hold. Figure~\ref{FigSec5} provides, for three bivariate probability measures~$P$, the plots of~$t\mapsto H\!D^{\rm sc}_{P}(\Sigma_t)$ and~$t\mapsto H\!D^{\rm sc}_P(\tilde{\Sigma}_t)$, where~$\Sigma_t=(1-t)\Sigma_a+t\Sigma_b$ is the linear path from~$\Sigma_a=I_2$ to~$\Sigma_b={\rm diag}(0.001,20)$ and where~$\tilde{\Sigma}_t=\Sigma_a^{1/2}(\Sigma_a^{-1/2} \Sigma_b \Sigma_a^{-1/2})^t \Sigma_a^{1/2}$ is the corresponding geodesic path. The three distributions considered are (i) the bivariate normal with location zero and scatter~$I_2$, (ii) the bivariate distribution with independent Cauchy marginals, and (iii) the empirical distribution associated with a random sample of size~$n=200$ from the  bivariate mixture distribution~$P=\frac{1}{2}P_1+\frac{1}{4}P_2+\frac{1}{4}P_3$, where~$P_1$ is the standard normal, $P_2$ is the normal with mean~$(0,4)'$ and covariance matrix~$\frac{1}{10}I_2$, and~$P_3$ is the normal with mean~$(0,-4)'$ and covariance matrix~$\frac{1}{10}I_2$. Figure~\ref{FigSec5} illustrates that (linear) quasi-concavity of scatter halfspace depth always holds, but that geodesic quasi-concavity may fail to hold. Despite this counterexample, extensive numerical experiments led us to think that geodesic quasi-concavity is the rule rather than the exception.

\section{Concentration halfspace depth}
\label{Secconcentration}


In various setups, the parameter of interest is the \emph{concentration matrix}~$\Gamma:=\Sigma^{-1}$ rather than the scatter matrix~$\Sigma$. For instance, in Gaussian graphical models, the~$(i,j)$-entry of~$\Gamma$ is zero if and only if the~$i$th and $j$th marginals are conditionally independent given all other marginals. It may then be useful to define a depth for inverse scatter matrices. The scatter halfspace depth in~(\ref{def}) naturally leads to defining the \emph{$T$-concentration halfspace depth} of~$\Gamma$ with respect to~$P$ as
$$
H\!D^{\rm conc}_{P,T}(\Gamma):=H\!D^{\rm sc}_{P,T}(\Gamma^{-1})
$$ 
and the corresponding \emph{$T$-concentration halfspace depth regions} as
$
R^{\rm conc}_{P,T}(\alpha)
\linebreak
:=
\big\{\Gamma\in\mathcal{P}_k : H\!D^{\rm conc}_{P,T}(\Gamma)\geq \alpha \big\}
$,
$\alpha\geq 0
.
$ 
As indicated by an anonymous referee, the definition of $T$-concentration halfspace depth alternatively results, through the use of ``innovated transformation" (see, e.g., \citealp{HallJin2010}, \citealp{FJY2013}, or \citealp{FanLv2016}), from the concept of (an affine-invariant) $T$-scatter halfspace depth. 

Concentration halfspace depth and concentration halfspace depth regions inherit the properties of their scatter antecedents, sometimes with subtle modifications. The former is affine-invariant and the latter are affine-equivariant as soon as they are based on an affine-equivariant~$T$. Concentration halfspace depth is upper $F$- and $g$-semicontinuous for any probability measure~$P$ (so that the regions~$R^{\rm conc}_{P,T}(\alpha)$ are $F$- and $g$-closed) and $F$- and~$g$-continuous if~$P$ is smooth at~$T_P$. While the regions~$R^{\rm conc}_{P,T}(\alpha)$ are still $g$-bounded (hence also, $F$-bounded) for~$\alpha > \alpha_{P,T}$, the outer regions~$R^{\rm conc}_{P,T}(\alpha)$, $\alpha\leq \alpha_{P,T}$, here may fail to be $F$-bounded (this is because implosion of~$\Sigma$, under which scatter halfspace depth may fail to go below~$\alpha_{P,T}$, is associated with explosion of~$\Sigma^{-1}$). Finally, uniform consistency and existence of a concentration halfspace deepest matrix are guaranteed under the same conditions on~$P$ and~$T$ as for scatter halfspace depth. 

Quasi-concavity of concentration halfspace depth and convexity of the corresponding regions require more comments. The linear path~$t\mapsto (1-t)\Gamma_a+t \Gamma_b$ between the concentration matrices~$\Gamma_a=\Sigma_a^{-1}$ and~$\Gamma_b=\Sigma_b^{-1}$ determines a \emph{harmonic path}~$t\mapsto\Sigma_t^*:=((1-t)\Sigma_a^{-1}+t \Sigma_b^{-1})^{-1}$ between the corresponding scatter matrices~$\Sigma_a$ and~$\Sigma_b$. In line with the definitions adopted in the previous sections, we will say that~$f:\mathcal{P}_k\to\R$ is  \emph{harmonic quasi-concave} if~$f(\Sigma_t^*)\geq \min(f(\Sigma_a),f(\Sigma_b))$ for any~$\Sigma_a,\Sigma_b\in\mathcal{P}_k$ and~$t\in [0,1]$, and that a subset~$R$ of~$\mathcal{P}_k$ is \emph{harmonic convex} if~$\Sigma_a,\Sigma_b\in R$ implies that~$\Sigma_t^*\in R$ for any~$t\in [0,1]$. 
Clearly, concentration halfspace depth is quasi-concave 
if and only if scatter halfspace depth is harmonic quasi-concave, 
which turns out to be the case in the elliptical and independent Cauchy cases. We thus have the following result. 

\begin{Theor}
\label{propharmonicquasiconc}
Let~$P$ be an elliptical probability measure over~$\R^k$ or the $k$-variate probability measure with independent Cauchy marginals, and let $T$ be an affine-equivariant location functional. Then, (i) $\Gamma\mapsto H\!D_{P,T}^{\rm conc}(\Gamma)$ is quasi-concave, so that (ii) $R^{\rm conc}_{P,T}(\alpha)$ is convex for any~$\alpha\geq 0$. 
\end{Theor}

However, concentration halfspace depth may fail to be quasi-concave, since, as we show by considering the mixture example in Figure~\ref{FigSec5}, scatter halfspace depth may fail to be harmonic quasi-concave. The figure, that also plots scatter halfspace depth along harmonic paths, confirms that, while scatter halfspace depth is harmonic quasi-concave for the Gaussian and independent Cauchy examples there, it is not in the mixture example. In this mixture example, thus, concentration halfspace depth fails to be quasi-concave and the corresponding depth regions fail to be convex. This is not a problem per se --- recall that famous (location) depth functions, like, e.g., the simplicial depth from \cite{Liu1990}, may provide non-convex depth regions. 

For completeness, we present the following result which shows that some form of quasi-concavity for concentration halfspace depth survives.

\begin{Theor}
\label{Propharmonicquasiconcavity}
Let~$P$ be a probability measure over~$\R^k$ and $T$ be a location functional. Then, (i) $\Gamma\mapsto H\!D^{\rm conc}_{P,T}(\Gamma)$ is harmonic quasi-concave, so that (ii) $R^{\rm conc}_{P,T}(\alpha)$ is harmonic convex for any~$\alpha\geq 0$.
\end{Theor}

Since concentration halfspace depth is harmonic quasi-concave if and only if scatter halfspace depth is quasi-concave, the result is a direct corollary of Theorem~\ref{Propquasiconcavity}. Quasi-concavity and harmonic quasi-concavity clearly are dual concepts, relative to scatter and concentration halfspace depths (which justifies the~$^*$ notation in the path~$\Sigma_t^*$, dual to~$\Sigma_t$). Interestingly, $\Gamma\mapsto H\!D^{\rm conc}_{P,T}(\Gamma)$ is geodesic quasi-concave if and only if~$\Sigma\mapsto H\!D^{\rm sc}_{P,T}(\Sigma)$ is, so that concentration halfspace depth regions are geodesic convex if and only if scatter halfspace depth regions are.


\section{Shape halfspace depth}
\label{Secshape}


In many multivariate statistics problems (PCA, CCA, sphericity testing, etc.), it is sufficient to know the scatter matrix~$\Sigma$ up to a positive scalar factor. In PCA, for instance, all scatter matrices of the form~$c\Sigma$, $c>0$, indeed provide the same unit eigenvectors~$v_\ell(c\Sigma)$, $\ell=1,\ldots,k$, hence the same principal components. Moreover, when it comes to deciding how many principal components to work with, a common practice is to look at the proportions of explained variance~$\sum_{\ell=1}^m \lambda_\ell(c\Sigma)/\sum_{\ell=1}^k \lambda_\ell(c\Sigma)$, $m=1,\ldots,k-1$, which do not depend on~$c$ either. In PCA, thus, the parameter of interest is a \emph{shape matrix}, that is, a normalized version, $V$ say, of the scatter matrix~$\Sigma$.

The generic way to normalize a scatter matrix~$\Sigma$ into a shape matrix~$V$ is based on a scale functional~$S$, that is, on a mapping~$S:\mathcal{P}_k\to\R^+_0$ satisfying (i)~$S(I_k)=1$ and (ii)~$S(c\Sigma)=cS(\Sigma)$ for any~$c>0$ and~$\Sigma\in\mathcal{P}_k$. In this paper, we will further assume that~(iii) \label{pagescalefunctcond} if~$\Sigma_1,\Sigma_2\in\mathcal{P}_k$ satisfy~$\Sigma_2 \geq \Sigma_1$ (in the sense that~$\Sigma_2-\Sigma_1$ is positive semidefinite), then $S(\Sigma_2)\geq S(\Sigma_1)$. Such a scale functional leads to factorizing~$\Sigma(\in\mathcal{P}_k)$ into
$
\Sigma
=
\sigma^2_S V_S
,	
$
where~$\sigma^2_S:=S(\Sigma)$ is the \emph{scale} of~$\Sigma$ and~$V_S:=\Sigma/S(\Sigma)$ is its \emph{shape matrix} (in the sequel, we will drop the subscript~$S$ in~$V_S$ to avoid overloading the notation). The resulting collection of shape matrices~$\VS$ will be denoted as~$\mathcal{P}_k^{S}$. Note that the constraint~$S(I_k)=1$ ensures that, irrespective of the scale functional~$S$ adopted, $I_k$ is a shape matrix. Common scale functionals satisfying~(i)-(iii) are 
(a)
$S_{\rm tr}(\Sigma)=({\rm tr}\,\Sigma)/k$,
%
(b)
$S_{\rm det}(\Sigma)= (\det \Sigma)^{1/k}$,
%
(c)
$S_{\rm tr}^*(\Sigma)=k/({\rm tr}\,\Sigma^{-1})$, and
%
(d)
$S_{11}(\Sigma) = \Sigma_{11}$; 
we refer to \cite{PaiVanB2014} for references where the scale functionals~(a)-(d) are used.
The corresponding shape matrices~$\VS$ are then normalized in such a way that (a) ${\rm tr}[\VS]=k$, \mbox{(b) $\det \VS=1$,} (c) ${\rm tr}[\VS^{-1}]=k$, or (d)~$\VS_{11}=1$.  

In this section, we propose a concept of halfspace depth for shape matrices. More precisely, for a probability measure~$P$ over~$\R^k$, we define the \emph{$(S,T)$-shape halfspace depth} of~$\VS(\in \mathcal{P}_k^{S})$ with respect to~$P$ as
\begin{equation}
\label{shapedepthdef}
H\!D^{{\rm sh},S}_{P,T}(\VS)
:=
\sup_{\sigma^2>0}
H\!D^{\rm sc}_{P,T}(\sigma^2 \VS)
, 	
\end{equation}
where $H\!D^{\rm sc}_{P,T}(\sigma^2 \VS)$ is the $T$-scatter halfspace depth of~$\sigma^2 \VS$ with respect to~$P$. The corresponding depth regions are defined as
$$
R^{{\rm sh},S}_{P,T}(\alpha):=\{\VS\in\mathcal{P}_k^{S}:H\!D^{{\rm sh},S}_{P,T}(\VS)\geq \alpha\}
$$
(alike scatter, we will drop the index~$T$ in~$H\!D^{{\rm sh},S}_{P,T}(\VS)$ and~$R^{{\rm sh},S}_{P,T}(\alpha)$ whenever~$T$ is the Tukey median). The halfspace deepest shape (if any) is obtained by maximizing the ``profile depth" in~(\ref{shapedepthdef}), in the same way a profile likelihood approach would be based on the maximization of a (shape) profile likelihood of the form~$L^{\rm sh}_{\VS}=\sup_{\sigma^2>0} L_{\sigma^2 \VS}$. To the best of our knowledge, such a profile depth construction has never been considered in the literature.  

\begin{figure}[H] 
\includegraphics[width=.95\textwidth]{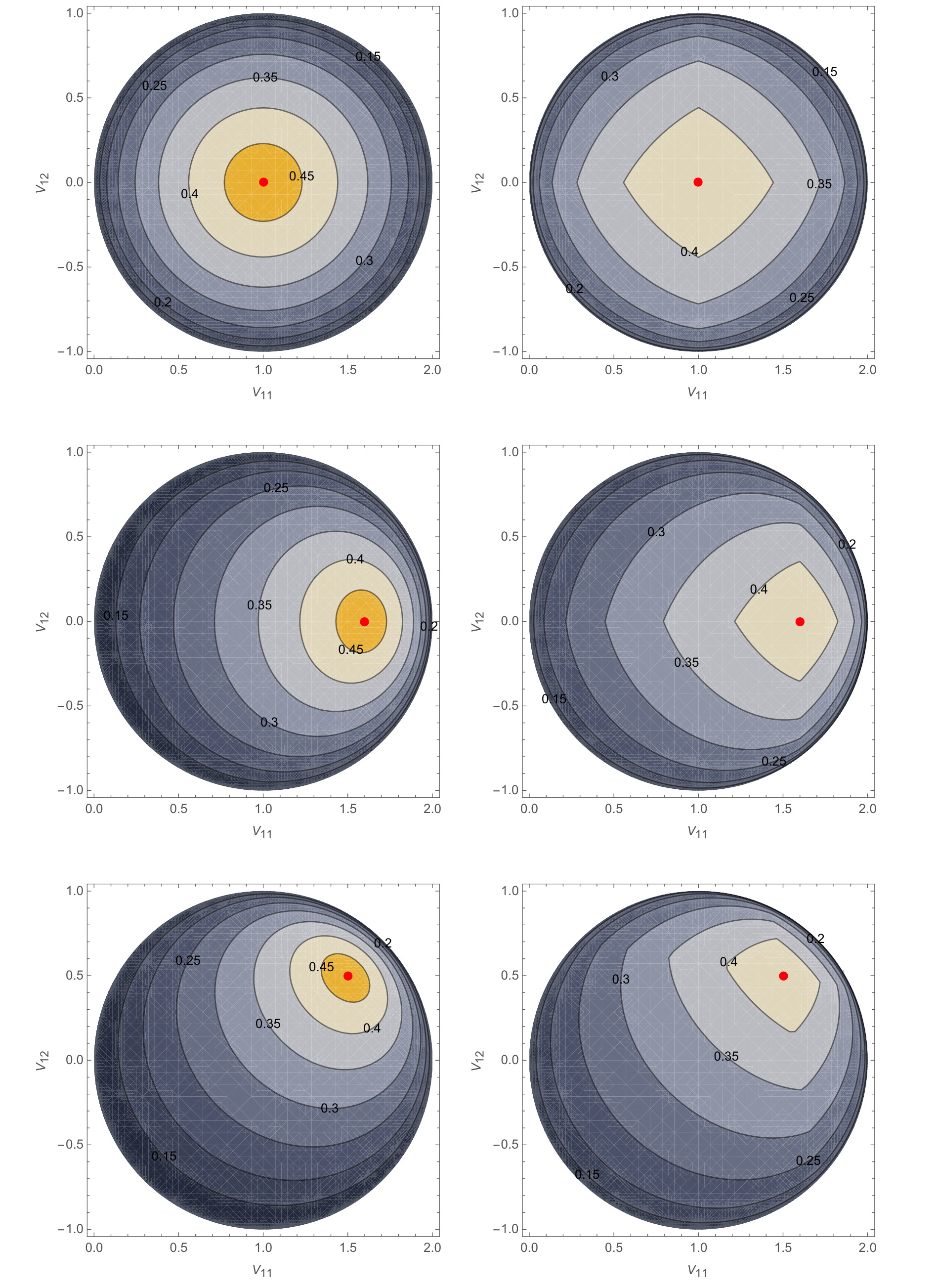}  
\caption{Contour plots of $(V_{11},V_{12})\mapsto H\!D^{{\rm sh},S_{\rm tr}}_{P,T}(V)$, for several bivariate probability measures~$P$ and an arbitrary affine-equivariant location functional~$T$, 
 where~$H\!D^{{\rm sh},S_{\rm tr}}_{P,T}(\VS)$ is the shape halfspace 
depth, with respect 
\vspace{-.7mm}
 to~$P$, of $V={\hspace{-1mm} V_{11}\ \hspace{1mm} V_{12}\, \choose \ V_{12}\ \hspace{0mm} \, 2-V_{11}\,}$.
\vspace{.1mm}
Letting $\Sigma_A={1\ 0\, \choose \,0\ 1\,}$, $\Sigma_B={4\ 0\, \choose \,0\ 1\,}$ and $\Sigma_C={3\ 1\, \choose \,1\ 1\,},$ 
\vspace{.3mm}
the probability measures~$P$ considered are those associated (i) with the bivariate normal distributions with location zero and scatter~$\Sigma_A$, $\Sigma_B$ and~$\Sigma_C$ (top, middle and bottom left), and (ii) with the distributions of~$\Sigma_A^{1/2}Z$, $\Sigma_B^{1/2}Z$ and~$\Sigma_C^{1/2}Z$, where~$Z$ has mutually independent Cauchy marginals (top, middle and bottom right). In each case, the ``true" $S_{\rm tr}$-shape matrix is marked in red.}
\label{FigSec7a} 
\end{figure}

We start the study of shape halfspace depth by considering our running, Gaussian and independent Cauchy, examples. For the $k$-variate normal with location~$\theta_0$ and scatter~$\Sigma_0$ (hence, with $S$-shape matrix~$V_{0}=\Sigma_0/S(\Sigma_0)$), 
$$
\sigma^2
\mapsto 
H\!D_{P,T}^{\rm sc}(\sigma^2 \VS)
=
2
\min\!\bigg(\!
\Phi\bigg(\frac{b\sigma\lambda^{1/2}_k(V_0^{-1}\VS)}{\sqrt{S(\Sigma_0)}}  \bigg) - {\frac{1}{2}}
,
1-\Phi\bigg(\frac{b\sigma\lambda^{1/2}_1(V_0^{-1}\VS)}{\sqrt{S(\Sigma_0)}}  \bigg)
\!\bigg)
$$ 
(see~(\ref{Gaussiandepth})) will be uniquely maximized at the~$\sigma^2$-value for which both arguments of the minimum are equal. It follows that
$$
H\!D^{{\rm sh},S}_{P,T}(\VS)
=
2\Phi\big(c(V_0^{-1}\VS) \lambda^{1/2}_k(V_0^{-1}\VS)\big) - 1
,
$$
where~$c(\Upsilon)$ is the unique solution of  
$
\Phi\big( c(\Upsilon) \lambda^{1/2}_k(\Upsilon)  \big) - {\textstyle\frac{1}{2}}
= 
1-\Phi\big( c(\Upsilon) \lambda^{1/2}_1(\Upsilon)  \big)
.
$
At the $k$-variate distribution with independent Cauchy marginals, we still have that (with the same notation as in~(\ref{Cauchytheo})) 
$$
H\!D^{\rm sc}_{P,T}(\sigma^2 \VS) 
=
2
\min
\Big(
\Psi\big(  \sigma/\max_s (s'\VS^{-1}s)^{1/2}\, \big) - 
{\textstyle\frac{1}{2}}
,
1 - \Psi\big( \sigma \sqrt{\max({\rm diag}(\VS))}\,  \big) 
\Big)
$$
is maximized for fixed~$V$ when both arguments of the minimum are equal, that is, when
$
\sigma^2
=
\big(
\max_s (s'\VS^{-1}s)/\max({\rm diag}(\VS))
\big)^{1/2}
.
$
Therefore, 
\begin{eqnarray*}
H\!D^{{\rm sh},S}_{P,T}(\VS)
&\!\!=\!\!&
2 \,\Psi\Big( \big( \max_s (s'V^{-1}s) \max({\rm diag}(V))\big)^{-1/4} \, \Big) - 1
\\[2mm]
&\!\!=\!\!&
\frac{2}{\pi}
\,
 \arctan \Big( \big( \max_s (s'\VS^{-1}s) \max({\rm diag}(\VS))\big)^{-1/4} \, \Big) 	
.
\end{eqnarray*}
Figure~\ref{FigSec7a} draws, for six probability measures~$P$ and any affine-equivariant~$T$, contour plots
\vspace{-.5mm}
  of~$(V_{11},V_{12})\mapsto H\!D^{{\rm sh},S_{\rm tr}}_{P,T}(V)$, where~$H\!D^{{\rm sh},S_{\rm tr}}_{P,T}(\VS)$ is the shape halfspace depth of $V={\hspace{-1mm} V_{11}\ \hspace{1mm} V_{12}\, \choose \ V_{12}\ \hspace{0mm} \, 2-V_{11}\,}$ with respect to~$P$. Letting
$\Sigma_A
={1\ 0\, \choose \,0\ 1\,}$,
$\Sigma_B
={4\ 0\, \choose \,0\ 1\,}$
and
$\Sigma_C
={3\ 1\, \choose \,1\ 1\,},$
the probability measures~$P$ considered are those associated (i) with the bivariate normal distributions with location zero and scatter~$\Sigma_A$, $\Sigma_B$ and~$\Sigma_C$, and (ii) with the distributions of~$\Sigma_A^{1/2}Z$, $\Sigma_B^{1/2}Z$ and~$\Sigma_C^{1/2}Z$, where~$Z$ has independent Cauchy marginals. Note that the maximal depth is larger in the Gaussian cases than in the Cauchy ones, that depth monotonically decreases along any ray originating from the deepest shape matrix and that it goes to zero if and only if the shape matrix converges to the boundary of the parameter space.  Shape halfspace depth contours are smooth in the Gaussian cases but not in the Cauchy ones.

In both the Gaussian and independent Cauchy examples above, the supremum in~(\ref{shapedepthdef}) is a maximum. For empirical probability measures~$P_n$, this will always be the case since~$H\!D^{\rm sc}_{P_n,T}(\sigma^2 V)$ then takes its values in~$\{\ell/n:\ell=0,1,\ldots,n\}$. The following result implies in particular that a sufficient condition for this supremum to be a maximum is that~$P$ is smooth at~$T_P$ (which is the case in both our running examples above). 

\begin{Theor}
\label{maxsigma2shape}
Let~$P$ be a probability measure over~$\R^k$ and~$T$ be a location functional. Fix~$\VS\in\mathcal{P}_k^{S}$ such that~$c \VS\in R^{\rm sc}_{P,T}(\alpha_{P,T})$ for some~$c>0$. Then, $H\!D^{{\rm sh},S}_{P,T}(\VS)
=H\!D^{\rm sc}_{P,T}(\sigma^2_{\VS} \VS)$ for some~$\sigma^2_{\VS}>0$. 
\end{Theor}

The following affine-invariance/equivariance and uniform consistency results are easily obtained from their scatter antecedents.

\begin{Theor}
\label{thaffinvshape}
Let~$T$ be an affine-equivariant location functional.  Then, 
(i) shape halfspace depth is affine-invariant in the sense that, for any probability measure~$P$ over~$\R^k$,  $V\in\mathcal{P}_k^{S}$, $A\in GL_k$ and~$b\in\R^k$, we have
$
H\!D^{{\rm sh},S}_{P_{A,b},T}(A V\! A' /
 S(A V\! A'))
=
H\!D^{{\rm sh},S}_{P,T}(V)
,
$
where~$P_{A,b}$ is as defined on page~\pageref{pagedefPAb}. Consequently, (ii) shape halfspace depth regions are affine-equivariant, in the sense that
$
R^{{\rm sh},S}_{P_{A,b},T}(\alpha)
=
\big\{
AV\! A'/
 S(A V\! A') : V\in R^{{\rm sh},S}_{P,T}(\alpha)
\big\}
$
for any probability measure~$P$ over~$\R^k$, $\alpha\geq 0$, $A\in GL_k$ and~$b\in\R^k$. 
\end{Theor}


\begin{Theor}
\label{Consistencyshape}
Let~$P$ be a smooth probability measure over~$\R^k$ and $T$ be a location functional. Let~$P_n$ denote the empirical probability measure associated with a random sample of size~$n$ from~$P$ and assume that~$T_{P_n}\to T_P$ almost surely as~$n\to\infty$. Then $\sup_{V\in\mathcal{P}_k^{S}}
|H\!D^{{\rm sh},S}_{P_n,T}(V)-H\!D^{{\rm sh},S}_{P,T}(V)| \to 0$ almost surely as~$n\to\infty$. 
\end{Theor}


Shape halfspace depth inherits the $F$- and $g$-continuity properties of scatter halfspace depth (Theorems~\ref{Propcontinuity} and~\ref{geodcontinuity}, respectively), at least for a smooth~$P$. More precisely, we have the following result. 

\begin{Theor}
\label{Propcontinuityshape}
Let~$P$ be a probability measure over~$\R^k$ and $T$ be a location functional. Then, 
\vspace{-.8mm}
 (i) $V\mapsto H\!D^{{\rm sh},S}_{P,T}(V)$ is upper $F$- and $g$-semicontinuous on~$R^{{\rm sh},S}_{P,T}(\alpha_{P,T})$, so that 
 \vspace{-.5mm}
 (ii) for any~$\alpha\geq \alpha_{P,T}$, the depth region~$R^{{\rm sh},S}_{P,T}(\alpha)$ is $F$- and $g$-closed. (iii) If~$P$ is smooth at~$T_P$, then~$V\mapsto H\!D^{{\rm sh},S}_{P,T}(V)$ is $F$- and $g$-continuous.  
\end{Theor}
 
The $g$-boundedness part of the following result will play a key role when proving the existence of a halfspace deepest shape. 

\begin{Theor}
\label{boundednessshape}
Let~$P$ be a probability measure over~$\R^k$ and $T$ be a location functional. Then, for any~$\alpha>\alpha_{P,T}$, $R^{{\rm sh},S}_{P,T}(\alpha)$ is $F$- and $g$-bounded, hence $g$-compact. If~$s_{P,T}\geq 1/2$, then this result is trivial in the sense that $R^{{\rm sh},S}_{P,T}(\alpha)$ is empty for~$\alpha>\alpha_{P,T}
$.  
\end{Theor}
 
Comparing with the scatter result in Theorem~\ref{boundedness}, the shape result for \mbox{$F$-boundedness} requires the additional condition~$\alpha>\alpha_{P,T}$ (for $g$-bounded\-ness, this condition was already required in Theorem~\ref{geodboundedness}). This condition is actually necessary for scale functionals~$S$ for which implosion of a shape matrix~$V$ cannot be obtained without explosion, as it is the case, e.g., for~$S_{\rm det}$ (the product of the eigenvalues of an $S_{\rm det}$-shape matrix being equal to one, the smallest eigenvalue of~$V$ cannot go to zero without the largest going to infinity). We illustrate this on the bivariate discrete example discussed below Theorem~\ref{geodboundedness}, still with an arbitrary centro-equivariant~$T$. The sequence of scatter matrices~$\Sigma_n={\rm diag}(\frac{1}{n},1)$ there defines a sequence of $S_{\rm det}$-shape matrices~$V_n={\rm diag}(\frac{1}{\sqrt{n}},\sqrt{n})$, 
\vspace{-1mm}
 that is neither $F$- nor $g$-bounded. Since~$H\!D^{{\rm sh},S_{\rm det}}_{P,T}(V_n)\geq H\!D^{\rm sc}_{P,T}(\Sigma_n)\geq 1/3=\alpha_{P,T}$ for any~$n$, we conclude that~$R^{{\rm sh},S_{\rm det}}_{P,T}(\alpha_{P,T})$ is both $F$- and $g$-unbounded. Note also that $F$-bounded\-ness of~$R^{{\rm sh},S}_{P,T}(\alpha)$ depends on~$S$. In particular, it is easy to check that the condition~$\alpha>\alpha_{P,T}$ for~$F$-boundedness is not needed for the scale functional~$S_{\rm tr}^*$ (that is, $R_{P,T}^{{\rm sh},S_{\rm tr}^*}(\alpha)$ is $F$-bounded for any~$\alpha>0$). Finally, one trivially has that all $R_{P,T}^{{\rm sh},S_{\rm tr}}(\alpha)$'s are $F$-bounded since the corresponding collection of shape matrices, $\mathcal{P}_k^{S_{\rm tr}}$, itself is $F$-bounded. Unlike $F$-boundedness, $g$-boundedness results are homogeneous in~$S$, which further suggests that the $g$-topology is the most appropriate one to study scatter/shape depths.

As announced, the $g$-part of Theorem~\ref{boundednessshape} allows to show that a halfspace deepest shape exists under mild conditions. More precisely, we have the following result.

\begin{Theor}
\label{geodmaxdepthshape}
Let~$P$ be a probability measure over~$\R^k$ and $T$ be a location functional. Assume 
\vspace{-.8mm}
 that~$R^{{\rm sh},S}_{P,T}(\alpha_{P,T})$ is non-empty. Then, $\alpha^S_{*P,T}:=\sup_{V\in\mathcal{P}_k^{S}} H\!D^{{\rm sh},S}_{P,T}(V)=H\!D^{{\rm sh},S}_{P,T}(V_*)$ for some~$V_*\in\mathcal{P}_k^{S}$. 
\end{Theor}

Alike scatter, a sufficient condition for the existence of a halfspace deepest shape is thus that~$P$ is smooth at~$T_P$.
In particular, a halfspace deepest shape exists in the Gaussian and independent Cauchy examples. In the $k$-variate independent Cauchy case, it readily follows from Theorem~\ref{propmaxCauchy} that, irrespective of the centro-equivariant~$T$ 
\vspace{-.8mm}
used,~$H\!D^{{\rm sh},S}_{P,T}(V)$ is uniquely maximized at~$V_*=I_k$, with corresponding maximal depth~$\frac{2}{\pi} \arctan\big( k^{-1/4} \big)$. The next Fisher-consistency result states that, in the elliptical case, the halfspace deepest shape coincides with the ``true" shape matrix. 

\begin{Theor}
\label{Fishconstshape}
Let~$P$ be an elliptical probability measure over~$\R^k$ with location~$\theta_0$ and scatter~$\Sigma_0$, hence with $S$-shape matrix~$V_0=\Sigma_0/S(\Sigma_0)$, and let~$T$ be an affine-equivariant location functional. Then, 
\vspace{-.7mm}
(i) 
$
H\!D^{{\rm sh},S}_{P,T}(V_0)
\geq
H\!D^{{\rm sh},S}_{P,T}(V)
$ for any~$V\in\mathcal{P}_k^{S}$;
(ii) if~$\mathcal{I}_{\rm MSD}[Z_1]$ is
\vspace{-1.4mm}
 a singleton 
 $($equivalently, if~$\mathcal{I}_{\rm MSD}[Z_1]=\{1\})$, where~$Z=(Z_1,\ldots,Z_k)'\stackrel{\mathcal{D}}{=}\Sigma_0^{-1/2}(X-\theta_0)$, then  $
V\mapsto H\!D^{{\rm sh},S}_{P,T}(V)$ is uniquely maximized at~$V_0$. 
\end{Theor}

We conclude this section by considering quasi-concavity properties of shape halfspace depth and convexity properties of the corresponding depth regions. It should be 
 noted that, for some scale functionals~$S$, the collection~$\mathcal{P}_k^{S}$ of $S$-shape matrices is not convex; for instance,
 \vspace{-1mm}
 neither~$\mathcal{P}_k^{S_{\det}}$ nor~$\mathcal{P}_k^{S_{\rm tr}^*}$ is convex, so that it does not make sense 
 \vspace{-.8mm}
to investigate whether or not~$V\mapsto H\!D^{{\rm sh},S}_{P,T}(V)$ is quasi-concave for these scale functionals. It does, however, for~$S_{\rm tr}$ and~$S_{11}$, and we have the following result.

\begin{Theor}
\label{Propquasiconcavityshape}
Let~$P$ be a probability measure over~$\R^k$ and $T$ be a location functional. Fix~$S=S_{\rm tr}$ or~$S=S_{11}$. Then, (i) $V\mapsto H\!D^{{\rm sh},S}_{P,T}(V)$ is quasi-concave, that is, for any~$V_a,V_b\in\mathcal{P}_k^{S}$ and~$t\in [0,1]$, $H\!D^{{\rm sh},S}_{P,T}(V_t)\geq \min(H\!D^{{\rm sh},S}_{P,T}(V_a),H\!D^{{\rm sh},S}_{P,T}(V_b))$, where we let~$V_t:=(1-t)V_a+t V_b$; (ii) for any~$\alpha\geq 0$, $R^{{\rm sh},S}_{P,T}(\alpha)$ is convex.
\end{Theor}

As mentioned above, neither~$\mathcal{P}_k^{S_{\det}}$ nor~$\mathcal{P}_k^{S_{\rm tr}^*}$ are convex in the usual sense (unlike for~$S_{\rm tr}$ and~$S_{11}$, thus, a unique halfspace deepest shape could not be defined through barycenters but would rather require a center-of-mass approach as in~(\ref{centerofmass})). However, $\mathcal{P}_k^{S_{\det}}$ is geodesic convex, which justifies studying the possible geodesic convexity of~$R^{\rm sh}_{P,S_{\det}}(\alpha)$ (this provides a parametric framework for which the shape version of~(Q3) in Section~\ref{SecAxiom} cannot be considered and for which it is needed to adopt the corresponding Property~($\widetilde{\mathrm{Q}3}$) instead). Similarly, $\mathcal{P}_k^{S_{\rm tr}^*}$ is harmonic
  \vspace{-.7mm}
 convex, so that it makes sense to investigate the harmonic convexity of~$R^{{\rm sh},S_{\rm tr}^*}_{P,T}(\alpha)$. We have the following results.

\begin{Theor}
\label{propgeodesicquasiconcshape}
Let $T$ be an affine-equivariant location functional and~$P$ be an arbitrary probability measure over~$\R^k$ for which $T$-scatter halfspace depth is geodesic quasi-concave.
\vspace{-1mm}
Then, (i) $V\mapsto H\!D^{{\rm sh},S_{\det}}_{P,T}(V)$ is geodesic 
 quasi-concave, so that~(ii) $R^{{\rm sh},S_{\det}}_{P,T}(\alpha)$ is geodesic convex for any~$\alpha\geq 0$.
\end{Theor}

\begin{Theor}
\label{Propquasiconcavityharmonicshape}
Let $T$ be an affine-equivariant location functional and~$P$ be an arbitrary probability measure over~$\R^k$ for which $T$-scatter halfspace depth is harmonic quasi-concave.
\vspace{-1mm}
Then, (i) $V\mapsto H\!D^{{\rm sh},S_{\rm tr}^*}_{P,T}(V)$ is harmonic
 quasi-concave, so that~(ii) $R^{{\rm sh},S_{\rm tr}^*}_{P,T}(\alpha)$ is harmonic convex for any~$\alpha\geq 0$.
\end{Theor}

An illustration of Theorems~\ref{Propquasiconcavityshape}-\ref{Propquasiconcavityharmonicshape} is provided in Appendix~\ref{SecSupShape}.


\section{A real-data application}
\label{Secrealdata}

In this section, we analyze the returns of the Nasdaq Composite and S$\&$P500 indices from February 1st, 2015 to February 1st, 2017. 
During that period, for each trading day and for each index, we collected returns every 5 minutes (that is, the difference between the index at a given time and $5$ minutes earlier, when available), resulting in usually $78$ bivariate observations per day. 
Due to some missing values, the exact number of returns per day varies, and only days with at least 70 observations were considered.  
The resulting dataset comprises a total of $38489$ bivariate returns distributed over $D=478$ trading days.  

The goal of this analysis is to determine which days, during the two-year period, exhibit an atypical behavior. In line with the fact that the main focus in finance is on volatily, atypicality here will refer to deviations from the ``global" behavior either in scatter (i.e., returns do not follow the global dispersion pattern) or in scale only (i.e., returns show a usual shape but their overall size is different). Atypical days will be detected by comparing intraday estimates of scatter and shape with a global version. 

Below, $\hat \Sigma_{\rm full}$ will denote the minimum covariance determinant (MCD) scatter estimate on the empirical distribution~$P_{\rm full}$ of the returns over the two-year period, and~$\hat V_{\rm full}$ will stand for the resulting shape estimate~$\hat V_{\rm full}=\hat\Sigma_{\rm full}/S_{\det}(\hat\Sigma_{\rm full})$. For any~$d=1,\ldots,D$, $\hat \Sigma_d$ and~$\hat V_d$ will denote the corresponding estimates on the empirical distribution~$P_d$ on day~$d$.

The rationale behind the choice of MCD rather than standard covariance as an estimation method for scatter/shape is twofold. 
First, the former will naturally deal with outliers inherently arising in the data (the first few returns after an overnight or weekend break are famously more volatile and their importance should be downweighted in the estimation procedure). 
Second, as hinted above, the global estimate will provide a baseline to measure the atypicality of any given day, which will be done, among others, using its intraday depth. It would be natural to use halfspace deepest scatter/shape matrices on~$P_{\rm full}$ as global estimates for scatter/shape. While locating the exact maxima is a  non-trivial task, the MCD shape estimator has already a high depth value ($H\!D^{{\rm sh},S_{\det}}_{P_{\rm full}}(\hat V_{\rm full})=0.481$), which makes it a very good proxy for the halfspace deepest shape. 
For the same reason, the scaled MCD estimator $\bar\Sigma_{\rm full}=\sigma^2_{\rm full} V_{\rm full}\textrm{ with }\sigma^2_{\rm full}=\textrm{argmax}_{\sigma^2} H\!D^{\rm sc}_{P_{\rm full}}(\sigma^2 V_{\rm full})$ (that, obviously, satisfies $H\!D^{\rm sc}_{P_{\rm full}}(\bar\Sigma_{\rm full})=0.481$) is similarly an excellent proxy for the halfspace deepest scatter. 
In contrast, the shape estimate associated with the standard covariance matrix (resp., the deepest scaled version of the covariance matrix) has a global shape (resp., scatter) depth of only $0.426$. 

For each day, the following measures of (a)typicality (three for scatter, three for shape) are
\vspace{-1mm}
  computed: (i) the scatter depth $H\!D^{\rm sc}_{P_d}(\bar\Sigma_{\rm full})$ of $\bar\Sigma_{\rm full}$ in day $d$, (ii) the 
  \vspace{-1mm}
shape depth $H\!D^{{\rm sh},S_{\det}}_{P_d}(\hat V_{\rm full})$ of $\hat V_{\rm full}$ in day $d$, (iii) the scatter Frobenius distance $d_F(\hat \Sigma_d,\hat\Sigma_{\rm full})$, (iv) the shape Frobenius distance $d_F(\hat V_d,\hat V_{\rm full})$, (v) the scatter geodesic distance $d_g(\hat \Sigma_d,\hat\Sigma_{\rm full})$, and (vi) the shape geodesic distance $d_g(\hat V_d,\hat V_{\rm full})$.
Of course, low depths or high distances point to atypical days. Practitioners might be tempted to base the distances in (iii)-(vi) on standard covariance estimates, which would actually provide poorer performances in the present outlier detection exercise (due to the masking effect resulting from using a non-robust global dispersion measure as a baseline). Here, we rather use  MCD-based estimates to ensure a fair comparison with the depth-based methods in (i)-(ii). 
 
\begin{figure}[h!] 
\includegraphics[width=\textwidth]{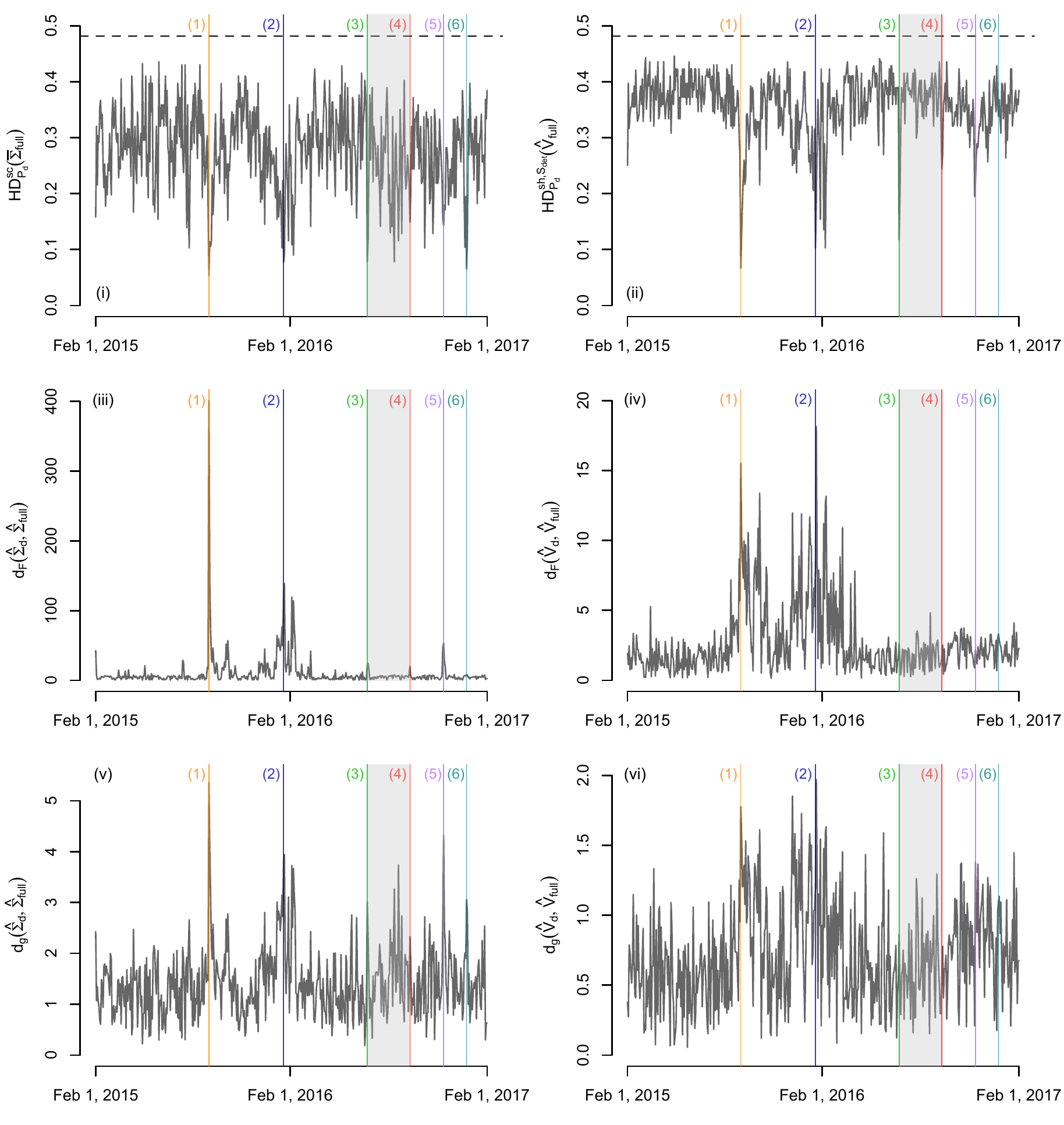}  
 \caption{Plots of 
 \vspace{-.7mm}
 (i) $H\!D^{\rm sc}_{P_{d}}(\bar \Sigma_{\rm full})$, (ii) $H\!D^{{\rm sh},S_{\det}}_{P_d}(\hat V_{\rm full})$, (iii) $d_F(\hat \Sigma_d,\hat\Sigma_{\rm full})$, (iv) $d_F(\hat V_d,\hat V_{\rm full})$, (v) $d_g(\hat \Sigma_d,\hat\Sigma_{\rm full})$ and (vi) $d_g(\hat V_d,\hat V_{\rm full})$, as a function of $d$, for the MCD scatter and shape estimates described in Section~\ref{Secrealdata}. The horizontal dotted lines in~(i)-(ii) correspond to the global depths $H\!D^{\rm sc}_{P_{\rm full}}(\bar\Sigma_{\rm full})$ and $H\!D^{{\rm sh},S_{\det}}_{P_{\rm full}}(\hat V_{\rm full})$, respectively. All depths make use of the Tukey median as a location functional. Vertical lines mark the six events listed in Section~\ref{Secrealdata}.}
\label{FigSec8}
\end{figure}

Figure~\ref{FigSec8} provides the plots of the quantities in~(i)-(vi) above as a function of $d$, $d=1,\dots,D$. Major events affecting the returns during the two years are marked there. They are (1) the Black Monday on August 24th, 2015 (orange) when world stock markets went down substantially, (2) the crude oil crisis on January 20th, 2016 (dark blue) when oil barrel prices fell sharply, (3) the Brexit vote aftermath on June 24th, 2016 (dark green), (4) the end of the low volatility period on September 13th, 2016 (red), (5) the Donald Trump election on November 9th, 2016 (purple), and (6) the announcement and aftermath of the federal rate hikes on December 14th, 2016 (teal). 

Detecting atypical events was achieved by flagging outliers in either collections of scatter or shape depth values. This was conducted by constructing box-and-whisker plots of those collections and marking events with depth value below 1.5 IQR of the first quartile. This procedure flagged events~(1), (2) and (6) as outlying in scatter and 21 days --- including events (1), (2), (3) and (5) --- as atypical in shape. Most of the resulting 22 outlying days can be associated (that is, are temporally close) to one of the \mbox{events~(1)-(6)} above. For example, 9 days are flagged within the period extending from January 20th, 2016 to February 9th, 2016, during which continuous slump in oil prices rocked the marked strongly, with biggest loss for S$\&$P 500 index on February 9th. Remarkably, out of the 22 flagged outliers, only two (namely October 1st, 2015 and December 14th, 2015) could not be associated with major events. Event (4), although not deemed outlying, was added to mark the end of the low volatility period.

Events (1) and (2) are noticeably singled out by all outlyingness measures, displaying low depth values and high Frobenius and geodesic distances, but the four remaining events tell a very different story. In particular, event~(6) exhibits a low scatter depth but a relatively high shape depth, which means that this day shows a shape pattern that is in line with the global one but is very atypical in scale (that is, in volatility size). Quite remarkably, the four distances considered fail to flag this day as an atypical one. A similar behavior appears throughout the two-month period spanning July, August and early September 2016 (between events~(3) and~(4)), during which the markets have seen a historical streak of small volatility. This period presents widely varying scatter depth values together with stable and high shape depth values, which is perfectly in line with what has been seen on the markets, where only the volatility of the indices was low in days that were otherwise typical. Again, the four distance plots are blind to this relative behavior of scatter and shape in the period. 
 
 Events (3) to (5) are picked up by depth measures and scatter distances, though more markedly by the former. This is particularly so for event~(3), which sticks out sharply in both depths. The fact that the scatter depth is even lower than the shape depth suggests that event~(3) is atypical not only in shape but also in scale. Interestingly, distance measures fully miss the shape outlyingness of this event. Actually, shape distances do not assign large values to any of the events~(3) to~(6) and, from March 2016 onwards, these distances stay in the same range --- particularly so for the Frobenius ones in~(iv). In contrast, the better ability of shape depth to spot outlyingness may be of particular importance in cases where one wants to discard the overall  volatility size to rather focus on the shape structure of the returns. 
 
To summarize, the detection of atypical patterns in the dispersion of intraday returns can more efficiently be performed with scatter/shape depths than on the basis of distance measures. Arguably, the fact that the proposed depths use all observations and not a sole estimate of scatter/shape allows to detect deviations from global behaviors more sharply. As showed above, comparing scatter and shape depth values provides a tool that permits the distinction between shape and scale outliers. 

\section{Final comments and perspectives}
\label{Secfinal}


In this work, we thoroughly investigated the structural properties of a concept of scatter halfspace depth linked to those proposed in~\cite{Zha2002} and~\cite{Chenetal16}. While we tried doing so under minimal assumptions, alternative scatter halfspace depth concepts may actually require even weaker assumptions, but they typically would make the computational burden heavier in the sample case. As an example, one might alternatively define the scatter halfspace depth of~$\Sigma(\in\mathcal{P}_k)$ with respect to~$P$ as
\begin{equation}
\label{scalt}	
H\!D^{\rm sc,alt}_P(\Sigma)
=
\sup_{\theta\in\R^k}
H\!D^{{\rm sc}}_{P,\theta}(\Sigma)
,
\end{equation}
where~$H\!D^{{\rm sc}}_{P,\theta}(\Sigma)$ is the scatter halfspace depth associated with the constant location functional at~$\theta$. This alternative scatter depth concept satisfies a uniform consistency result such as the one in Theorem~\ref{Consistency} without any condition on~$P$, whereas the scatter halfspace depth~$H\!D^{\rm sc}_{P,T}(\Sigma)$ in~(\ref{def}) 
 requires that~$P$ is smooth (see Theorem~\ref{Consistency}). In the sample case, however, evaluation of~$H\!D^{\rm sc,alt}_{P_n}(\Sigma)$ is computationally much more involved than~$H\!D^{\rm sc}_{P_n,T}(\Sigma)$. Alternative concentration and shape halfspace depth concepts may be defined along the same lines and will show the same advantages/disadvantages. compared to those proposed in this paper.      

Another possible concept of scatter halfspace depth bypasses the need to choose a location functional~$T$ by exploiting a pairwise difference approach; see \cite{Zha2002} and \cite{Chenetal16}. In our notation, the resulting scatter depth of~$\Sigma$ with respect to~$P=P^X$ is
\begin{equation}
\label{scalt2}	
H\!D^{\rm sc,U}_P(\Sigma)
=
H\!D^{{\rm sc}}_{P^{X-\tilde{X}},0}(\Sigma)
,
\end{equation}
where~$\tilde{X}$ is an independent copy of~$X$ and where $0$ denotes the origin of~$\R^k$. On one hand, the sample version of~(\ref{scalt2}) is a $U$-statistic of order two, which will increase the computational burden compared to the sample version of~(\ref{def}). On the other hand, uniform consistency results for~(\ref{scalt2}) (which here follow from Glivenko-Cantelli results for $U$-processes, such as the one in Corollary~3.3 from \citealp{ArcGin1993}) will again hold without any assumption on~$P$, which is due to the fact that, as already mentioned, the smoothness assumption in Theorem~\ref{Consistency} is superfluous when a constant location functional~$T$ is used. At first sight, thus, the pros and cons for~(\ref{scalt2}) are parallel to those for~(\ref{scalt}), that is, weaker distributional assumptions are obtained at the expense of computational ease. However, (\ref{scalt2}) suffers from a major disadvantage: it does not provide Fisher consistency at the elliptical model (see~(Q2) in Section~\ref{SecAxiom}). This results 
\vspace{-.4mm}
from the fact that if~$P=P^X$ is elliptical with location~$\theta$ and scatter~$\Sigma$, then~$P^{X-\tilde{X}}$ is elliptical with location~$0$ and scatter~$c_P\Sigma$, where the scalar factor~$c_P$ depends on the type of elliptical distribution: for multinormal and Cauchy elliptical distributions, e.g., $c_P=2$ and~$4$, respectively, so that if one replaces~$X-\tilde{X}$ with~$(X-\tilde{X})/\sqrt{2}$ to achieve Fisher consistency at the multinormal, then Fisher consistency will still not hold at the Cauchy. Actually, the maximizer of~$H\!D^{\rm sc,U}_P(\Sigma)$ is useless as a measure of scatter for the original probability measure~$P$, as its interpretation requires knowing which type of elliptical distribution~$P$ is. This disqualifies the pairwise difference scatter depth, as well as the companion concentration depth. Note, however, that the corresponding shape depth will not suffer from this Fisher consistency problem since the normalization of scatter matrices into shape matrices will get rid of the scalar factor~$c_P$.  

As both previous paragraphs suggest and as it is often the case with statistical depth, computational aspects are key for the application of the proposed depths. Evaluating (good approximations of) the scatter halfspace depth~$H\!D^{\rm sc}_{P_n,T}(\Sigma)$ of a 
\vspace{-.4mm}
 given~$\Sigma$ can of course be done for very small dimensions~$k=2$ or~$3$ by simply sampling the unit sphere~$\mathcal{S}^{k-1}$.
Even for such small dimensions, however, computing the halfspace deepest scatter is non-trivial: while scatter halfspace depth relies on a low-dimensional (that is, $k$-dimensional) projection-pursuit approach, identifying the halfspace deepest scatter indeed requires exploring the collection of scatter matrices~$\mathcal{P}_k$, that is of higher dimension, namely of dimension~$k(k+1)/2$. Fortunately, the fixed-location scatter halfspace depth --- hence, also its $T$-version proposed in this paper, after appropriate centering of the observations --- can be computed in higher dimensions through the algorithm proposed in~\cite{Chenetal16}, where the authors performed simulations requiring to compute the deepest scatter matrix for dimensions and sample sizes as large as~$10$ and~$2000$, respectively. Their  implementation of this algorithm is available as an R package at \url{https://github.com/ChenMengjie/DepthDescent}.

As pointed by an anonymous referee, the concept of scatter halfspace depth also makes sense when the parameter space is the compactification of~$\mathcal{P}_k$, that is, is the collection~$\overline{\mathcal{P}}_k$ of $k\times k$ symmetric positive \emph{semi-}definite matrices. Interestingly, it is actually easier to investigate the properties of scatter halfspace depth over~$\overline{\mathcal{P}}_k$ than over~$\mathcal{P}_k$. The $F$-continuity and \mbox{$F$-boundedness} results in Theorems~3.1-3.2 extend, mutatis mutandis, to~$\overline{\mathcal{P}}_k$. Unlike~$(\mathcal{P}_k,d_F)$, the metric space~$(\overline{\mathcal{P}}_k,d_F)$ is complete, so that the regions~$R^{\rm sc}_{P,T}(\alpha)$ are then $F$-compact for any~$\alpha>0$. Consequently, a trivial adaptation of the proof of Theorem~4.3 allows to show that there \emph{always} exists a halfspace deepest scatter matrix in~$\overline{\mathcal{P}}_k$. It is fortunate that these neat results can be established by considering the $F$-distance only, as the geodesic distance, that is unbounded on~$\overline{\mathcal{P}}_k\times\overline{\mathcal{P}}_k$, could not have been considered here. Of course, in many applications,~$\mathcal{P}_k$ remains the natural parameter space since many multivariate statistics procedures will require inverting scatter matrices. In such applications, it will be of little help to practitioners that the deepest halfspace scatter matrix belongs to~$\overline{\mathcal{P}}_k\setminus \mathcal{P}_k$, which explains why our detailed investigation focusing on~$\mathcal{P}_k$ is of key importance.

Perspectives for future research are rich and diverse. The proposed halfspace depth concepts for scatter, concentration and shape can be extended to other scatter functionals of interest. In particular, halfspace depths that are relevant for PCA
could result from the ``profile depth" approach in Section~\ref{Secshape}.
For instance, the \emph{\!$\,T$-``first principal direction" halfspace depth} of~$\beta(\in\mathcal{S}^{k-1})$  with respect to the probability measure~$P$ over~$\R^k$ can be defined as 
$$
H\!D^{1^{\rm st}{\rm pd}}_{P,T}(\beta)
=
\sup_{\Sigma\in \mathcal{P}_{k,1,\beta}}
H\!D^{{\rm sc}}_{P,T}(\Sigma)
,
\ 
\textrm{ with}
\
\mathcal{P}_{k,1,\beta}
:=
\{\Sigma\in\mathcal{P}_k:\Sigma\beta=\lambda_1(\Sigma)\beta\}
.
$$
The halfspace deepest first principal direction is a promising robust estimator of the true underlying first principal direction, at least under ellipticity. Obviously, the depth of any other principal direction, or the depth of any eigenvalue, can be defined accordingly. Another direction of research is to explore inferential applications of the proposed depths. Clearly, point estimation is to be based on halfspace deepest scatter, concentration or shape matrices; \cite{Chenetal16} partly studied this already for scatter in high dimensions. Hypothesis testing is also of primary interest. In particular, a natural test for~$\mathcal{H}_0:\Sigma=\Sigma_0$, where~$\Sigma_0\in\mathcal{P}_k$ is fixed, would reject the null for small values of~$H\!D^{{\rm sc}}_{P_n,T}(\Sigma_0)$. For shape matrices, 
\vspace{-1mm}
 a test of sphericity would similarly reject the null for small values of~$H\!D^{{\rm sh},S}_{P_n,T}(I_k)$. These topics, obviously, are beyond the scope of the present work.


\appendix



\section{Numerical illustrations}
\label{SecSupNum}


\subsection{Validation of the scatter halfspace depth expressions in~(\ref{Gaussiandepth})-(\ref{Cauchytheo})}
\label{SecSupMonteCarlo}

\begin{figure}[h!]
\includegraphics[width=\textwidth]{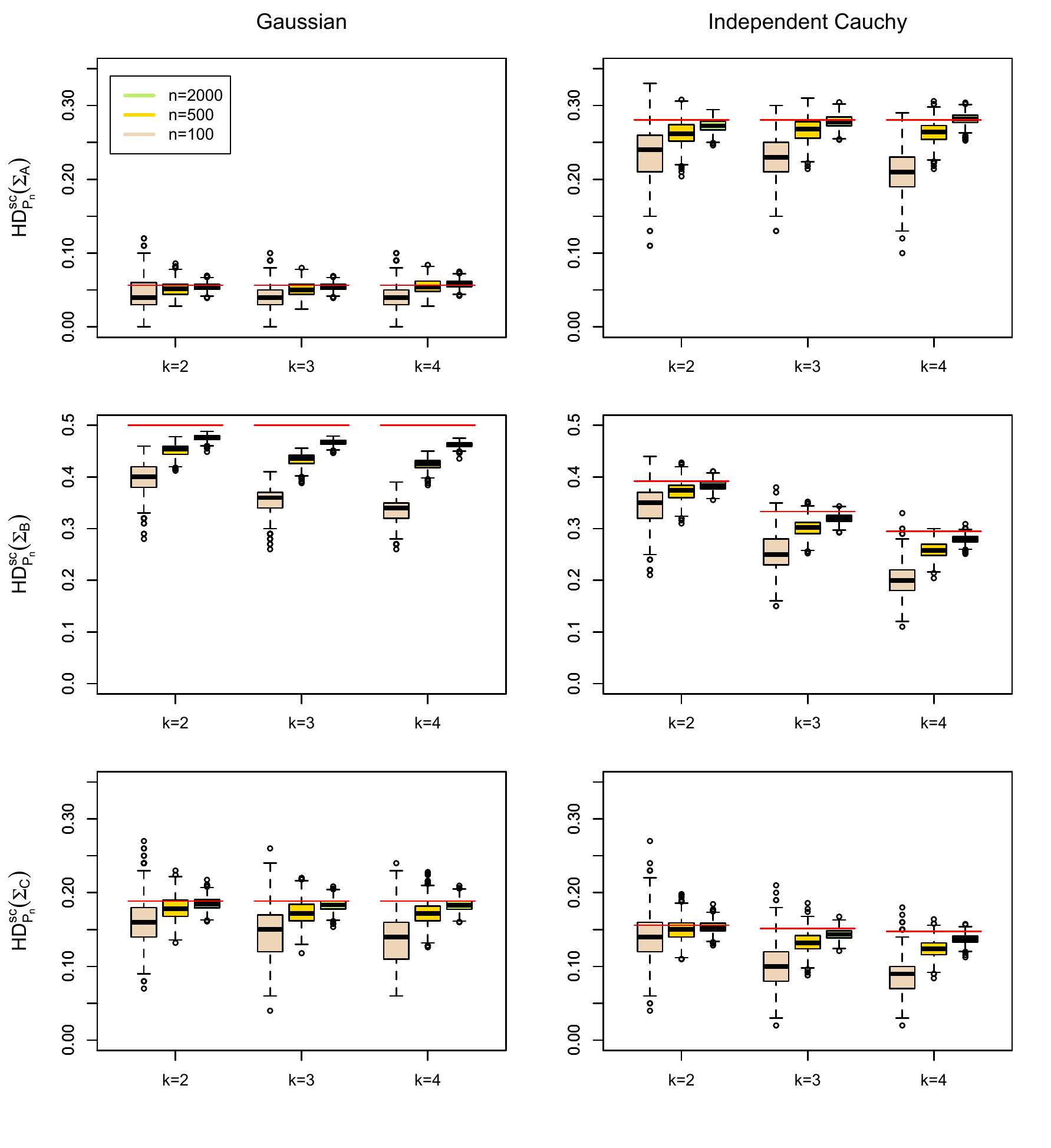}  
\vspace{-10mm}
\caption{Boxplots, for various values of~$n$ and $k$, of~$H\!D^{\rm sc}_{P_n}(\Sigma_A)$ (top row),~$H\!D^{\rm sc}_{P_n}(\Sigma_B)$ (middle row) and~$H\!D^{\rm sc}_{P_n}(\Sigma_C)$ (bottom row) based on $M=1000$ independent random samples of size~$n$ from the $k$-variate multinormal distribution with location zero and scatter~$I_k$ (left) or from the $k$-variate distribution with independent Cauchy marginals (right); we refer to Section~\ref{SecSupMonteCarlo} for the expressions of~$\Sigma_A$, $\Sigma_B$ and~$\Sigma_C$.}
\label{FigSec2b}  
\end{figure}

The Monte Carlo exercise we use to validate these expressions is the following. For each possible combination of~$n\in\{100,500,2000\}$ and~$k\in\{2,3,4\}$, we generated $M=1000$ independent random samples of size~$n$ (i) from the $k$-variate normal distribution with location zero and scatter~$I_k$ and (ii) from the $k$-variate distribution with independent Cauchy marginals. Letting~$\Lambda^A_k
:={\rm diag}(8,I_{k-1})$,
$\Lambda^B_k:=I_{k}$,
$\Lambda^C_k:={\rm diag}( {\textstyle{\frac{1}{8}}} ,I_{k-1})$,
and 
$$  
O_k 
:= 
{\rm diag}
\Bigg( 
{\,\frac{1}{\sqrt{2}}\  \frac{1}{\sqrt{2}}\, \choose \ \ \frac{-1}{\sqrt{2}}\ \ \frac{1}{\sqrt{2}}\,}
,
I_{k-2} 
\Bigg)
,
$$
we evaluated, in each sample, the depths~$H\!D^{\rm sc}_{P_n}(\Sigma_\ell)$ of~$\Sigma_\ell=O_k \Lambda^\ell_k  O_k'$, $\ell=A,B,C$, where~$P_n$ denotes the empirical probability measure associated with the $k$-variate sample of size~$n$ at hand (each evaluation of~$H\!D^{\rm sc}_{P_n}(\cdot)$ is done by approximating the infimum in~$u\in\mathcal{S}^{k-1}$ by a minimum over $N=10000$ directions randomly sampled from the uniform distribution over~$\mathcal{S}^{k-1}$). For each~$n$, $k$, and each underlying distribution (multinormal or independent Cauchy), Figure~\ref{FigSec2b} reports boxplots of the corresponding~$M$ values of~$H\!D^{\rm sc}_{P_n}(\Sigma_\ell)$, $\ell=A,B,C$. Clearly, the results support the theoretical depth expressions obtained in~(\ref{Gaussiandepth})-(\ref{Cauchytheo}), as well as the consistency result in Theorem~\ref{Consistency} (the bias for~$H\!D^{\rm sc}_{P_n}(\Sigma_B)$ in the Gaussian case is explained by the fact that, as we have seen above,~$\Sigma_B$ maximizes~$H\!D_P^{\rm sc}(\Sigma)$, with a maximal depth value equal to~$1/2$).    

\subsection{Illustrations of Theorem~\ref{Propquasiconcavity}}
\label{SecSupQuasi}

\begin{figure}[h!]
\includegraphics[width=\textwidth]{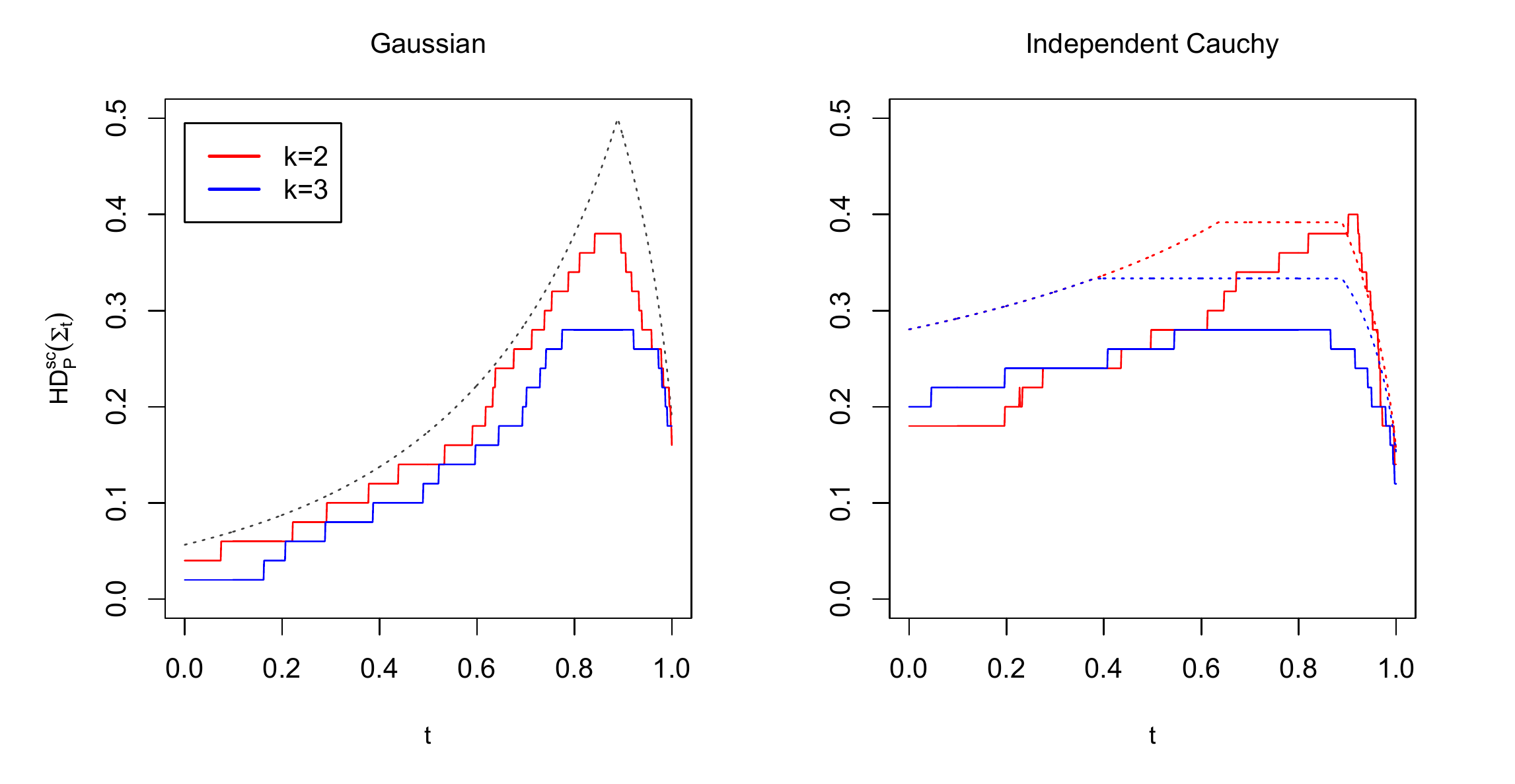}   
\vspace{-7mm}
 \caption{The dotted curves are the graphs of~$t\mapsto H\!D^{\rm sc}_P(\Sigma_t)$, where~$\Sigma_t=(1-t)\Sigma_A+t \Sigma_C$ involves the scatter matrices~$\Sigma_A$ and~$\Sigma_C$ considered in Figure~\ref{FigSec2b}, for the $k$-variate normal distribution with location zero and scatter~$I_2$ (left) and for the $k$-variate distribution with independent Cauchy marginals (right). The solid curves are associated with a single random sample of size~$n=50$ from the corresponding distributions. Red and blue correspond to the bivariate and trivariate cases, respectively; in the population Gaussian case, the graph of~$t\mapsto H\!D^{\rm sc}_P(\Sigma_t)$ is the same for~$k=2$ and~$k=3$, hence is plotted in black.}
\label{FigSec3} 
\end{figure}

Figure~\ref{FigSec3} plots, for~$k=2,3$, the graphs of~$t\mapsto H\!D^{\rm sc}_P(\Sigma_t)$ for $\Sigma_t:=(1-t)\Sigma_A+t \Sigma_C$, where~$\Sigma_A$ and~$\Sigma_C$ are the scatter matrices considered in the numerical exercise performed in Section~\ref{SecSupMonteCarlo} and where~$P$ is either the $k$-variate normal distribution with location zero and scatter~$I_k$ or the $k$-variate distribution with independent Cauchy marginals. The same figure also provides the corresponding sample plots, based on a single random sample of size~$n=50$ drawn from each of these two distributions. All plots are compatible with the quasi-concavity result in Theorem~\ref{Propquasiconcavity}. Figure~\ref{FigSec3} also illustrates the continuity of~$t\mapsto H\!D^{\rm sc}_{P,T}(\Sigma_t)$ for smooth probability measures~$P$ (Theorem~\ref{Propcontinuity}) and shows that continuity may fail to hold in the sample case.

\subsection{Illustrations of Theorems~\ref{Propquasiconcavityshape}-\ref{Propquasiconcavityharmonicshape}}
\label{SecSupShape}

\begin{figure}[h!] 
\vspace{6mm}
\includegraphics[width=.98\textwidth]{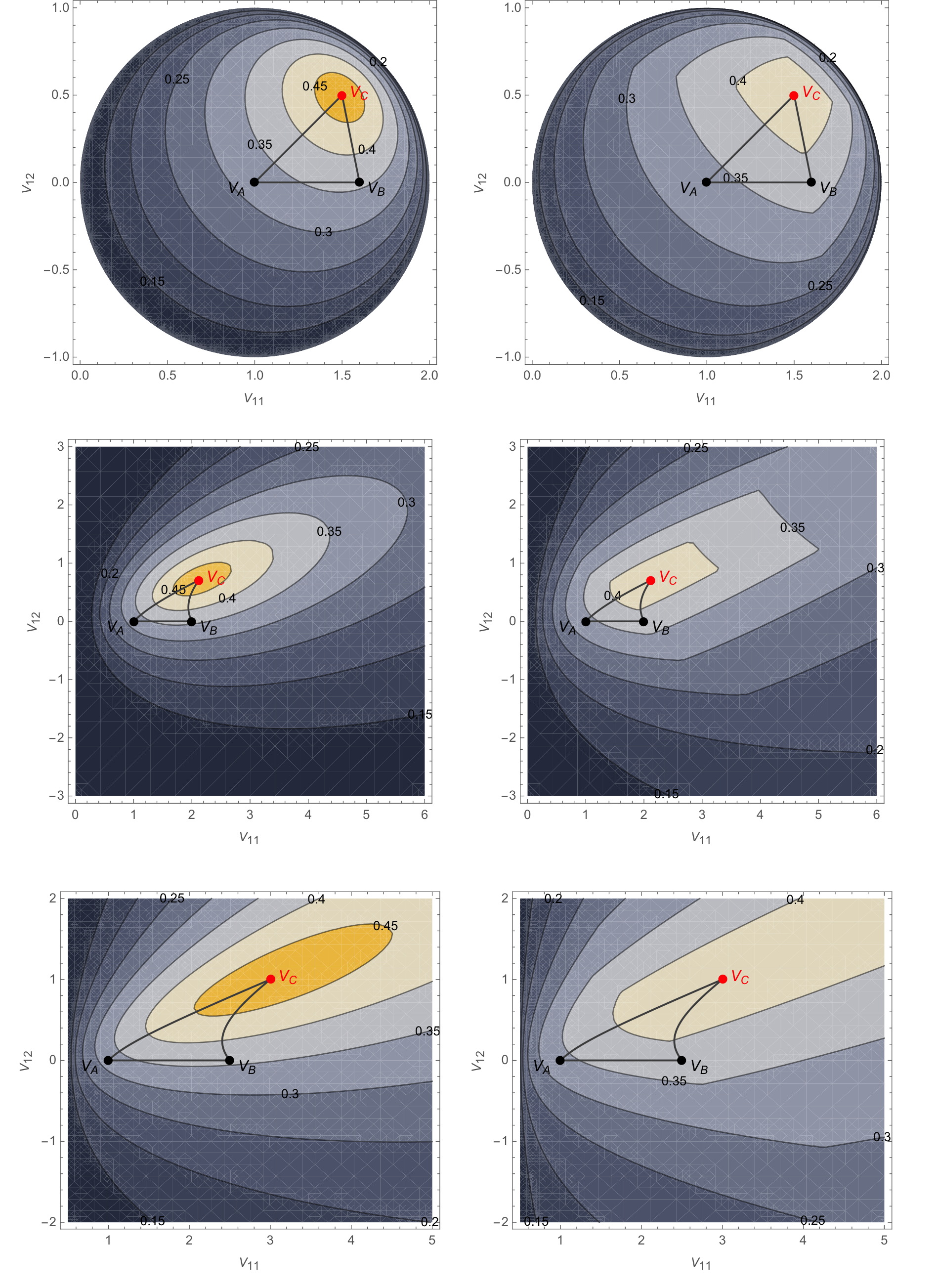}  
\vspace{-3mm}
\caption{Contour plots
\vspace{-.8mm}
  of $(V_{11},V_{12})\mapsto H\!D^{{\rm sh},S}_{P,T}(V_S)$, for $S_{\rm tr}$ (top), $S_{\det}$ (middle) and $S_{\rm tr}^*$ (bottom), where $H\!D^{{\rm sh},S}_{P,T}(\VS_S)$ is the shape halfspace 
depth, with respect to~$P$, of the $S$-shape~$V_S(\in\mathcal{P}_2^S)$ with upper-left entry~$V_{11}$ and upper-right entry~$V_{12}$, for an arbitrary affine-equivariant location functional~$T$. The probability measures~$P$ considered are those associated (i) with the 
\vspace{-.8mm}
 bivariate normal distribution with location zero and scatter~$\Sigma_C={3\ 1\, \choose \,1\ 1\,}$ (left)  and (ii) with the distribution of~$\Sigma_C^{1/2}Z$, where~$Z$ has mutually independent Cauchy marginals (right). In each case, the (``true") $S$-shape associated with~$\Sigma_C$ is marked in red and those associated with~$\Sigma_A$ and~$\Sigma_B$ from Figure~\ref{FigSec7a} are marked in black. Linear paths (top), geodesic paths (middle) and harmonic paths (bottom) between these three shapes are drawn.}
\label{FigSec7b} 
\end{figure}

Figure~\ref{FigSec7b} draws, for an arbitrary affine-equivariant location functional~$T$, contour plots of $(V_{11},V_{12})\mapsto H\!D^{{\rm sh},S}_{P,T}(V_S)$, for the scale functionals~$S_{\rm tr}$, $S_{\det}$ and $S_{\rm tr}^*$, where $H\!D^{{\rm sh},S}_{P,T}(\VS_S)$ is the shape halfspace 
depth, with respect to~$P$, of the $S$-shape~$V_S(\in\mathcal{P}_2^S)$ with upper-left entry~$V_{11}$ and upper-right entry~$V_{12}$. The probability measures~$P$ considered are those associated (i) with the bivariate normal distribution with location zero and scatter~$\Sigma_C={3\ 1\, \choose \,1\ 1\,}$, and (ii) with the distribution of~$\Sigma_C^{1/2}Z$, where~$Z$ has mutually independent Cauchy marginals. For~$S_{\rm tr}$, $S_{\det}$ and $S_{\rm tr}^*$, the figure also shows the linear, geodesic and harmonic paths, respectively, linking the (``true") $S$-shape associated with~$\Sigma_C$ and those associated with~$\Sigma_A
={1\ 0\, \choose \,0\ 1\,}$ and $\Sigma_B={4\ 0\, \choose \,0\ 1\,}$. The results illustrate the convexity of the regions~$R^{{\rm sh},S_{\rm tr}}_{P,T}(\alpha)$, along with the geodesic (resp., harmonic) convexity of the regions~$R^{{\rm sh},S_{\det}}_{P,T}(\alpha)$ (resp., $R^{{\rm sh},S_{\rm tr}^*}_{P,T}(\alpha)$).

\section{Proofs}
\label{SecSupProofs}


In the proofs below, we will often use the fact that
\begin{eqnarray*}
\lefteqn{
\hspace{-2mm} 
H\!D_{P,T}^{\rm sc}(\Sigma)
=
\inf_{u} 
\min
\big( 
P\big[ |u'(X-T_P)| \leq \sqrt{u'\Sigma u}\, \big]
,
P\big[ |u'(X-T_P)| \geq \sqrt{u'\Sigma u}\, \big]
\big)
}
\\[2mm]
& & 
\hspace{-3mm} 
=
\min\Big(\! 
\inf_{u} 
P\big[ |u'(X-T_P)| \leq \sqrt{u'\Sigma u}\, \big]
,
\inf_{u} 
P\big[ |u'(X-T_P)| \geq \sqrt{u'\Sigma u}\, \big]
\Big)
,
\end{eqnarray*}
where all infima are over the unit sphere~$\mathcal{S}^{k-1}$ of~$\R^k$. 
\vspace{1mm}

\subsection{Proofs from Section~2}
\label{AppSecScatter}

 \
 \vspace{2mm}
  
  

{\sc Proof of Theorem~\ref{thaffinv}.} 
Fix~$A\in GL_k$ and~$b\in\R^k$. By using the affine-equivariance of~$T$ (that is,~$T_{P_{A,b}}=AT_P+b$) and letting~$u_A:=A'u/\|A'u\|$, we obtain
\begin{eqnarray*}
\lefteqn{
H\!D^{\rm sc}_{P_{A,b},T}(A\Sigma A')
}
\\[2mm]
& &
\hspace{-5mm} 
=
\inf_{u\in\mathcal{S}^{k-1}}
\min\big( 
 P\big[ |u'A(X-T_P)| \leq \sqrt{u' A\Sigma A' u}\, \big]
,
 P\big[ |u'A(X-T_P)| \geq \sqrt{u' A\Sigma A' u}\, \big]
\big)
\\[2mm]
& &
\hspace{-5mm} 
=
\inf_{u\in\mathcal{S}^{k-1}}
\min\big( 
P\big[ |u_A'(X-T_P)| \leq 
{\textstyle{\sqrt{u_A' \Sigma u_A}}}
\, \big]
,
P\big[ |u_A'(X-T_P)| \geq 
{\textstyle{\sqrt{u_A' \Sigma u_A}}}
\, \big]
\big)
\\[2mm]
& &
\hspace{-5mm} 
=
H\!D^{\rm sc}_{P,T}(\Sigma)
,
\end{eqnarray*}
where the last equality follows from the fact that the mapping~$u\mapsto u_A$ is a one-to-one transformation of~$\mathcal{S}^{k-1}$.  
\cqfd
\vspace{3mm}

We now establish the explicit scatter halfspace depth expression~(\ref{Cauchytheo}) in the independent Cauchy case. For that purpose, consider again a random vector~$X=(X_1,\ldots,X_k)'$ with independent Cauchy marginals. Using the fact that, for any non-negative real numbers~$a_\ell$, $\ell=1,\ldots,k$, the random variables~$\sum_{\ell=1}^k a_{\ell} X_{\ell}$ and~$(\sum_{\ell=1}^k a_{\ell}) X_1$ then share the same distribution, we obtain that, denoting by~$\|x\|_1=\sum_{\ell=1}^k |x_{\ell}|$ the $L_1$-norm of~$x=(x_1,\ldots,x_k)'$,  
$$
P\big[ |u'X| \leq \sqrt{u'\Sigma u}\, \big]
=
P\big[ \|u\|_1\, |X_1| \leq \sqrt{u' \Sigma u}\, \big]
=
2 \Psi\bigg(  \frac{\sqrt{u' \Sigma u}}{\|u\|_1} \bigg) - 1
,
$$
where~$\Psi$ is the Cauchy cumulative distribution function. Therefore, if~$T$ is centro-equiva\-riant, we obtain
\begin{eqnarray*}
H\!D_{P,T}^{\rm sc}(\Sigma)
&\!\!=\!\!&
\min
\bigg(
2 \Psi\bigg(  \inf_{u\in\mathcal{S}^{k-1}} \frac{\sqrt{u' \Sigma u}}{\|u\|_1} \bigg) - 1
,
2 - 2 \Psi\bigg( \sup_{u\in\mathcal{S}^{k-1}}  \frac{\sqrt{u' \Sigma u}}{\|u\|_1} \bigg) 
\bigg)
\\[2mm]
&\!\!=\!\!&
2
\min
\bigg(
\Psi\big(   \min_v \sqrt{v' \Sigma v}\, \big) - {\textstyle\frac{1}{2}}
,
1 - \Psi\big(  \max_v \sqrt{v' \Sigma v} \, \big) 
\bigg)
,
\end{eqnarray*}
where the minimum and maximum in~$v$ are over the unit $L_1$-sphere~$\{ v\in\R^k : \|v\|_1=1\}$. The formula~(\ref{Cauchytheo}) then follows from the following  result.

\begin{Lem} 
\label{Cauchylem}
For any~$\Sigma\in\mathcal{P}_k$, the maximal and minimal values of~$v' \Sigma v$ when~$v$ runs over the 
$L_1$-sphere~$\{ v\in\R^k : \|v\|_1=\sum_{i=1}^k |v_i|=1\}$ are $\max({\rm diag}(\Sigma))$ and~$1/\max_s (s'\Sigma^{-1}s)$, respectively, where~$\max_s$ is the maximum over~$s=(s_1,\ldots,s_{k})\in\{-1,1\}^{k}$. 
\end{Lem}

{\sc Proof of Lemma~\ref{Cauchylem}.} 
We start with the following considerations. In dimension~$k$, the~$L_1$-sphere can be parametrized as
$$
v_{s,t}
=
(s_1t_1,s_2t_2,\ldots,s_{k-1} t_{k-1},s_{k}(1-t_1-\ldots-t_{k-1}))'
,
$$
where~$s=(s_1,\ldots,s_{k})\in\{-1,1\}^{k}$ and~$(t_1,\ldots,
t_{k-1})\in {\rm Simpl}_{k-1}:=\{ (t_1,\ldots,
\linebreak
t_{k-1}) : t_1 \geq 0,\ldots,t_{k-1}\geq 0,t_1+\ldots+t_{k-1}\leq 1\}$. 
Clearly, symmetry of the $L_1$-sphere and of the function to be maximized/minimized allows to restrict to~$s_{k}=1$. Now, for any given~$s$ of the form~$(s_1,\ldots,s_{k-1},1)$, consider the function 
\begin{eqnarray*}
f_s :\hspace{13mm}  {\rm Simpl}_{k-1} \hspace{3mm}  & \to & \R
\\[2mm]
t=(t_1,\ldots,t_{k-1})  & \mapsto & v_{s,t}' \Sigma v_{s,t}	
.
\end{eqnarray*}
This function is twice differentiable, with a gradient~$\nabla f_{s}(t)$ whose~$i$th component is
$$
\frac{\partial}{\partial t_i} f_s(t)
=
(0,...,0,s_i,0,\ldots,-1) \Sigma v_{s,t} 
+
v_{s,t}' \Sigma (0,...,0,s_i,0,\ldots,-1)'
$$
$$
=
2 (0,...,0,s_i,0,\ldots,-1) \Sigma v_{s,t} 
$$
and a Hessian matrix~$H_{s}(t)$ whose~$(i,j)$-entry is 
$$
\frac{\partial^2}{\partial t_it_j} f_s(t)
=
2
(0,...,0,s_i,0,\ldots,-1) \Sigma 
(0,...,0,s_j,0,\ldots,-1)'
.
$$
Now, for any~$z\in\R^{k-1}$, 
$$
z'H_s(t)z
=
2
(s_1z_1,...,s_kz_k,-z_1-\ldots-z_k) 
\Sigma 
(s_1z_1,...,s_kz_k,-z_1-\ldots-z_k) 
'
\geq 
0
,
$$
with equality if and only if~$z=0$. Therefore, $f_s$ is strictly convex               over~${\rm Simpl}_{k-1}$.  

Let us start with the maximum. For a given~$s$, strict convexity of~$f_s$ implies that the maximum of~$f_s$ can only be achieved at~$e_{i,k-1}$, $i=1,\ldots,k-1$, where~$e_{i,\ell}$ stands for the $i$th vector of the canonical basis of~$\R^{\ell}$, or at~$0(\in\R^{k-1})$. Since~$f_s(e_{i,k-1})=\Sigma_{ii}$, $i=1,\ldots,k-1$, and~$f_s(0)=\Sigma_{kk}$, it follows that the maximal value of~$f_s$ over~${\rm Simpl}_{k-1}$ is~$\max({\rm diag}(\Sigma))$. Since this is the case for any~$s$, we conclude that the maximum of~$v'\Sigma v$ over the~$L_1$-sphere is itself~$\max({\rm diag}(\Sigma))$.  

Let us then turn to the minimum. The minimum of~$f_s$, when extended into a (still strictly convex) function defined on~$\R^{k-1}$, is the solution of the gradient conditions
$$
(0,...,0,s_i,0,\ldots,-1) \Sigma v_{s,t} 
=
0
,
\qquad
i=1,\ldots,k-1
.
$$
Writing simply~$e_k=e_{k,k}$ for the $k$th vector of the canonical basis of~$\R^k$ and letting
$$
S_s
=
\left(
\begin{array}{cccc}
s_1 & & & -1 \\[1mm]
& \ddots & &\vdots \\[1mm]
& & s_{k-1} &-1 
\end{array}
\right)
,
$$
these gradient conditions rewrite
$
S_s \Sigma 
(e_{k} + S_s't )
=
0
$
(recall that we restricted to~$s_k=1$),
and their unique solution in~$\R^{k-1}$ is
$
t_{s}^{\rm min}
:=
- (S_s \Sigma S_s')^{-1} S_s \Sigma e_{k}
.
$

It will be useful below to have a more explicit expression of~$t_{s}^{\rm min}$. Note that the gradient conditions above state that~$\Sigma(e_{k}+S_s' t_{s}^{\rm min})$ is in the null space of~$S_s$. Since this null space is easily checked to be~$\{\lambda s:\lambda\in\R\}$, this implies that~$e_{k}+S_s' t_{s}^{\rm min}=\lambda\Sigma^{-1}s$ for some~$\lambda\in\R$. Premultiplying both sides of this equation by~$s'$, we obtain $1=\lambda s'\Sigma^{-1}s$, which yields~$\lambda=1/(s'\Sigma^{-1}s)$. Thus, 
\begin{equation}
	\label{touseforminval}
	e_{k}+S_s' t_{s}^{\rm min}
=
\frac{
1
}
{
s'\Sigma^{-1}s
}
\,
\Sigma^{-1}s
.
\end{equation}
In the first $k-1$ components, this yields (after multiplication by~$s_{i}$)
\begin{equation}
\label{ti}
(t_{s}^{\rm min})_i
=
\frac{
s_{i}e_i'\Sigma^{-1}s
}
{
s'\Sigma^{-1}s
}
,
\qquad
i=1,\ldots,k-1
,
\end{equation}
while the $k$th component provides (still with~$s_k=1$)
\begin{equation}
\label{tk}
1 - \sum_{i=1}^{k-1}
(t_{s}^{\rm min})_i
=
\frac{
e_k'\Sigma^{-1}s
}
{
s'\Sigma^{-1}s
}
=
\frac{
s_{k} e_k'\Sigma^{-1}s
}
{
s'\Sigma^{-1}s
}
\cdot
\end{equation}
Strict convexity of~$f_s$ implies that its minimal value over~$\R^{k-1}$ is~$f_{s}(t_{s}^{\rm min})$. By using~(\ref{touseforminval}), this minimal value takes the form
$$
f_s(t_{s}^{\rm min})
=
v_{s,t_{s}^{\rm min}}' \Sigma v_{s,t_{s}^{\rm min}}
=
(e_{k} + S_s't_{s}^{\rm min} )' \Sigma (e_{k} + S_s't_{s}^{\rm min} )
=
\frac{1}{s' \Sigma^{-1} s}
\cdot
$$

Now, consider an arbitrary sign $k$-vector~$s_*$ that maximizes~$s' \Sigma^{-1} s$ among the~$2^{k-1}$ corresponding sign vectors~$s$ to be considered (the last component of~$s$ is still fixed to one). Assume for a moment that~$t_{s_*}^{\rm min}$ is in the interior of ${\rm Simpl}_{k-1}$. Since it is the minimal value of~$f_{s*}$ over~$\R^{k-1}$, $f_{s_*}(t_{s_*}^{\rm min})=1/(s_*' \Sigma^{-1} s_*)$ of courses minimizes~$f_{s*}$ over~${\rm Simpl}_{k-1}$. Pick then another sign vector~$s$. By construction, $1/(s_*' \Sigma^{-1} s_*)$ is smaller than or equal to~$f_{s}(t_{s}^{\rm min})=1/(s' \Sigma^{-1} s)$, which, as the minimal value of~$f_s$ when extended to~$\R^{k-1}$, can only be smaller than or equal to the minimal value of~$f_s$ over~${\rm Simpl}_{k-1}$. Therefore, $1/(s_*' \Sigma^{-1} s_*)$ is then the minimal value of~$v'\Sigma^{-1}v$ over the unit $L_1$-sphere. 

It thus remains to show that~$t_{s_*}^{\rm min}$ indeed belongs to the interior of~${\rm Simpl}_{k-1}$. Equivalently (in view of~(\ref{ti})-(\ref{tk})), it remains to show that 
$
s_{*i}e_i'\Sigma^{-1}s_* > 0
$
for~$i=1,\ldots, k$. Assume then that~$s_{*\ell}e_\ell'\Sigma^{-1}s_*\leq 0$ for some~$\ell$. Defining~$s_{**}$ as the vector obtained from~$s_*$ by only changing the sign of its $\ell$th component, that is, putting~$s_{**}:=s_* - 2s_{*\ell}e_\ell$, we have
\begin{eqnarray*}
\lefteqn{
s_{**}'\Sigma^{-1}s_{**}
-
s_*'\Sigma^{-1}s_*
=
(s_* - 2s_{*\ell}e_\ell)'\Sigma^{-1}(s_* - 2s_{*\ell}e_\ell)
-
s_*'\Sigma^{-1}s_*
}	
\\[2mm]
& & 
\hspace{13mm} 
=
-4 s_{*\ell} e_\ell' 
\Sigma^{-1}
(s_* - s_{*\ell}e_\ell)
=
-4 s_{*\ell} e_\ell' \Sigma^{-1} s_*
+4 e_\ell' \Sigma^{-1} e_\ell
>
0
,
\end{eqnarray*}
which contradicts the maximality property of~$s_*$. 
\cqfd
\vspace{3mm}


The proof of Theorem~\ref{Consistency} requires Lemmas~\ref{LemConsistency}-\ref{LemConsistency3} below. Before proceeding, we introduce the inner and outer ``slabs"
$$
H^{\rm in}_{u,c}:=\{x\in\R^k:|u'x|\leq c\}
\quad
\textrm{ and }
\quad
H^{\rm out}_{u,c}:=\{x\in\R^k:|u'x|\geq c\}.
$$ 
Henceforth, the superscript ``in/out" is to be read as ``in (resp., out)". We will further write~$B(\theta,r)$ for the ball~$\{x\in\R^k: \|x-\theta\|< r \}$ and~$\bar{B}(\theta,r)$ for its closure.

\begin{Lem}
\label{LemConsistency}
Let~$P,Q$ be two probability measures over~$\R^k$. Define 
\linebreak
$H\!D^{\rm in/out}_{P,T}(\Sigma):=\inf_{u\in\mathcal{S}^{k-1}} P[T_P+\Sigma^{1/2} H_{u}^{\rm in/out}]$, with~$H_{u}^{\rm in/out}:=H_{u,1}^{in/out}$. Then, for any~$\Sigma\in\mathcal{P}_k$,
\begin{eqnarray*}
\lefteqn{	
|H\!D^{\rm in/out}_{P,T}(\Sigma) - H\!D^{\rm in/out}_{Q,T}(\Sigma)|
}
\\[3mm]
& & 
\hspace{5mm} 
\leq
\sup_{C\in\mathcal{C}^{\rm in/out}} 
| 
P[ C ]
 - 
Q[ C ]
|
+
\sup_{C\in\mathcal{C}_0^{\rm in/out}} 
|
P[T_Q+C]
- 
P[T_P+C]
|
,
\end{eqnarray*}
where we let~$\mathcal{C}^{\rm in/out}:=\{ \theta+H_{u,c}^{\rm in/out} : (\theta,u,c)\in \R^k\times \mathcal{S}^{k-1} \times \R^+_0\}$ and~$\mathcal{C}_0^{\rm in/out}:=\{H_{u,c}^{\rm in/out} : (u,c)\in \mathcal{S}^{k-1}\times\R^+_0\}$. 
\end{Lem}

{\sc Proof of Lemma~\ref{LemConsistency}.}
We prove only the ``in" result, since the proof of the ``out" result is entirely similar. First assume that~$H\!D^{\rm in}_{P,T}(\Sigma) \geq  H\!D^{\rm in}_{Q,T}(\Sigma)$. Then, for any~$\varepsilon>0$, there exists~$u_0=u_0(\Sigma,Q,\varepsilon)$ such that~$Q[T_Q+\Sigma^{1/2}H_{u_0}^{\rm in}] \leq H\!D^{\rm in}_{Q,T}(\Sigma)+\varepsilon$, so that
\begin{eqnarray*}
\lefteqn{
\hspace{-1mm} 
|H\!D^{\rm in}_{P,T}(\Sigma) - H\!D^{\rm in}_{Q,T}(\Sigma)|
=
H\!D^{\rm in}_{P,T}(\Sigma)- H\!D^{\rm in}_{Q,T}(\Sigma)
}
\\[2mm]
& & 
\hspace{7mm} 
\leq
P[T_P+\Sigma^{1/2}H_{u_0}^{\rm in}]
- 
Q[T_Q+\Sigma^{1/2}H_{u_0}^{\rm in}]
+\varepsilon
\\[2mm]
& & 
\hspace{7mm} 
\leq
P[T_Q+\Sigma^{1/2}H_{u_0}^{\rm in}]
- 
Q[T_Q+\Sigma^{1/2}H_{u_0}^{\rm in}]
\\[2mm]
& & 
\hspace{23mm} 
+
P[T_P+\Sigma^{1/2}H_{u_0}^{\rm in}]
- 
P[T_Q+\Sigma^{1/2}H_{u_0}^{\rm in}]
+
\varepsilon
\\[2mm]
& & 
\hspace{7mm} 
\leq
\sup_{C\in\mathcal{C}^{\rm in}} 
| 
P[C]
 - 
Q[C]
|
+
\sup_{C\in\mathcal{C}_0^{\rm in}} 
|
P[T_Q+C]
- 
P[T_P+C]
|
+
\varepsilon
.
\end{eqnarray*}
Similarly, if~$H\!D^{\rm in}_{P,T}(\Sigma) \leq  H\!D^{\rm in}_{Q,T}(\Sigma)$, then, for any~$\varepsilon>0$, there exists~$u_1=u_1(\Sigma,P,\varepsilon)$ such that~$P[T_P+\Sigma^{1/2}H_{u_1}^{\rm in}] \leq H\!D^{\rm in}_{P,T}(\Sigma)+\varepsilon$, so that
\begin{eqnarray*}
\lefteqn{
\hspace{-1mm} 
|H\!D^{\rm in}_{P,T}(\Sigma) - H\!D^{\rm in}_{Q,T}(\Sigma)|
=
H\!D^{\rm in}_{Q,T}(\Sigma)- H\!D^{\rm in}_{P,T}(\Sigma)
}
\\[2mm]
& & 
\hspace{7mm} 
\leq
Q[T_Q+\Sigma^{1/2}H_{u_1}^{\rm in}]
- 
P[T_P+\Sigma^{1/2}H_{u_1}^{\rm in}]
+\varepsilon
\\[2mm]
& & 
\hspace{7mm} 
\leq
Q[T_Q+\Sigma^{1/2}H_{u_1}^{\rm in}]
- 
P[T_Q+\Sigma^{1/2}H_{u_1}^{\rm in}]
\\[2mm]
& & 
\hspace{23mm} 
+
P[T_Q+\Sigma^{1/2}H_{u_1}^{\rm in}]
- 
P[T_P+\Sigma^{1/2}H_{u_1}^{\rm in}]
+
\varepsilon
\\[2mm]
& & 
\hspace{7mm} 
\leq
\sup_{C\in\mathcal{C}^{\rm in}} 
| 
P[C]
 - 
Q[C]
|
+ 
\sup_{C\in\mathcal{C}_0^{\rm in}} 
|
P[T_Q+C]
- 
P[T_P+C] 
|
+
\varepsilon
.
\end{eqnarray*}
Since, in both cases, the result holds for any~$\varepsilon>0$, the result is proved.
\cqfd
\vspace{3mm}

\begin{Lem}
\label{LemConsistency2}
Let~$P$ be a probability measure over~$\R^k$ and~$K$ be a compact subset of~$\R^k$. For any~$c>0$, let
\vspace{-.8mm}
 $
 s^K_{P}(c):=\sup_{(\theta,u)\in K \times \mathcal{S}^{k-1}} P[ |u'(X-\theta)| \leq  c ]$ and write~$s^K_{P}:=s^K_{P}(0)$.
Then
(i)
$s^K_{P}(c)\to s^K_{P}$  as~$c\stackrel{>}{\to} 0$
and 
(ii)~$s^K_{P}=P[ u_0'(X-\theta_0)=0 ]$ for some~$(\theta_0,u_0)\in K \times \mathcal{S}^{k-1}$. 
\end{Lem}

{\sc Proof of Lemma~\ref{LemConsistency2}}.
Clearly, $s^K_{P}(c)$ is increasing in~$c$ over~$[0,\infty)$, which guarantees that
$
\tilde{s}^K_{P}
:=
\lim_{c\stackrel{>}{\to} 0} 
s^K_{P}(c)
$
exists and satisfies~$\tilde{s}^K_{P}\geq s^K_{P}$. 
Now, fix an arbitrary decreasing sequence~$(c_n)$ converging to~$0$ and consider a sequence~$((\theta_n,u_n))$ in~$K \times \mathcal{S}^{k-1}$ such that
$$
P[|u_n'(X-\theta_n)\leq c_n]\geq s^K_{P}(c_n)-(1/n).
$$
Compactness of $K \times \mathcal{S}^{k-1}$ guarantees the existence of a subsequence~$((\theta_{n_\ell},u_{n_\ell}))$ that converges in~$K \times \mathcal{S}^{k-1}$, to~$(\theta_0,u_0)$ say. Clearly, without loss of generality, we can assume that~$(u_{n_\ell}' u_0)$ is an increasing sequence and that~$\|\theta_{n_\ell}-\theta_0\|$ is a decreasing sequence (if that is not the case, one can always extract a further subsequence which meets these monotonicity properties). Let then~$\bar{B}_\ell:=\bar{B}(\theta_0,\|\theta_{n_\ell}-\theta_0\|)$ and~$C_\ell:=\{ u\in\mathcal{S}^{k-1} : u'u_0 \geq u_{n_\ell}' u_0 \}$. Clearly,~$\bar{B}_\ell$  and~$C_\ell$ are decreasing sequences of sets, with~$\cap_\ell \bar{B}_\ell=\{\theta_0\}$ and~$\cap_\ell C_\ell=\{u_0\}$. Therefore,
\begin{eqnarray*}
	\lim_{\ell\to\infty} 
r_{\ell}
&\!\!:=\!\!&
\lim_{\ell\to\infty} 
P[ X\in \cup_{\theta\in \bar{B}_\ell}\cup_{u\in C_\ell} \{ x=\theta+y: |u'y|\leq c_{n_\ell} \} ] 
\\[2mm]
&\!\!=\!\!&
P[ u_0' (X-\theta_0)=0 ] 
.
\end{eqnarray*}
Now, for any~$\ell$, 
$
r_{\ell}
\geq  
P[ |u_{n_\ell}'(X-\theta_{n_\ell})|\leq c_{n_\ell} ]  
\geq 
s^K_{P}(c_{n_\ell})-(1/n_\ell)
$,
which implies that~$s^K_{P} \geq P[ |u_0' (X-\theta_0)|=0 ]  \geq \tilde{s}^K_{P}$. Therefore, $\tilde{s}^K_{P} = s^K_{P}=P[ |u_0' (X-\theta_0)|=0 ]$. 
\cqfd

\begin{Lem}
\label{LemConsistency3}
Let~$P$ be a smooth probability measure over~$\R^k$ and fix~$\theta_0\in\R^k$. Then
$$
\sup_{(u,c)\in \mathcal{S}^{k-1}\times\R^+_0}
|
P[\theta+H_{u,c}^{\rm in/out}]
- 
P[\theta_0+H_{u,c}^{\rm in/out}]
|
\to
0
$$
as~$\theta\to\theta_0$.
\end{Lem}

{\sc Proof of Lemma~\ref{LemConsistency3}}.
We start with the ``in'' result. Fix $\varepsilon>0$. 
Pick~$c_1>0$ large enough to have~$P[B(\theta_0,c_1/2)]\geq 1-(\varepsilon/2)$. Pick then~$c_0\in(0,c_1)$ such that $s^{\bar{B}(\theta_0,2c_1)}_{P}(c_0)<\varepsilon/2$ for any~$c\in(0,c_0]$ (existence of such a~$c_0$ is guaranteed by Lemma~\ref{LemConsistency2} and the smoothness assumption on~$P$). Then,
$$
\sup_{(u,c)\in \mathcal{S}^{k-1}\times\R^+_0}
|
P[\theta+H_{u,c}^{\rm in}]
- 
P[\theta_0+H_{u,c}^{\rm in}]
|
=
\max
\big(
Q_{0,c_0}^{\rm in},Q_{c_0,c_1}^{\rm in},Q_{c_1,\infty}^{\rm in}
\big)
,
$$
where
$Q_{0,c_0}^{\rm in}$, $Q_{c_0,c_1}^{\rm in}$ and $Q_{c_1,\infty}^{\rm in}$ are the suprema of 
$|
P[\theta+H_{u,c}^{\rm in}]
- 
P[\theta_0+H_{u,c}^{\rm in}]
|$ over $u\in \mathcal{S}^{k-1}$ and, respectively, $c\in(0,c_0]$, $c\in(c_0,c_1]$ and $c\in(c_1,\infty)$. Now, fix~$\delta\in (0,\min(c_0,c_1/2))$ and~let~$\theta\in B(\theta_0,\delta)$. 

(i) The choice of $c_0$ implies
\begin{eqnarray}
	Q_{0,c_0}^{\rm in}
	&\leq &
	\sup_{u\in \mathcal{S}^{k-1}}  P[\theta+H_{u,c_0}^{\rm in}]
	+
	\sup_{u\in \mathcal{S}^{k-1}}  P[\theta_0+H_{u,c_0}^{\rm in}]
	\nonumber
	\\[2mm]
	&\leq &
	 2 \sup_{(\theta,u)\in \bar{B}(\theta_0,2c_1)\times \mathcal{S}^{k-1}} P[\theta+H_{u,c_0}^{\rm in}]
	\leq 
	\varepsilon
	.
	\label{jobpartQ1}
\end{eqnarray}
	
(ii) Let $u\in \mathcal{S}^{k-1}$ and $c\geq c_0$. Assume, without loss of generality, that $u'(\theta-\theta_0)\geq 0$ (the case $u'(\theta-\theta_0)\leq 0$ proceeds similarly). The set $\theta+H_{u,c}^{\rm in}$ rewrites
$$
	\theta+H_{u,c}^{\rm in}=\{\theta+x:|u'x|\leq c\}=\{\theta_0+y:|u'y-(u'(\theta-\theta_0))|\leq c\}.
$$
Therefore,
\begin{eqnarray}
\lefteqn{	
\hspace{0mm} 
|
P[\theta+H_{u,c}^{\rm in}]
- 
P[\theta_0+H_{u,c}^{\rm in}]
|
}
\nonumber
\\[2mm]
& &
\hspace{10mm} 
\leq
P\Big[\{\theta_0+y:-c\leq u'y\leq -c+u'(\theta-\theta_0)\}\Big]
\nonumber
\\[2mm]
& &
\hspace{40mm} 
 +
P\Big[\{\theta_0+y:c\leq u'y\leq c+u'(\theta-\theta_0)\}\Big]
\nonumber
\\[2mm]
& &
\hspace{10mm} 
\leq
P\Big[\{\theta_0+y:-c\leq u'y\leq -c+\delta\}\Big]
\nonumber
\\[2mm]
& &
\hspace{40mm} 
 +
P\Big[\{\theta_0+y:c\leq u'y\leq c+\delta\}\Big].
\label{crucialinbound}
\end{eqnarray}
(iia) For~$u\in \mathcal{S}^{k-1}$ and $c_0< c\leq c_1$, set $\theta_1=\theta_0-cu+\delta u/2$ and $\theta_2=\theta_0+cu+\delta u/2$. It holds
 $\{\theta_0+y:-c\leq u'y\leq -c+\delta\}=\{\theta_1+x:|u'x|\leq \delta/2\}$ and $\{\theta_0+y:c\leq u'y\leq c \delta\}=\{\theta_2+x:|u'x|\leq \delta/2\}$. Since~$s^{\bar{B}(\theta_0,2c_1)}_P(\delta/2)\leq s^{\bar{B}(\theta_0,2c_1)}_P(c_0)<\varepsilon/2$ and $\|\theta_\ell-\theta_0\|\leq c+\delta/2<2c_1$ ($\ell=1,2$), (\ref{crucialinbound}) yields
\begin{equation}
	Q_{c_0,c_1}^{\rm in}
	\leq 
	2\sup_{(\theta,u)\in \bar{B}(\theta_0,2c_1)\times \mathcal{S}^{k-1}} P\Big[\{\theta+x:|u'x|\leq \delta/2\}\Big]\leq \varepsilon.
		\label{jobpartQ2}
\end{equation}
(iib) For~$u\in \mathcal{S}^{k-1}$ and $c>c_1$, the sets $\{\theta_0+y:-c\leq u'y\leq -c+\delta\}$ and $\{\theta_0+y:c\leq u'y\leq c+\delta\}$ lie outside the ball $B(\theta_0,c_1/2)$ since $-c+\delta<-c+(c_1/2)<-c_1/2$ and~$c>c_1>c_1/2$, respectively. Therefore, it follows from~(\ref{crucialinbound}) that
\begin{equation}
Q_{c_1,\infty}^{\rm in}
\leq 
2 P\big[ \R^k\setminus B(\theta_0,c_1/2)]
\leq \varepsilon
.	
		\label{jobpartQ3}
\end{equation}
According to~(\ref{jobpartQ1}), (\ref{jobpartQ2}) and (\ref{jobpartQ3}), all three quantities $Q_{0,c_0}^{\rm in}$, $Q_{c_0,c_1}^{\rm in}$ and $Q_{c_1,\infty}^{\rm in}$ are bounded by $\varepsilon$ as soon as $\|\theta-\theta_0\|<\delta$, which concludes the proof of the ``in" result.  

The proof of the ``out" result proceeds similarly. For the same choices of $c_0$, $c_1$ and $\delta$, and the respective suprema~$Q_{0,c_0}^{\rm out}$, $Q_{c_0,c_1}^{\rm out}$ and $Q_{c_1,\infty}^{\rm out}$, it holds
\begin{eqnarray*}
	Q_{c_1,\infty}^{\rm out}
	&\leq &
	\sup_{u\in \mathcal{S}^{k-1}} P[\theta+H_{u,c_1}^{\rm out}]
	+
	\sup_{u\in \mathcal{S}^{k-1}}  P[\theta_0+H_{u,c_1}^{\rm out}]
	\nonumber
	\\[2mm]
	&\leq&
	 2 \sup_{(\theta,u)\in \bar{B}(\theta_0,c_1/2)\times \mathcal{S}^{k-1}} P\big[
	 \R^k\setminus B(\theta_0,c_1/2)
	 \big]
	\leq 
	\varepsilon
	.
\end{eqnarray*}
Moreover, the inequality $Q_{0,c_0}^{\rm out}\leq \varepsilon$ follows from the fact that
$$
|
P[\theta+H_{u,c}^{\rm out}]
- 
P[\theta_0+H_{u,c}^{\rm out}]
|
\leq 
P[\theta+H_{u,c_0}^{\rm in}]
+ 
P[\theta_0+H_{u,c_0}^{\rm in}]
$$
for~$c\leq c_0$. Finally, it can be proved that~$Q_{c_0,c_1}^{\rm out}\leq \varepsilon$ along the exact same lines as above. The ``out" result follows.
\cqfd
\vspace{3mm}

{\sc Proof of Theorem~\ref{Consistency}}.
The collection~$\mathcal{H}$ of all halfspaces in~$\R^k$ is a Vapnik-Chervo\-nenkis class; see, e.g., page~152 of~\cite{WelVan1996}. Hence, Lemma~2.6.17 of the same implies that~$\mathcal{H}\sqcap \mathcal{H}:=\{ H_1 \cap H_2: H_1,H_2\in \mathcal{H}\}$ and~$\mathcal{H}\sqcup \mathcal{H}:=\{ H_1 \cup H_2: H_1,H_2\in \mathcal{H}\}$ are also Vapnik-Chervonenkis classes. Consequently, using henceforth the notation from Lemma~\ref{LemConsistency}, $\mathcal{C}^{\rm in}(\subset \mathcal{H}\sqcap \mathcal{H}$) and~$\mathcal{C}^{\rm out}(\subset \mathcal{H}\sqcup \mathcal{H})$ are themselves Vapnik-Chervonenkis classes,
which implies that  
\begin{equation}
\label{quasic}
\sup_{C\in\mathcal{C}^{\rm in}} 
| 
P_n[ C ]
 - 
P[ C ]
|
\to
0
\quad
\textrm{and}
\quad
\sup_{C\in\mathcal{C}^{\rm out}} 
| 
P_n[ C ]
 - 
P[ C ]
|
\to
0
\end{equation}
almost surely as~$n\to\infty$.
Also, since~$T_{P_n}\to T_P$ almost surely as~$n\to\infty$, Lemma~\ref{LemConsistency3} entails
\begin{equation}
\label{quasic2}
\sup_{C\in\mathcal{C}_0^{\rm in}} 
|
P[T_{P_n}+C]
- 
P[T_P+C]
|
\to
0
\quad
\textrm{and}
\quad
\sup_{C\in\mathcal{C}_0^{\rm out}} 
|
P[T_{P_n}+C]
- 
P[T_P+C]
|
\to
0
\end{equation}
almost surely as~$n\to\infty$.

Now, by using Lemma~\ref{LemConsistency}, we obtain that, for any~$\Sigma\in\mathcal{P}_k$, 
\begin{eqnarray*}
\lefteqn{
\hspace{2mm} 
|H\!D^{\rm sc}_{P_n,T}(\Sigma)- H\!D^{\rm sc}_{P,T}(\Sigma)|
}
\\[3mm]
& &
\hspace{3mm} 
=
|
\min(H\!D_{P_n,T}^{\rm in}(\Sigma),H\!D^{\rm out}_{P_n,T}(\Sigma))
-
\min(H\!D^{\rm in}_{P,T}(\Sigma),H\!D^{\rm out}_{P,T}(\Sigma))
|
\\[3mm]
& &
\hspace{3mm} 
\leq
\max
\big(
|H\!D^{\rm in}_{P_n,T}(\Sigma)-H\!D^{\rm in}_{P,T}(\Sigma)|
,
|H\!D^{\rm out}_{P_n,T}(\Sigma)-H\!D^{\rm out}_{P,T}(\Sigma)|
\big)
\\[3mm]
& &
\hspace{3mm} 
\leq
\max
\Big(
\sup_{C\in\mathcal{C}^{\rm in}} 
| 
P_n[ C ]
 - 
P[ C ]
|
+
\sup_{C\in\mathcal{C}_0^{\rm in}} 
|
P[T_{P_n}+C]
- 
P[T_P+C]
|
,
\\[1mm]
& &
\hspace{23mm} 
\sup_{C\in\mathcal{C}^{\rm out}} 
| 
P_n[ C ]
 - 
P[ C ]
|
+
\sup_{C\in\mathcal{C}_0^{\rm out}} 
|
P[T_{P_n}+C]
- 
P[T_P+C]
|
\Big)
.
\end{eqnarray*}
Consequently, 
\begin{eqnarray*}
\lefteqn{
\sup_{\Sigma\in\mathcal{P}_k}
|H\!D^{\rm sc}_{P_n,T}(\Sigma)-H\!D^{\rm sc}_{P,T}(\Sigma)|
}
\\[2mm]
& &
\hspace{3mm} 
\leq 
\max
\Big(
\sup_{C\in\mathcal{C}^{\rm in}} 
| 
P_n[ C ]
 - 
P[ C ]
|
+
\sup_{C\in\mathcal{C}_0^{\rm in}} 
|
P[T_{P_n}+C]
- 
P[T_P+C]
|
,
\\[1mm]
& &
\hspace{18mm} 
\sup_{C\in\mathcal{C}^{\rm out}} 
| 
P_n[ C ]
 - 
P[ C ]
|
+
\sup_{C\in\mathcal{C}_0^{\rm out}} 
|
P[T_{P_n}+C]
- 
P[T_P+C]
|
\Big)
,
\end{eqnarray*}
which, in view of~(\ref{quasic}) and~(\ref{quasic2}), establishes the result.
\cqfd
\vspace{3mm}

 We close this section by proving that the Tukey median~$\theta_P$ is strongly consistent without any assumption on~$P$. 
 
\begin{Lem} 
\label{LemConvergence}
Let~$P$ be a probability measure over~$\R^k$ and~$P_n$ denote the empirical measure associated with a random sample of size~$n$ from~$P$. Then $\theta_{P_n}\to\theta_P$ almost surely as~$n\to\infty$. 
\end{Lem}

{\sc Proof of Lemma~\ref{LemConvergence}}.
For any~$\theta\in\R^k$ and any probability measure~$Q$ over~$\R^k$, denote by~$H\!D^{\rm loc}_Q(\theta)$ the location halfspace depth of~$\theta$ with respect to~$Q$. Recall that we defined~$\theta_Q$ as the barycentre of~$M^{\rm loc}_Q=\{\theta\in\R^k: H\!D^{\rm loc}_Q(\theta)=\max_{\eta \in\R^k}H\!D^{\rm loc}_Q(\eta)\}$. 
It is well known that~$\theta\mapsto H\!D^{\rm loc}_Q(\theta)$ is upper semicontinuous; see, e.g., Lemma~6.1 in \cite{DonGas1992}. In general, this function is not uniquely maximized at~$\theta_Q$. However, it is easy to define a modified depth function~$\theta\mapsto H\!D^{\rm loc}_{Q,{\rm mod}}(\theta)$ that is still upper semicontinuous, agrees with~$\theta\mapsto H\!D^{\rm loc}_Q(\theta)$ on~$\R^k/M^{\rm loc}_Q$, and for which~$\theta_Q$ is the unique maximizer. In view of the uniform consistency of location halfspace depth (see, e.g., (6.2) and~(6.6) in \citealp{DonGas1992}), the result then follows from Theorem~2.12 and Lemma 14.3 in \cite{Kos2008}.
\cqfd
\vspace{3mm}



\subsection{Proofs from Section~3}
\label{AppSecFrobenius}


\
\vspace{3mm}
 
{\sc Proof of Theorem~\ref{Propcontinuity}.} (i) 
Fix~$u\in\mathcal{S}^{k-1}$. Since~$H_u^{\rm in}:=\{ x\in\R^k : |u'x|\leq 1\}$ is a closed subset of~$\R^k$, the mapping~$P\mapsto P[ H_u^{\rm in} ]$ is upper semicontinuous for weak convergence. Now, Slutzky's lemma entails that, as~$d_F(\Sigma,\Sigma_0)\to 0$, the measure defined by~$B\mapsto P[ T_P + \Sigma^{1/2}B]$ converges weakly to the one defined by~$B\mapsto P[ T_P + \Sigma_0^{1/2}B]$. Therefore,~$\Sigma\mapsto P[ T_P + \Sigma^{1/2}H_u^{\rm in}]$ is upper $F$-semicontinuous at~$\Sigma_0$. Since~$H_u^{\rm out}:=\{ x\in\R^k : |u'x|\geq 1\}$ is also a closed subset of~$\R^k$, the same argument shows that~$\Sigma\mapsto P[ T_P + \Sigma^{1/2}H_u^{\rm out}]$ is upper $F$-semicontinuous at~$\Sigma_0$. Therefore
$$
\Sigma
\mapsto
H\!D_{P,T}^{\rm sc}(\Sigma)
=
\min\Big( 
\inf_{u\in\mathcal{S}^{k-1}} 
P[ T_P + \Sigma^{1/2} H_u^{\rm in} ]
,
\inf_{u\in\mathcal{S}^{k-1}} 
P[ T_P + \Sigma^{1/2} H_u^{\rm out} ]
\Big)
,
$$
is upper $F$-semicontinuous (recall that the infimum of a collection of upper semicontinuous functions is upper semicontinuous). 

(ii) The result directly follows from the fact that~$R^{\rm sc}_{P,T}(\alpha)$ is the inverse image of~$[\alpha,+\infty)$ by the upper $F$-semicontinuous function~$\Sigma\mapsto H\!D_{P,T}^{\rm sc}(\Sigma)$.  

(iii) Fix a sequence~$(\Sigma_n)$ in~$\mathcal{P}_k$ converging to~$\Sigma_0$ with respect to the Frobenius distance. With the same notation as in the proof of~(i), note that, for any~$\Sigma$, 
$$
H\!D_{P,T}^{\rm sc}(\Sigma)
\
=
\inf_{u\in\mathcal{S}^{k-1}} 
\min\Big( 
P[ T_P + \Sigma^{1/2} H_u^{\rm in} ]
,
P[ T_P + \Sigma^{1/2} H_u^{\rm out} ]
\Big)
.
$$
For any~$n$, pick then~$u_n(\in\mathcal{S}^{k-1})$ such that 
$$
\min\Big( 
P[T_P + \Sigma_n^{1/2} H_{u_n}^{\rm in} ]
,
P[T_P + \Sigma_n^{1/2} H_{u_n}^{\rm out} ]
\Big)
\leq 
H\!D^{\rm sc}_{P,T}(\Sigma_n)+\frac{1}{n}
\cdot 
$$
Compactness of~$\mathcal{S}^{k-1}$ implies that we can extract a subsequence~$(u_{n_\ell})$ of~$(u_n)$ that converges to~$u_0(\in\mathcal{S}^{k-1})$. Writing~$\mathbb{I}[C]$ for the indicator function of the set~$C$, the dominated convergence theorem then yields that
\begin{eqnarray*}
	\lefteqn{
P[T_P + \Sigma_{n_\ell}^{1/2} H_{u_{n_\ell}}^{\rm in} ]
-
P[T_P + \Sigma_0^{1/2} H_{u_0}^{\rm in} ]
}
\\[2mm]
& & 
\hspace{13mm} 
=
\int_{\R^k}
(
\mathbb{I}[T_P + \Sigma_{n_\ell}^{1/2} H_{u_{n_\ell}}^{\rm in}]
-
\mathbb{I}[T_P + \Sigma_0^{1/2} H_{u_0}^{\rm in}]
)
\,dP
\to 0
\end{eqnarray*}
as~$\ell\to\infty$ (the smoothness
\vspace{-1mm}
  assumption on~$P$ guarantees that~$\mathbb{I}[T_P + \Sigma_{n_\ell}^{1/2} H_{u_{n_\ell}}^{\rm in}]
-
\mathbb{I}[T_P + \Sigma_0^{1/2} H_{u_0}^{\rm in}]
\to 0$ $P$-almost everywhere). Proceeding in the same way, we obtain that~$P[T_P + \Sigma_{n_\ell}^{1/2} H_{u_{n_\ell}}^{\rm out} ]
-
P[T_P + \Sigma_0^{1/2} H_{u_0}^{\rm out} ]
\to 0
$
as~$\ell\to\infty$. Consequently,
\begin{eqnarray*}
\liminf_{n\to\infty} H\!D^{\rm sc}_{P,T}(\Sigma_n)
&\!\!=\!\!&
\liminf_{n\to\infty} 
\min\Big( 
P[T_P + \Sigma_n^{1/2} H_{u_n}^{\rm in} ]
,
P[T_P + \Sigma_n^{1/2} H_{u_n}^{\rm out} ]
\Big)
\\[2mm]
&\!\!=\!\!&
\liminf_{\ell\to\infty} 
\min\Big( 
P[T_P + \Sigma_{n_\ell}^{1/2} H_{u_{n_\ell}}^{\rm in} ]
,
P[T_P + \Sigma_{n_\ell}^{1/2} H_{u_{n_\ell}}^{\rm out} ]
\Big)
\\[2mm]
&\!\!=\!\!&
\min\Big( 
P[T_P + \Sigma_0^{1/2} H_{u_0}^{\rm in} ]
,
P[T_P + \Sigma_0^{1/2} H_{u_0}^{\rm out} ]
\Big)
\\[2mm]
&\!\!\geq\!\!&
H\!D^{\rm sc}_{P,T}(\Sigma_0)
.
\end{eqnarray*}
We conclude that, if~$P$ is smooth at~$T_p$, then~$\Sigma\to H\!D_{P,T}^{\rm sc}(\Sigma)$ is also lower $F$-semicontinuous, hence $F$-continuous.
\cqfd
\vspace{3mm}


{\sc Proof of Theorem~\ref{boundedness}.}
Fix~$\alpha>0$. Note that~$\lambda_1(\Sigma)\geq \|\Sigma\|_F/\sqrt{k} \geq  (\|\Sigma-I_k\|_F- \|I_k\|_F)/\sqrt{k}$ for any~$\Sigma\in\mathcal{P}_k$.  
Therefore, denoting by~$v_1(\Sigma)$ an arbitrary unit eigenvector associated with~$\lambda_1(\Sigma)$, we have that, for any~$\Sigma\notin B_F(I_k,r)$,
\begin{eqnarray*}
\lefteqn{	
H\!D_{P,T}^{\rm sc}(\Sigma)
\leq
\inf_{u\in\mathcal{S}^{k-1}}
P\big[ |u'(X-T_P)| \geq {\textstyle{\sqrt{u'\Sigma u}}}\, \big]
}
\\[2mm]
& & 
\hspace{13mm} 
\leq 
P[ |v'_1(\Sigma)(X-T_P)| \geq  \sqrt{\lambda_1(\Sigma)}]
\leq 
P\Big[ \|X-T_P\| \geq  {\textstyle \frac{(r-1)^{1/2}}{k^{1/4}} } \Big]
,
\end{eqnarray*}
which can be made strictly smaller than $\alpha$ for~$r$ large enough.  
This confirms that, for~$r$ large enough, $R^{\rm sc}_{P,T}(\alpha)$ is included in the ball~$B_F(I_k,r)$, hence is $F$-bounded.  
\cqfd
\vspace{3mm}
  

{\sc Proof of Theorem~\ref{Propquasiconcavity}.} 
(i) With~$\Sigma_t=(1-t)\Sigma_a+t \Sigma_b$, we clearly have that, for any~$u\in\mathcal{S}^{k-1}$, $\min(u'\Sigma_a u,u' \Sigma_b u)\leq u' \Sigma_t u \leq \max(u'\Sigma_a u,u' \Sigma_b u)$. This entails that, for any~$u\in\mathcal{S}^{k-1}$,
\begin{eqnarray*}
\lefteqn{
P\big[ |u'(X-T_P)| \leq 
{\textstyle{ \sqrt{u' \Sigma_t u} }}
\, \big]
}
\\[2mm]
& &
\hspace{13mm} 
\geq
\min(
P\big[ |u'(X-T_P)| \leq 
{\textstyle{ \sqrt{u' \Sigma_a u} }}
\, \big]
,
P\big[ |u'(X-T_P)| \leq 
{\textstyle{ \sqrt{u' \Sigma_b u} }}
\, \big]
)
\\[2mm]
& &
\hspace{13mm} 
\geq
\min(H\!D^{\rm sc}_{P,T}(\Sigma_a),H\!D^{\rm sc}_{P,T}(\Sigma_b))
\end{eqnarray*}
and
\begin{eqnarray*}
\lefteqn{
P\big[ |u'(X-T_P)| \geq 
{\textstyle{ \sqrt{u' \Sigma_t u} }}
\, \big]
}
\\[2mm]
& &
\hspace{13mm} 
\geq
\min(
P\big[ |u'(X-T_P)| \geq 
{\textstyle{ \sqrt{u' \Sigma_a u} }}
\, \big]
,
P\big[ |u'(X-T_P)| \geq 
{\textstyle{ \sqrt{u' \Sigma_b u} }}
\, \big]
)
\\[2mm]
& &
\hspace{13mm} 
\geq
\min(H\!D^{\rm sc}_{P,T}(\Sigma_a),H\!D^{\rm sc}_{P,T}(\Sigma_b))
.
\end{eqnarray*}
The result follows. 
(ii) If both~$\Sigma_a,\Sigma_b\in R^{\rm sc}_{P,T}(\alpha)$, then Part~(i) of the result entails that, for any~$t\in[0,1]$, $H\!D^{\rm sc}_{P,T}(\Sigma_t)\geq \min(H\!D^{\rm sc}_{P,T}(\Sigma_a),H\!D^{\rm sc}_{P,T}(\Sigma_b))\geq \alpha$, so that~$\Sigma_t\in R^{\rm sc}_{P,T}(\alpha)$.  
\cqfd
\vspace{3mm}


\subsection{Proofs from Section~4}
\label{AppSecGeod}


For the sake of completeness, we prove the following result.

\begin{Lem}
\label{totalboundlemma}
Let~$R$ be a $g$-bounded subset of~$\mathcal{P}_k$. Then~$R$ is totally $g$-bounded, that is, for any~$\varepsilon$, there exist~$\Sigma_i$, $i=1,\ldots,m=m(\varepsilon)$ such that~$R\subset \cup_{i=1}^m B_g(\Sigma_i,\varepsilon)$. 
\end{Lem}

{\sc Proof of Lemma~\ref{totalboundlemma}.}
 As a mapping from the metric space~$(\mathcal{S}_k,d_F)$ (recall that~$d_F$ denotes the Frobenius distance) to the metric space~$(\mathcal{P}_k,d_g)$, $A\mapsto \exp(A)$ is continuous; see the proof of Proposition~10 in \cite{BhaHol06}. Denoting, for any~$A\in\mathcal{S}_k$, as~${\rm vech}(A)$ the vector obtained by stacking the upper-diagonal entries of~$A$ on top of each other, the mapping~$v\mapsto {\rm vech}^{-1}(v)$ from~$(\R^{k(k+1)/2},d_E)$ (equipped with the usual Euclidean distance~$d_E$) to~$(\mathcal{S}_k,d_F)$ is trivially continuous, so that the mapping~$f: (\R^{k(k+1)/2},d_E)\to (\mathcal{P}_k,d_g): v \mapsto  f(v):= \exp({\rm vech}^{-1}(v))$ is also continuous. 
 
Now, fix~$\varepsilon>0$, pick~$r>0$ such that~$R$ is included in the closed ball~$\bar{B}:=\bar{B}_g(I_k,r):=\{ \Sigma\in \mathcal{P}_k: d_g(\Sigma,I_k)\leq r\}$, and consider the resulting open covering~$\{B_g(\Sigma,\varepsilon):\Sigma \in \bar{B}\}$ of~$\bar{B}$. From continuity, $\mathcal{C}:=\{f^{-1}(B_g(\Sigma,\varepsilon)):\Sigma \in \bar{B}\}$ is an open covering of the closed set~$f^{-1}(\bar{B})$ in~$\R^{k(k+1)/2}$. It is easy to check that, for any~$\Sigma\in\bar{B}$, $\lambda_1(\Sigma)\leq \exp(r/\sqrt{k})$, so that~$f^{-1}(\bar{B})$ is bounded, hence compact. Therefore, a finite subcovering~$\{f^{-1}(B_g(\Sigma_i,\varepsilon)):i=1,\ldots,m\}$ of~$f^{-1}(\bar{B})$ can be extracted from~$\mathcal{C}$, which provides the desired finite covering~$\{B_g(\Sigma_i,\varepsilon):i=1,\ldots,m\}$ of~$\bar{B}$, hence of~$R$, with open $g$-balls of radius~$\varepsilon$.   
\cqfd
\vspace{3mm}


{\sc Proof of Theorem~\ref{geodboundedness}.}
Assume first that~$s_{P,T}<1/2$ and fix~$\varepsilon>0$. We will then prove that~$R^{\rm sc}_{P,T}(s_{P,T}+\varepsilon)$ is $g$-bounded by showing that, for~$r>0$ large enough, it is included in the $g$-ball~$B_g(I_k,r)$. To do so, first note that~(\ref{defindist}) entails
\begin{eqnarray}
\lefteqn{
\hspace{-10mm}
d_g(\Sigma,I_k)
=
\sqrt{\sum_{i=1}^k (\log \lambda_i(\Sigma))^2}
\leq
\sqrt{k}\, \max( |\log \lambda_1(\Sigma)| , |\log \lambda_k(\Sigma)|)
}
\nonumber
\\[2mm]
& & 
\hspace{20mm}
=
\sqrt{k}\, \max( \log \lambda_1(\Sigma) , \log \lambda^{-1}_k(\Sigma))
.
\label{distmaj}
\end{eqnarray}
Therefore, $\Sigma\notin B_g(I_k,r)$ implies that 
(i) 
$
\lambda_1(\Sigma)  > \exp(r/\sqrt{k})
$
or
(ii) 
$
\lambda_k(\Sigma)  < \exp(-r/\sqrt{k})
$
(or both). In case (i), 
\begin{eqnarray*}
H\!D_{P,T}^{\rm sc}(\Sigma)
&\!\! \leq \!\! &
\inf_{u\in\mathcal{S}^{k-1}}
P\big[ |u'(X-T_P)| \geq 
{\textstyle{\sqrt{u'\Sigma u}}}\, \big]
\\[1mm]
&\!\! \leq \!\! &
P[ |v'_1(\Sigma)(X-T_P)| \geq  \sqrt{\lambda_1(\Sigma)}]
\\[2mm]
&\!\! \leq \!\! &
P[ |v'_1(\Sigma)(X-T_P)| \geq  \exp(r/2\sqrt{k})]
\\[2mm]
&\!\! \leq \!\! &
P[ \| X-T_P \| \geq  \exp(r/2\sqrt{k})]
,
\end{eqnarray*}
which can be made smaller than $\varepsilon$ (hence, smaller than~$s_{P,T}+\varepsilon$) for~$r$ large enough.  
In case~(ii), we have that, using the notation~$s_P^K(\cdot)$ from Lemma~\ref{LemConsistency2},
\begin{eqnarray*}
H\!D_{P,T}^{\rm sc}(\Sigma)
&\!\! \leq \!\! &
\inf_{u\in\mathcal{S}^{k-1}}
P\big[ |u'(X-T_P)| \leq 
{\textstyle{\sqrt{u'\Sigma u}}}\, \big]
\\[2mm]
&\!\! \leq \!\! &
P\big[ |v_k'(\Sigma)(X-T_P)| \leq \lambda_k^{1/2}(\Sigma)\big]
\\[3mm]
&\!\! \leq \!\! &
 s_P^{\{T_P\}}(\lambda^{1/2}_k(\Sigma))
\leq s_P^{\{T_P\}}(\exp(-r/2\sqrt{k}))
,
\end{eqnarray*}
which, in view of Lemma~\ref{LemConsistency2}(i), can be made smaller than~$s_{P,T}+\varepsilon$ for~$r$ large enough. We conclude that, for~$\alpha>s_{P,T}$, $R^{\rm sc}_{P,T}(\alpha)$ is $g$-bounded, hence also (Lemma~\ref{totalboundlemma}) totally $g$-bounded. Since it is also $g$-closed (which follows from Theorem~\ref{geodcontinuity}(ii)), it is $g$-compact (recall from Section~\ref{SecGeod} that, in a complete metric space, any closed and totally bounded set is compact). 

Finally, if~$s_{P,T}\geq 1/2$, then, with~$u_0\in\mathcal{S}^{k-1}$ such that~$P[|u_0' (X-T_P)|=0]=s_{P,T}$ (existence is guaranteed in Lemma~\ref{LemConsistency2}(ii); take~$K=\{T_P\}$ there), we have
$
H\!D_{P,T}^{\rm sc}(\Sigma) 
\leq
P[ |u_0'(X-T_P)| \geq \sqrt{u_0'\Sigma u_0} ]
\leq 
P[ |u_0'(X-T_P)| >0 ]
=
1-
s_{P,T}
,
$
so that~$R^{\rm sc}_{P,T}(\alpha)$ is empty for any~$\alpha>1-s_{P,T}=\alpha_{P,T}$.
\cqfd
\vspace{3mm}
 

{\sc Proof of Theorem~\ref{geodmaxdepth}.}
By assumption, $R^{\rm sc}_{P,T}(\alpha_{P,T})$ is non-empty, so that~$\alpha_{*P,T}=\sup_{\Sigma\in\mathcal{P}_k} H\!D_{P,T}^{\rm sc}(\Sigma)\geq \alpha_{P,T}$.
If~$\alpha_{*P,T}=\alpha_{P,T}$, then the result holds since the maximal depth~$\alpha_{*P,T}$ is achieved at any scatter matrix in the non-empty set~$R^{\rm sc}_{P,T}(\alpha_{P,T})$. Assume then that~$\alpha_{*P,T}>\alpha_{P,T}$. Fix~$\delta>0$ such that~$\alpha_{*P,T}-\delta>\alpha_{P,T}$, so that~$R^{\rm sc}_{P,T}(\alpha_{*P,T}-\delta)$ is $g$-compact (Theorem~\ref{geodboundedness}). For any positive integer~$n$, it is possible to pick a scatter matrix~$\Sigma_n$ with~$H\!D^{\rm sc}_{P,T}(\Sigma_n)\geq \alpha_{*P,T}-(\delta/n)$. The $g$-compactness of~$R^{\rm sc}_{P,T}(\alpha_{*P,T}-\delta)$ implies that there exists a subsequence~$(\Sigma_{n_\ell})$ that $g$-converges in~$R^{\rm sc}_{P,T}(\alpha_{*P,T}-\delta)$, to~$\Sigma_*$ say. For any~$\varepsilon\in(0,\delta)$, all terms of~$(\Sigma_{n_\ell})$ are eventually in the $g$-closed set~$R^{\rm sc}_{P,T}(\alpha_{*P,T}-\varepsilon)$, so that its $g$-limit~$\Sigma_*$ must also belong to~$R^{\rm sc}_{P,T}(\alpha_{*P,T}-\varepsilon)$. For any such~$\varepsilon$, we thus have $\alpha_{*P,T}-\varepsilon\leq H\!D^{\rm sc}_{P,T}(\Sigma_*)\leq \alpha_{*P,T}$, which proves that~$H\!D^{\rm sc}_{P,T}(\Sigma_*)=\alpha_{*P,T}$.  
\cqfd
\vspace{3mm}


The proof of Theorem~\ref{propmaxCauchy} requires the following preliminary result.

\begin{Lem}
\label{LemmaxCauchy}
Fix~$\Sigma\in\mathcal{P}_k$ with~$\max({\rm diag}(\Sigma))\leq 1$. Then $\max_s s'\Sigma^{-1}s\geq k$ $($where~$\max_s$ is the maximum over~$s=(s_1,\ldots,s_{k})\in\{-1,1\}^{k})$, with equality if and only if~$\Sigma=I_k$. 	
\end{Lem}
 
{\sc Proof of Lemma~\ref{LemmaxCauchy}.}
We prove the result by induction. Clearly, the result holds for~$k=1$. Assume then that the result holds for~$k$. Writing
$$
\Sigma
=
\left(
\begin{array}{cc}
	\Sigma_{-} & v \\[1mm]
	v' & \Sigma_{k+1,k+1}
\end{array}
\right)
\quad
\textrm{ and }
\quad
s
=
\bigg(
\begin{array}{c}
s_- \\[0mm]
s_{k+1}
\end{array}
\bigg)
,
$$
the classical formula for the inverse of a block partitioned matrix yields
\begin{eqnarray}
s' \Sigma^{-1} s
& \!\!=\!\! &
s_-' \Sigma_-^{-1} s_-
+
\frac{
(s_-'\Sigma_{-}^{-1} v- s_{k+1})^2
}{\Sigma_{k+1,k+1}-v'\Sigma_{-}^{-1}v}
\nonumber
\\[2mm]
& \!\!=\!\! &
s_-' \Sigma_-^{-1} s_-
+
1
+
\frac{
(s_-'\Sigma_{-}^{-1} v- s_{k+1})^2-(\Sigma_{k+1,k+1}-v'\Sigma_{-}^{-1}v)
}{\Sigma_{k+1,k+1}-v'\Sigma_{-}^{-1}v}
\cdot		
\label{jslehq}
\end{eqnarray}
By induction assumption, there exists~$s_-$ such that
\begin{equation}
\label{toshowlemmax}
s' \Sigma^{-1} s
\geq 
k
+
1
+
\frac{
(s_-'\Sigma_{-}^{-1} v- s_{k+1})^2-(\Sigma_{k+1,k+1}-v'\Sigma_{-}^{-1}v)
}{\Sigma_{k+1,k+1}-v'\Sigma_{-}^{-1}v}
\cdot
\end{equation}
Now, irrespective of~$s_-$, choosing~$s_{k+1}=-{\rm sign}(s_-'\Sigma_{-}^{-1} v)$ yields
\begin{eqnarray*}
\lefteqn{
(s_-'\Sigma_{-}^{-1} v- s_{k+1})^2
- (\Sigma_{k+1,k+1}-v'\Sigma_{-}^{-1}v)
}
\\[2mm]
& & 
\hspace{3mm} 
=
(s_-'\Sigma_{-}^{-1} v)^2 + 2 |s_-'\Sigma_{-}^{-1} v| + v'\Sigma_{-}^{-1}v 
+ 1 - \Sigma_{k+1,k+1}
\geq 0
,
\end{eqnarray*}
since~$\Sigma_{k+1,k+1}\leq 1$. Jointly with~(\ref{toshowlemmax}), this provides~$\max_s s'\Sigma^{-1}s\geq k+1$. 

Now, assume that~$\max_s s'\Sigma^{-1}s=k+1$. We consider two cases. (a) $\max_{s_-} s_-'\Sigma_-^{-1} s_- >k$. Pick an arbitrary~$s_{*-}$ such that $s_{*-}'\Sigma_-^{-1} s_{*-} >k$. Then, with~$s_*=(s_{*-}',s_{*,k+1})'$, we have 
$$
s_*' \Sigma^{-1} s_*
>
k
+
1
+
\frac{
(s_{*-}'\Sigma_{-}^{-1} v- s_{*,k+1})^2-(\Sigma_{k+1,k+1}-v'\Sigma_{-}^{-1}v)
}{\Sigma_{k+1,k+1}-v'\Sigma_{-}^{-1}v}
\cdot
$$
Choosing again~$s_{*,k+1}=-{\rm sign}(s_{*-}'\Sigma_{-}^{-1} v)$ makes the third term of the righthand side non-negative, which implies that~$\max_s s'\Sigma^{-1}s>k+1$, a contradiction. (b) $\max_{s_-} s_-'\Sigma_-^{-1} s_- 
= k$. By induction assumption, we must then have~$\Sigma_-=I_k$. For any~$s=(s_-',s_{k+1})'$, (\ref{jslehq}) thus yields 
$$
s' \Sigma^{-1} s
=
s_-' s_-
+
\frac{
(s_-' v- s_{k+1})^2
}{\Sigma_{k+1,k+1}-v'v}
=
k
+
\frac{
(s_-' v- s_{k+1})^2
}{\Sigma_{k+1,k+1}-v'v}
\cdot
$$
Since~$\max_s s'\Sigma^{-1}s=k+1$, we must have that
$$
1
=
\max_s
\frac{
(s_-' v- s_{k+1})^2
}{\Sigma_{k+1,k+1}-v'v}
=
\frac{
(1+\sum_{\ell=1}^k |v_\ell|)^2
}{\Sigma_{k+1,k+1}-v'v}
=:
\frac{c}{d}
\cdot
$$
Since~$c\geq 1$ and~$d\leq 1$ (recall that~$\Sigma_{k+1,k+1}\leq 1$), this imposes that~$c=d=1$, which leads to~$v=0$ and~$\Sigma_{k+1,k+1}=1$. Jointly with~$\Sigma_-=I_k$, this shows that we must have~$\Sigma=I_{k+1}$, which establishes the result. 
\cqfd
\vspace{3mm}

{\sc Proof of Theorem~\ref{propmaxCauchy}.}
First note that, with~$\Sigma_*=\sqrt{k}I_k$, (\ref{Cauchytheo}) yields
$$
H\!D_{P,T}^{\rm sc}(\Sigma_*) 
=
2
\min
\Big(
\Psi\big(  k^{-1/4} \big) - \textstyle{\frac{1}{2}}
,
1 - \Psi\big( k^{1/4} \big) 
\Big)
=
\frac{2}{\pi} \arctan\big( k^{-1/4} \big)
$$
and fix an arbitrary~$\Sigma\in\mathcal{P}_k$. If~$\max({\rm diag}(\Sigma))>\sqrt{k}$, then 
\begin{equation}
\label{majo1}
H\!D_{P,T}^{\rm sc}(\Sigma)
\leq 
2
\big(1 - \Psi\big( \sqrt{\max({\rm diag}(\Sigma))} \big)\big)
<
2
\big(1 - \Psi\big( k^{1/4} \big)\big)
=
H\!D_{P,T}^{\rm sc}(\Sigma_*)
.	
\end{equation}
If~$\max({\rm diag}(\Sigma))\leq\sqrt{k}$, then Lemma~\ref{LemmaxCauchy} yields
$$
\max_s s' \Sigma^{-1} s
=
k^{-1/2} \max_s s' (k^{-1/2}\Sigma)^{-1} s
\geq
k^{1/2}
,
$$
so that
\begin{eqnarray}
\lefteqn{
H\!D_{P,T}^{\rm sc}(\Sigma)
\leq 
2
\big(\Psi\big( 1/ \max_s \sqrt{s' \Sigma^{-1} s}\, \big) - \textstyle{\frac{1}{2}} \big)
}
\nonumber
\\[2mm]
& & 
\hspace{28mm} 
\leq 
2
\big(\Psi\big( k^{-1/4} \big)-\textstyle{\frac{1}{2}} \big)
=
H\!D_{P,T}^{\rm sc}(\Sigma_*)
.
	\label{majo2}
\end{eqnarray}
We conclude that~$H\!D_{P,T}^{\rm sc}(\Sigma_*)\geq H\!D_{P,T}^{\rm sc}(\Sigma)$ for any~$\Sigma\in\mathcal{P}_k$. Now, assume that~$H\!D_{P,T}^{\rm sc}(\Sigma)=H\!D_{P,T}^{\rm sc}(\Sigma_*)$ for some~$\Sigma\in\mathcal{P}_k$. If~$\max({\rm diag}(\Sigma))>\sqrt{k}$, we can only have~$H\!D_{P,T}^{\rm sc}(\Sigma)<H\!D_{P,T}^{\rm sc}(\Sigma_*)$, as showed in~(\ref{majo1}). Thus we must have~$\max({\rm diag}(\Sigma))\leq\sqrt{k}$, and by assumption, all inequalities in~(\ref{majo2}) should be equalities. In view of Lemma~\ref{LemmaxCauchy}, this implies that~$k^{-1/2}\Sigma=I_k$, which establishes the result. 
\cqfd
\vspace{3mm}


\subsection{Proofs from Sections~5 and~6}
\label{AppSecAxiom}


\
\vspace{3mm}

{\sc Proof of Theorem~\ref{Fishconst}.} 
We start with the case~$\theta_0=0$ and~$\Sigma_0=I_k$, for which
\begin{eqnarray*}
\lefteqn{
\min( 
P\big[ |u'X| \leq 
{\textstyle{\sqrt{u'\Sigma u}}}
\,\big]
,
P\big[ |u'X| \geq 
{\textstyle{\sqrt{u'\Sigma u}}}
\,\big]
)
}
\\[2mm]
& & 
\hspace{13mm} 
=
\min( 
P[ |X_1| \leq \sqrt{u'\Sigma u}  ]
,
P[ |X_1| \geq \sqrt{u'\Sigma u}  ]
)
\end{eqnarray*}
for any~$u\in\mathcal{S}^{k-1}$, so that (note that the affine equivariance of~$T_P$ entails that~$T_P=0$)
\begin{eqnarray*}
\lefteqn{
H\!D_{P,T}^{\rm sc}(\Sigma)
=
\inf_{z\in {\rm Sp}(\Sigma)}
\min( 
P[ X_1^2 \leq z  ]
,
P[ X_1^2 \geq z  ]
)
}
\\[2mm]
& & 
\hspace{18mm} 
\leq
\min( 
P[ X_1^2 \leq 1  ]
,
P[ X_1^2 \geq 1  ]
)
=
H\!D_{P,T}^{\rm sc}(I_k)
,
\end{eqnarray*}
where the equality holds if and only if~${\rm Sp}(\Sigma)\subset \mathcal{I}_{\rm MSD}[X_1]$. The result for a general location~$\theta_0$ and scatter~$\Sigma_0$ readily follows from affine invariance and the identity~${\rm Sp}(AB)={\rm Sp}(BA)$.  
(ii) By definition, if~$\mathcal{I}_{\rm MSD}[Z_1]$ is a singleton, then this singleton must be~$\{1\}$. Consequently, if~$H\!D_{P,T}^{\rm sc}(\Sigma)=H\!D_{P,T}^{\rm sc}(\Sigma_0)$, then Part~(i) of the result entails that~$\lambda_k(\Sigma_0^{-1}\Sigma)=\lambda_1(\Sigma_0^{-1}\Sigma)=1$. This implies that~$\Sigma_0^{-1}\Sigma=I_k$, which establishes the result. 
\cqfd
\vspace{3mm}
  

We turn to the proofs of Theorems~\ref{propgeodesicquasiconc} and~\ref{propharmonicquasiconc}, that require the following preliminary results for the elliptical case (Lemma~\ref{lemgeodesicharmonicquasiconc}) and for the independent Cauchy case (Lemma~\ref{LemCauchyquasi}).

\begin{Lem}
\label{lemgeodesicharmonicquasiconc}
For any~$\Sigma_a,\Sigma_b\in\mathcal{P}_k$ and~$t\in [0,1]$, let $\tilde{\Sigma}_t:=\Sigma_a^{1/2} \big( \Sigma_a^{-1/2} \Sigma_b 
\linebreak
\Sigma_a^{-1/2} \big)^t \Sigma_a^{1/2}$
and~$\Sigma_t^*:=((1-t) \Sigma_a^{-1} + t \Sigma_b^{-1})^{-1}$. Then 
(i) $\lambda_1(\Sigma_t^*) \leq \lambda_1(\tilde{\Sigma}_t) \leq \max( \lambda_1(\Sigma_a) , \lambda_1(\Sigma_b))$ 
and
(ii) $\lambda_k(\tilde{\Sigma}_t) \geq \lambda_k(\Sigma_t^*) \geq \min( \lambda_k(\Sigma_a) , \lambda_k(\Sigma_b))$.
\end{Lem}

{\sc Proof of Lemma~\ref{lemgeodesicharmonicquasiconc}}.
(i) With the usual order on positive semidefinite matrices ($A\leq B$ iff~$B-A$ is positive semidefinite), the (weighted) 
harmonic-geometric-arithmetic inequality (see, e.g., Lemma~2.1(vii) in \citealp{LawLim13})
\begin{equation}
\label{GA} 
\Sigma_t^*
 \leq 
\tilde{\Sigma}_t
\leq 
\Sigma_t
:=
(1-t) \Sigma_a + t \Sigma_b
\end{equation}
holds for any~$t\in[0,1]$. This implies that
\begin{equation}
\label{specineq}
\lambda_1(\Sigma_t^*)
\leq
\lambda_1(\tilde{\Sigma}_t)
\leq 
\lambda_1(\Sigma_t)
. 
\end{equation}
Indeed, if, e.g., the second inequality in~(\ref{specineq}) does not hold (the argument for the first inequality is strictly the same), then, denoting as~$\tilde{v}_{1t}$ an arbitrary eigenvector associated with~$\lambda_1(\tilde{\Sigma}_t)$, we have~$\tilde{v}_{1t}' \tilde{\Sigma}_t \tilde{v}_{1t}
=
\lambda_1(\tilde{\Sigma}_t)
>
\lambda_1(\Sigma_t)
\geq
\tilde{v}_{1t}' \Sigma_t \tilde{v}_{1t}
$,
which contradicts~(\ref{GA}). Hence,~(\ref{specineq}) holds and provides
\begin{eqnarray*}
\lambda_1(\Sigma_t^*)
\leq 
\lambda_1(\tilde{\Sigma}_t)
&\!\!\! \leq\!\!\! & 
\max_{u\in\mathcal{S}^{d-1}} 
u'((1-t) \Sigma_a + t \Sigma_b)u
\\[2mm]
&\!\!\! \leq\!\!\! & 
(1-t)
\max_{u\in\mathcal{S}^{d-1}} u'  \Sigma_a u
+t
\max_{u\in\mathcal{S}^{d-1}} u'  \Sigma_b u
\\[2mm]
&\!\!\! =\!\!\! & 
(1-t)
\lambda_1(\Sigma_a)
+t
\lambda_1(\Sigma_b)
\\[2mm]
&\!\!\! \leq\!\!\! & 
\max(\lambda_1(\Sigma_a),\lambda_1(\Sigma_b))
,
\end{eqnarray*}
as was to be showed. (ii) Proceeding in a similar way as above, it is readily showed that~(\ref{GA}) implies that $\lambda_k(\tilde{\Sigma}_t)\geq \lambda_k(\Sigma_t^*)$. 
Using this, we obtain
\begin{eqnarray*}
\lefteqn{
\lambda_k(\tilde{\Sigma}_t)
\geq 
\lambda_k(\Sigma_t^*)
= 
\lambda_1^{-1}((\Sigma_t^*)^{-1})
}
\\[3mm]
& &
\hspace{-3mm} 
=
\Big(
\max_{u\in\mathcal{S}^{d-1}} 
u' ((1-t) \Sigma_a^{-1} + t \Sigma_b^{-1}) u
\Big)^{-1}
\geq 
\Big( 
(1-t)
\lambda_1(\Sigma_a^{-1})
+t
\lambda_1(\Sigma_b^{-1})
\Big)^{-1}
\\[2mm]
& & 
\hspace{-3mm} 
= 
\Big( 
(1-t)
\lambda_k^{-1}(\Sigma_a)
+t
\lambda_k^{-1}(\Sigma_b)
\Big)^{-1}
\geq 
\min(\lambda_k(\Sigma_a),\lambda_k(\Sigma_b))
,
\end{eqnarray*}
since any weighted harmonic mean of two real numbers is a convex linear combination of these.
\cqfd 
\vspace{3mm}


\begin{Lem}
\label{LemCauchyquasi}
For any~$\Sigma_a,\Sigma_b\in\mathcal{P}_k$ and~$t\in [0,1]$, let $\tilde{\Sigma}_t:=\Sigma_a^{1/2} \big( \Sigma_a^{-1/2} \Sigma_b 
\linebreak
\Sigma_a^{-1/2} \big)^t \Sigma_a^{1/2}$
and~$\Sigma_t^*:=((1-t) \Sigma_a^{-1} + t \Sigma_b^{-1})^{-1}$. Then, 
\begin{equation}
\label{cc1}	
\max({\rm diag}(\Sigma_t^*))\leq \max({\rm diag}(\tilde{\Sigma}_t))\leq \max( \max({\rm diag}(\Sigma_a)),\max({\rm diag}(\Sigma_b)))
\end{equation}
and 
\begin{equation}
\label{cc2}	
\max_s s' \tilde{\Sigma}_t^{-1} s
\leq
\max_s s' (\Sigma_t^*)^{-1} s
\leq
\max\big(
	\max_s s' \Sigma_a^{-1} s
,
 \max_s s' \Sigma_b^{-1} s
\big)
\end{equation}
$($where~$\max_s$ is the maximum over~$s=(s_1,\ldots,s_{k})\in\{-1,1\}^{k})$, so that both the mappings~$\Sigma\mapsto \max({\rm diag}(\Sigma))$ and $\Sigma\mapsto \max_s s' \Sigma^{-1} s$ are geodesic and harmonic quasi-convex.
\end{Lem}

{\sc Proof of Lemma~\ref{LemCauchyquasi}.}
The result in~(\ref{cc1}) readily follows from the fact that the weighted harmonic-geometric-arithmetic inequality~$\Sigma_t^*\leq\tilde{\Sigma}_t\leq (1-t) \Sigma_a + t \Sigma_b$ yields $(\Sigma_t^*)_{\ell\ell}\leq (\tilde{\Sigma}_t)_{\ell\ell}\leq (1-t) (\Sigma_a)_{\ell\ell} + t (\Sigma_b)_{\ell\ell} \leq \max( (\Sigma_a)_{\ell\ell},(\Sigma_b)_{\ell\ell})$ for any~$\ell=1,\ldots,k$. Turning to~(\ref{cc2}), the harmonic-geometric inequality implies that~$\tilde{\Sigma}_t^{-1}\leq (\Sigma_t^*)^{-1}$, which readily yields
\linebreak
$
\max_s s' \tilde{\Sigma}_t^{-1} s
\leq
\max_s s' (\Sigma_t^*)^{-1} s
$.
Consequently, it only remains to prove the second inequality in~(\ref{cc2}). To do so, choose an arbitrary~$s_*$ such that
$\max_s s' (\Sigma_t^*)^{-1} s
= s_*' (\Sigma_t^*)^{-1} s_*$. Then 
\begin{eqnarray}
\lefteqn{
\hspace{-13mm}
\max_s s' (\Sigma_t^*)^{-1} s
=
s_*' (\Sigma_t^*)^{-1} s_*
=
s_*' \big( (1-t) \Sigma_a^{-1}  + t \Sigma_b^{-1}\big) s_*
}
\nonumber
\\[2mm]
& & \hspace{3mm} 
\leq
(1-t) \max_s s' \Sigma_a^{-1} s  + t \max_s s' \Sigma_b^{-1} s
\nonumber
\\[2mm]
& & \hspace{3mm} 
\leq
\max\big( \max_s s' \Sigma_a^{-1} s
,
\max_s s' \Sigma_b^{-1} s
\big)
,
\label{thisforonehalf}
\end{eqnarray}
which establishes the result.  
\cqfd
\vspace{3mm}

 
We can now prove Theorems~\ref{propgeodesicquasiconc} and~\ref{propharmonicquasiconc}. 
\vspace{3mm}

{\sc Proof of Theorem~\ref{propgeodesicquasiconc}}. 
(i) We start by considering the case where~$P$ is an elliptical probability measure over~$\R^k$ with location~$\theta_0$ and scatter~$\Sigma_0$, where we first prove the result for~$\theta_0=0$ and~$\Sigma_0=I_k$. Then we have
\begin{eqnarray}
H\!D^{\rm sc}_{P,T}(\tilde{\Sigma}_t)
&\!\!=\!\!&
\inf_{u\in\mathcal{S}^{k-1}} 
\min
\big( 
P\big[ |u'X| \leq 
{\textstyle{\sqrt{u'\tilde{\Sigma}_t u}}}
\, \big]
,
P\big[ |u'X| \geq 
{\textstyle{\sqrt{u'\tilde{\Sigma}_t u}}}
\, \big]
\big)
\nonumber
\\[2mm]
&\!\!=\!\!&
\min
\Big( 
\inf_{u\in\mathcal{S}^{k-1}} 
P[ |X_1| \leq 
{\textstyle{\sqrt{u'\tilde{\Sigma}_t u}}}
 ]
,
\inf_{u\in\mathcal{S}^{k-1}}
P[ |X_1| \geq 
{\textstyle{\sqrt{u'\tilde{\Sigma}_t u}}}
 ]
\Big)
\nonumber
\\[3mm]
&\!\!=\!\!&
\min( 
P[ |X_1| \leq \lambda_k^{1/2}(\tilde{\Sigma}_t) ]
,
P[ |X_1| \geq \lambda_1^{1/2}(\tilde{\Sigma}_t) ]
)
.
\label{dgsk}
\end{eqnarray}
Since Lemma~\ref{lemgeodesicharmonicquasiconc} entails that
\begin{eqnarray*}
\lefteqn{
P[ |X_1| \leq \lambda_k^{1/2}(\tilde{\Sigma}_t) ]
 \geq 
P[ |X_1| \leq  \min( \lambda_k^{1/2}(\Sigma_a) , \lambda_k^{1/2}(\Sigma_b) ) ]
}
\\[1mm]
& & 
\hspace{3mm} 
 = 
\min( P[ |X_1| \leq \lambda_k^{1/2}(\Sigma_a)] , P[ |X_1| \leq \lambda_k^{1/2}(\Sigma_b)] )
\\[2mm]
& & 
\hspace{3mm} 
 \geq 
\min(
H\!D^{\rm sc}_{P,T}(\Sigma_a)
,
H\!D^{\rm sc}_{P,T}(\Sigma_b)
)
\end{eqnarray*}
and
\begin{eqnarray*}
\lefteqn{
P[ |X_1| \geq \lambda_1^{1/2}(\tilde{\Sigma}_t) ]
 \geq 
P[ |X_1| \geq  \max( \lambda_1^{1/2}(\Sigma_a) , \lambda_1^{1/2}(\Sigma_b) ) ]
}
\\[1mm]
& &
\hspace{3mm} 
=\min( P[ |X_1| \geq \lambda_1^{1/2}(\Sigma_a)] , P[ |X_1| \geq \lambda_1^{1/2}(\Sigma_b)] )
\\[2mm]
& &
\hspace{3mm} 
 \geq 
\min(
H\!D^{\rm sc}_{P,T}(\Sigma_a)
,
H\!D^{\rm sc}_{P,T}(\Sigma_b)
)
,
\end{eqnarray*}
the result for~$\theta_0=0$ and~$\Sigma_0=I_k$ follows from~(\ref{dgsk}). 

We now prove the result in the elliptical case with arbitrary values of~$\theta_0$ and~$\Sigma_0$. To this end, let~$A=\Sigma_0^{-1/2}$ and note that the square roots (in~$\mathcal{P}_k$) of $\Upsilon_a:=A \Sigma_a A'$ and $\Upsilon_b:=A \Sigma_b A'$ are of the form~$\Upsilon_a^{1/2}=A \Sigma_a^{1/2}O_a$ and~$\Upsilon_b^{1/2}=A \Sigma_b^{1/2}O_b$, for some $k\times k$ orthogonal matrices~$O_a,O_b$. Consequently, 
\begin{eqnarray*}
\lefteqn{
\tilde{\Upsilon}_t
:=
A \tilde{\Sigma}_t A'
=
A
\Sigma_a^{1/2} 
( \Sigma_a^{-1/2} \Sigma_b \Sigma_a^{-1/2} )^t 
\Sigma_a^{1/2}
A'
}
\\[2mm]
& &
\hspace{3mm} 
=
\Upsilon_a^{1/2} O_a'
( O_a\Upsilon_a^{-1/2}
\Upsilon_b 
\Upsilon_a^{-1/2} O_a' )^t 
O_a \Upsilon_a^{1/2}
=
\Upsilon_a^{1/2} 
( \Upsilon_a^{-1/2}
\Upsilon_b
\Upsilon_a^{-1/2} )^t 
\Upsilon_a^{1/2}
\end{eqnarray*}
describes a geodesic path from~$\Upsilon_a$ to~$\Upsilon_b$. 
Since the result holds at~$P_0=P_{A,-A\theta_0}$ (where the notation~$P_{A,b}$ was defined on page~\pageref{pagedefPAb}, affine invariance then entails that
\begin{eqnarray*}
H\!D^{\rm sc}_{P,T}(\tilde{\Sigma}_t)=H\!D^{\rm sc}_{P_0,T}(\tilde{\Upsilon}_t)
&\!\!\!\geq \!\!\!&
 \min(H\!D^{\rm sc}_{P_0,T}(\Upsilon_a),H\!D^{\rm sc}_{P_0,T}(\Upsilon_b))
\\[2mm]
&\!\!\!=\!\!\!&
\min(H\!D^{\rm sc}_{P,T}(\Sigma_a),H\!D^{\rm sc}_{P,T}(\Sigma_b))
,
\end{eqnarray*}
as was to be showed. 

We now turn to the case where the probability measure~$P$ over~$\R^k$ has independent Cauchy marginals. Fix~$\Sigma_a,\Sigma_b\in\mathcal{P}_k$ and consider the geodesic path~$\tilde{\Sigma}_t$, $t\in[0,1]$, from~$\Sigma_a$ to~$\Sigma_b$. Recall that 
$$
H\!D_{P,T}^{\rm sc}(\Sigma)
=
2
\min
\big(\Psi\big( 1/ \max_s \sqrt{s' \Sigma^{-1} s}\, \big) - \textstyle{\frac{1}{2}}
,
1 - \Psi\big( \sqrt{\max({\rm diag}(\Sigma))}\, 
\big) 
\big),
$$
where~$\Psi$ stands for the Cauchy cumulative distribution function; see~(\ref{Cauchytheo}). Lemma~\ref{LemCauchyquasi} readily entails that 
\begin{eqnarray} 
\lefteqn{
\hspace{4mm} 
2 - 2 \Psi\big( 
{\textstyle{\sqrt{\max({\rm diag}(\tilde{\Sigma}_t))}}}
\, \big) 
}
\nonumber
\\[2mm]
& & 
\hspace{10mm} 
\geq 
\min
\big(
2 - 2 \Psi\big( \sqrt{\max({\rm diag}(\Sigma_a))}\, \big)
,
2 - 2 \Psi\big( \sqrt{\max({\rm diag}(\Sigma_b))}\, \big)
\big)
\nonumber
\\[2mm]
& & 
\hspace{10mm} 
\geq 
\min(
H\!D^{\rm sc}_{P,T}(\Sigma_a)
,
H\!D^{\rm sc}_{P,T}(\Sigma_b)
)
.
\label{t1}
\end{eqnarray}
Lemma~\ref{LemCauchyquasi} also provides~$\max_s \textstyle{s' \tilde{\Sigma}_t^{-1} s} \leq \max( \max_s \textstyle{s' \Sigma_a^{-1} s},\max_s \textstyle{s' \Sigma_b^{-1} s})$, which rewrites 
$$
1/\max_s \textstyle{(s' \tilde{\Sigma}_t^{-1} s)^{1/2}} \geq \min( 1/\max_s \textstyle{(s' \Sigma_a^{-1} s)^{1/2}} , 1/\max_s \textstyle{(s' \Sigma_b^{-1} s)^{1/2}})
.
$$
This implies that
\begin{eqnarray} 
\lefteqn{
2 \Psi\big(  1/ \max_s \textstyle{(s' \tilde{\Sigma}_t^{-1} s)^{1/2}}\, \big) - 1
}
\nonumber
\\[2mm]
& & 
\hspace{8mm} 
\geq 
\min
\big(
2 \Psi\big(  1/ \max_s \textstyle{(s' \Sigma_a^{-1} s)^{1/2}}\, \big) - 1
,
2 \Psi\big(  1/ \max_s \textstyle{(s' \Sigma_b^{-1} s)^{1/2}}\, \big) - 1
\big)
\nonumber
\\[2mm]
& & 
\hspace{8mm}  
\geq 
\min(
H\!D^{\rm sc}_{P,T}(\Sigma_a)
,
H\!D^{\rm sc}_{P,T}(\Sigma_b)
)
.
\label{t2}
\end{eqnarray}
From~(\ref{t1})-(\ref{t2}), it readily follows that
$H\!D^{\rm sc}_{P,T}(\tilde{\Sigma}_t)
\geq \min(
H\!D^{\rm sc}_{P,T}(\Sigma_a)
,
\linebreak
H\!D^{\rm sc}_{P,T}(\Sigma_b)
)
$, which concludes the proof of Part~(i).

(ii) Let then~$P$ be an arbitrary probability measure over~$\R^k$ satisfying Part~(i) of the result. For any~$\Sigma_a,\Sigma_b\in R^{\rm sc}_{P,T}(\alpha)$, we then have $H\!D^{\rm sc}_{P,T}(\tilde{\Sigma}_t)\geq \min(H\!D^{\rm sc}_{P,T}(\Sigma_a),H\!D^{\rm sc}_{P,T}(\Sigma_b))\geq \alpha$, so that~$\tilde{\Sigma}_t\in R^{\rm sc}_{P,T}(\alpha)$.
\cqfd
\vspace{3mm}
 
 
{\sc Proof of Theorem~\ref{propharmonicquasiconc}}. In view of the remark given right before the statement of the theorem, it is sufficient to show that, at the probability measures~$P$ considered, $\Sigma\mapsto H\!D^{\rm sc}_{P,T}(\Sigma)$ is harmonic quasi-concave. The proof is then entirely similar to the proof of Theorem~\ref{propgeodesicquasiconc}. In the elliptical case, the affine-invariance argument is based on the identity~$
\Upsilon_t^*
:=
A \Sigma_t^* A'
=
((1-t)\Upsilon_a^{-1}+t\Upsilon_b^{-1})^{-1}
$, with~$\Upsilon_a:=A \Sigma_a A'$ and~$\Upsilon_b:=A \Sigma_b A'$.
\cqfd
\vspace{3mm}


\subsection{Proofs from Section~7}
\label{appshapeec}


\
\vspace{3mm}

{\sc Proof of Theorem~\ref{maxsigma2shape}.}
We may restrict to the
\vspace{-1mm}
  case where~$H\!D^{{\rm sh},S}_{P,T}(\VS)>\alpha_{P,T}$ (indeed, the 
  \vspace{-.6mm}
assumptions ensure that~$H\!D^{{\rm sh},S}_{P,T}(\VS)\geq \alpha_{P,T}$ and that the result holds if~$H\!D^{{\rm sh},S}_{P,T}(\VS)= \alpha_{P,T}$).
For any~$\alpha>\alpha_{P,T}$, consider then~$I_\alpha:=I_{\alpha,P,T}(\VS):=\{\sigma^2\in\R^+_0 : \sigma^2 \VS \in R^{\rm sc}_{P,T}(\alpha) \}=\{\sigma^2\in\R^+_0 : H\!D^{\rm sc}_{P,T}(\sigma^2 \VS)\geq \alpha\}$. The convexity of~$R^{\rm sc}_{P,T}(\alpha)$ (Theorem~\ref{Propquasiconcavity}(ii)) implies that~$I_\alpha$ is an interval. Since~$\alpha>\alpha_{P,T}$, Theorem~\ref{geodboundedness} shows that~$R^{\rm sc}_{P,T}(\alpha)$ is $g$-bounded, which implies there exist~$\eta_\alpha>0$ and~$M_\alpha>\eta_\alpha$ such that~$I_\alpha\subset[\eta_\alpha,M_\alpha]$. Since, moreover, Theorem~\ref{Propcontinuity} implies that~$\sigma^2 \mapsto H\!D^{\rm sc}_{P,T}(\sigma^2 \VS)$ is upper semicontinuous, $I_\alpha$ is also closed, hence (still for~$\alpha>\alpha_{P,T}$) compact.  

Now, fix~$\delta>0$ such that~$H\!D^{{\rm sh},S}_{P,T}(\VS)-\delta>\alpha_{P,T}$. For any~$n$, pick then~$\sigma^2_n$ in the (non-empty) interval~$I_{H\!D^{{\rm sh},S}_{P,T}(\VS)-(\delta/n)}$. The resulting sequence~$(\sigma^2_n)$ is in the compact set~$I_{H\!D^{{\rm sh},S}_{P,T}(\VS)-\delta}$, hence admits a subsequence~$(\sigma^2_{n_\ell})$ converging in~$\R^+_0$, to~$\sigma^2_{\VS}$, say. Fix then an arbitrary~$\varepsilon\in (0,\delta)$. For~$\ell$ large enough, all~$\sigma^2_{n_\ell}$ belong to the closed set~$I_{H\!D^{{\rm sh},S}_{P,T}(\VS)-\varepsilon}$, so that~$\sigma^2_{\VS}$ also belongs to~$I_{H\!D^{{\rm sh},S}_{P,T}(\VS)-\varepsilon}$. This shows that~$H\!D^{{\rm sh},S}_{P,T}(\VS)-\varepsilon\leq H\!D^{\rm sc}_{P,T}(\sigma^2_{\VS} \VS)\leq H\!D^{{\rm sh},S}_{P,T}(\VS)$. Since~$\varepsilon$ can be taken arbitrarily small, the result is proved.
\cqfd
\vspace{3mm}
   

{\sc Proof of Theorem~\ref{thaffinvshape}.} 
From Theorem~\ref{thaffinv}, we readily obtain
\begin{eqnarray*}
\lefteqn{
H\!D^{{\rm sh},S}_{P_{A,b},T}(A V\! A' / S(A V\! A'))
=
\sup_{\sigma^2>0}
H\!D^{\rm sc}_{P_{A,b},T}(\sigma^2 A V\! A'/ S(A V\! A'))
}
\\[2mm]
& &
\hspace{10mm}
=
\sup_{\sigma^2>0}
H\!D^{\rm sc}_{P_{A,b},T}(\sigma^2 A V\! A')
=
\sup_{\sigma^2>0}
H\!D^{\rm sc}_{P,T}(\sigma^2 V)
=
H\!D^{{\rm sh},S}_{P,T}(V)
,
\end{eqnarray*}
which establishes the result.  
\cqfd
\vspace{3mm}


{\sc Proof of Theorem~\ref{Consistencyshape}.} 
Consider arbitrary probability measures~$P,Q$ on~$\R^k$. Fix~$V\in\mathcal{P}_k^{S}$ and assume (without loss of generality) that~$H\!D^{{\rm sh},S}_{P,T}(V) \leq  H\!D^{{\rm sh},S}_{Q,T}(V)$. Then, for any~$\varepsilon>0$, there exists~$\sigma^2_{\varepsilon}>0$ such that~$H\!D^{{\rm sh},S}_{Q,T}(V) \leq H\!D^{\rm sc}_{Q,T}(\sigma^2_{\varepsilon} V)+\varepsilon$, so that
\begin{eqnarray*}
\lefteqn{
\hspace{-1mm} 
|H\!D^{{\rm sh},S}_{P,T}(V) - H\!D^{{\rm sh},S}_{Q,T}(V)|
=
H\!D^{{\rm sh},S}_{Q,T}(V) - H\!D^{{\rm sh},S}_{P,T}(V)
}
\\[2mm]
& & 
\hspace{-0mm} 
\leq
H\!D^{\rm sc}_{Q,T}(\sigma^2_{\varepsilon} V)
+\varepsilon
- 
H\!D^{\rm sc}_{P,T}(\sigma^2_{\varepsilon} V)
\leq
\sup_{\Sigma\in\mathcal{P}_k} 
| 
H\!D_{P,T}^{\rm sc}(\Sigma)
-
H\!D_{Q,T}^{\rm sc}(\Sigma)
|
+
\varepsilon
.
\end{eqnarray*}
Since this holds for any~$\varepsilon>0$ and since~$V$ is arbitrary, we have that
$$
\sup_{V\in\mathcal{P}_k^{S}} 
|H\!D^{{\rm sh},S}_{P,T}(V) - H\!D^{{\rm sh},S}_{Q,T}(V)|
\leq
\sup_{\Sigma\in\mathcal{P}_k} 
| 
H\!D^{\rm sc}_{Q,T}(\Sigma)
-
H\!D_{P,T}^{\rm sc}(\Sigma)
|
.
$$
The result then follows from Theorem~\ref{Consistency}.
\cqfd
\vspace{3mm}


{\sc Proof of Theorem~\ref{Propcontinuityshape}.} (i) Fix a shape matrix~$V_0\in R^{{\rm sh},S}_{P,T}(\alpha_{P,T})$ and assume, ad absurdum, that there exists a sequence~$(V_n)$ in $\mathcal{P}_k^{S}$ that $g$-converges (resp., $F$-converges) to~$V_0$ and such that 
$$
\lim\sup_{n\to\infty} H\!D^{{\rm sh},S}_{P,T}(V_n)>H\!D^{{\rm sh},S}_{P,T}(V_0).
$$
 Extracting a subsequence if necessary, we can fix~$\varepsilon>0$ small enough to have $H\!D^{{\rm sh},S}_{P,T}(V_0)+\varepsilon<H\!D^{{\rm sh},S}_{P,T}(V_n)$ for any~$n$. Fix then, for any~$n$, $\sigma_n^2>0$ such that $H\!D^{\rm sc}_{P,T}(\sigma_n^2V_n)>H\!D^{{\rm sh},S}_{P,T}(V_n)-\varepsilon/2$, which yields 
\begin{equation}
\label{toshusc}
	H\!D^{\rm sc}_{P,T}(\sigma_n^2V_n)>H\!D^{{\rm sh},S}_{P,T}(V_0)+\varepsilon/2
	\geq \alpha_{P,T}+\varepsilon/2
	.
	\end{equation}
 Now, we can assume without loss of generality that $V_n$ belongs to a neighborhood of~$V_0$ that is $g$-compact in~$\mathcal{P}_k^{S}(P)$. Since~(\ref{toshusc}) implies that~$\sigma_n^2V_n$ belongs, for any~$n$, to the $g$-bounded (Theorem~\ref{geodboundedness}) scatter depth region $R^{\rm sc}_{P,T}(\alpha_{P,T}+\varepsilon/2)$, the sequence~$(\sigma^2_n)$ then stays away from~$0$ and~$\infty$ (that is, the $\sigma_n^2$'s belong to a common compact set of~$\R^+_0$). Consequently, there exists a subsequence~$(\sigma^2_{n_\ell})$ such that~$(\sigma_{n_\ell}^2 V_{n_\ell})$ $g$-converges (resp., $F$-converges)  to~$\sigma_0^2 V_0$, say. In view of~(\ref{toshusc}), we therefore found~$\varepsilon>0$ such that, for any~$\ell$, 
$$
H\!D^{\rm sc}_{P,T}(\sigma_{n_\ell}^2 V_{n_\ell})
> 
H\!D^{{\rm sh},S}_{P,T}(\sigma^2_0 V_0)+\varepsilon/2
,
$$
where~$(\sigma_{n_\ell}^2 V_{n_\ell})$ $g$-converges (resp., $F$-converges) to~$\sigma_0^2 V_0$, which contradicts the scatter depth upper semicontinuity result in Theorem~\ref{geodcontinuity} (resp., in Theorem~\ref{Propcontinuity}).
	(ii) The result follows from the fact that~$R^{\rm sh}_{P,T}(\alpha)$ is the inverse image of~$[\alpha,+\infty)$ by the upper $F$- and $g$-semicontinuous function~$V\mapsto H\!D^{{\rm sh},S}_{P,T}(V)$. 
	(iii) Since the supremum of lower semicontinuous functions is a lower semicontinuous function, Theorems~\ref{Propcontinuity} and~\ref{geodcontinuity}(iii) yield that~$V\mapsto H\!D^{{\rm sh},S}_{P,T}(V)$ is lower-semi 
	\vspace{-1mm}
 continuous. The result then follows from Part~(i) and the fact that the smoothness of~$P$ at~$T_P$ implies that~$R^{{\rm sh},S}_{P,T}(\alpha_{P,T})=R^{{\rm sh},S}_{P,T}(0)=\mathcal{P}_k^{S}$. 
\cqfd 
\vspace{3mm}


The proof of Theorem~\ref{boundednessshape} requires the following lemma.

\begin{Lem}
\label{lemshapeboundedness}
Let~$S$ be a scale functional, that is a mapping from~$\mathcal{P}_k$ to~$\R^+_0$ that satisfies the properties~(i)-(iii) on page~\pageref{pagescalefunctcond}. Then, $\lambda_k(V)\leq 1\leq \lambda_1(V)$ for any~$V\in \mathcal{P}_k^{S}$.
\end{Lem}

{\sc Proof of Lemma~\ref{lemshapeboundedness}.}
(a) Factorize~$V$ into~$V
=
O\, {\rm diag}(\lambda_1(V),\ldots,
\linebreak
\lambda_k(V)) O'
$, where~$O$ is a $k\times k$ orthogonal matrix. Then,  
since
$
\lambda_k(V) I_k
= 
O\, {\rm diag}(\lambda_k(V),\ldots,\lambda_k(V)) O'
\leq
V
\leq 
 O\, {\rm diag}(\lambda_1(V),\ldots,\lambda_1(V)) O'
=
\lambda_1(V) I_k
$
(where $A\leq B$ still means that~$B-A$ is positive semidefinite), the properties of a scale functional yield~$\lambda_k(V)
=S(\lambda_k(V)I_k)\leq S(V)\leq S(\lambda_1(V) I_k)
=\lambda_1(V)$.  
\cqfd
\vspace{3mm}

  {\sc Proof of Theorem~\ref{boundednessshape}.} 
  We start with the proof of the result for~$s_{P,T}<1/2$ and $g$-boundedness. We fix~$\varepsilon>0$ and intend to prove that~$R^{{\rm sh},S}_{P,T}(s_{P,T}+\varepsilon)$ is $g$-bounded by showing that, for~$r>0$ large enough, it is included in the $g$-ball~$B_g(I_k,r)$. To do so, fix a shape matrix~$V(\in\mathcal{P}_k^{S})$ that does not belong to~$B_g(I_k,r)$ ($r$ is to be chosen later). In view of~(\ref{distmaj}), we then have 
(i) $\lambda_1(V)  > \exp(r/\sqrt{k})$
or
(ii) $\lambda_k(V)  < \exp(-r/\sqrt{k})$
(or both). 

We start with case~(i). Fix (so far, arbitrarily)~$\sigma^2_0>0$. Then for any~$\sigma^2\in(0,\sigma^2_0]$, Lemma~\ref{lemshapeboundedness} entails that (denoting by~$v_k(V)$ an arbitrary unit vector associated with~$\lambda_k(V)$)
\begin{eqnarray}
H\!D^{\rm sc}_{P,T}(\sigma^2 V)
&\!\! \leq\!\! &
\inf_{u\in\mathcal{S}^{k-1}}
P\big[ |u'(X-T_P)| \leq 
\sigma 
{\textstyle{\sqrt{u'V u}}}
 \,\big]
\nonumber
\\[2mm]
&\!\! \leq\!\! &
P\big[ |v_k'(V)(X-T_P)| \leq \sigma \lambda_k^{1/2}(V)\big]
\nonumber
\\[2mm]
&\!\! \leq\!\! &
P[ |v_k'(V)(X-T_P)| \leq \sigma_0]
\leq s_P^{\{T_P\}}(\sigma_0)
\label{bshap1}
,
\end{eqnarray}
where we used the notation~$s_P^K(\cdot)$ introduced in Lemma~\ref{LemConsistency2}.
By using this lemma, pick then~$\sigma_0^2>0$ such that~$s_P^{\{T_P\}}(\sigma_0)< s_P^{\{T_P\}}+(\varepsilon/2)=s_{P,T}+(\varepsilon/2)$. Denoting by~$v_1(V)$ an arbitrary unit vector associated with~$\lambda_1(V)$, we then have that, for any~$\sigma^2\in [\sigma^2_0,\infty)$,
\begin{eqnarray}
H\!D^{\rm sc}_{P,T}(\sigma^2 V)
&\!\! \leq\!\! &
\inf_{u\in\mathcal{S}^{k-1}}
P\big[ |u'(X-T_P)| \geq \sigma 
{\textstyle{\sqrt{u'V u}}}\, \big]
\nonumber
\\[1mm]
&\!\! \leq\!\! &
P\big[ |v'_1(V)(X-T_P)| \geq  \sigma \lambda_1^{1/2}(V)\,\big]
\nonumber
\\[2mm] 
&\!\! \leq\!\! &
P\big[ |v'_1(V)(X-T_P)| \geq \sigma_0 \exp(r/2\sqrt{k})\big]
\nonumber
\\[2mm] 
&\!\! \leq\!\! &
P\big[ \| X-T_P \| \geq \sigma_0 \exp(r/2\sqrt{k})\big]
,
\label{bshap2}
\end{eqnarray}
which, for~$r$ large enough, can be made smaller than $\varepsilon/2$ (hence, smaller than~$s_{P,T}+(\varepsilon/2)$). For~$r$ large enough, thus, (\ref{bshap1})-(\ref{bshap2}) guarantee that $H\!D^{{\rm sh},S}_{P,T}(V)=\sup_{\sigma^2>0} H\!D^{\rm sc}_{P,T}(\sigma^2 V)< s_{P,T}+\varepsilon$, as was to be showed. 

We then turn to case~(ii). By picking~$\sigma_0$ large enough, we have that, for any~$\sigma^2\in [\sigma^2_0,\infty)$, 
\begin{eqnarray}
\hspace{-17mm} 
H\!D^{\rm sc}_{P,T}(\sigma^2 V)
&\!\! \leq\!\! &
\inf_{u\in\mathcal{S}^{k-1}}
P\big[ |u'(X-T_P)| \geq \sigma 
{\textstyle{\sqrt{u'V u}}}\, \big]
\nonumber
\\[1mm]
&\!\! \leq\!\! &
P\big[ |v'_1(V)(X-T_P)| \geq  \sigma \lambda_1^{1/2}(V)\,\big]
\nonumber
\\[2mm]
&\!\! \leq\!\! &
P\big[ |v'_1(V)(X-T_P)| \geq \sigma_0\big]
\nonumber
\\[2mm]
&\!\! \leq\!\! &
P\big[ \| X-T_P \| \geq \sigma_0 \big]
< \varepsilon/2,
\label{bshap3}
\end{eqnarray}
where we used Lemma~\ref{lemshapeboundedness}. For any~$\sigma^2\in(0,\sigma^2_0]$, we then have
\begin{eqnarray}
\hspace{-0mm} 
H\!D^{\rm sc}_{P,T}(\sigma^2 V)
&\!\! \leq\!\! &
\inf_{u\in\mathcal{S}^{k-1}}
P\big[ |u'(X-T_P)| \leq 
\sigma 
{\textstyle{\sqrt{u'V u}}}
 \,\big]
\nonumber
\\[1mm]
&\!\! \leq\!\! &
P\big[ |v_k'(V)(X-T_P)| \leq \sigma \lambda_k^{1/2}(V)\big]
\nonumber
\\[2mm]
&\!\! \leq\!\! &
P[ |v_k'(V)(X-T_P)| \leq \sigma_0\exp(-r/2\sqrt{k})]
\nonumber
\\[2mm]
&\!\! \leq\!\! &
 s_P^{\{T_P\}}(\sigma_0 \exp(-r/2\sqrt{k})
\label{bshap4}
< s_{P,T}+(\varepsilon/2)
,
\end{eqnarray}
for~$r$ large enough. Thus, for~$r$ large enough, (\ref{bshap3})-(\ref{bshap4}) still yield that
$
H\!D^{{\rm sh},S}_{P,T}(V)=\sup_{\sigma^2>0} H\!D^{\rm sc}_{P,T}(\sigma^2 V) < s_{P,T}+\varepsilon,
$
 as was to be showed. We thus conclude that, for~$\alpha>s_{P,T}$, $R^{{\rm sh},S}_{P,T}(\alpha)$ is $g$-bounded (its $g$-compacity then follows from the same argument as in the proof of Theorem~\ref{geodboundedness}). 

The proof for $F$-boundedness (still for~$s_P<1/2$) follows along the same lines and is actually simpler since only one of both cases~(i)-(ii) above is to be considered. Recall indeed that, as seen in the proof of Theorem~\ref{boundedness}, $V\notin B_F(I_k,r)$ implies that~$\lambda_1(V)>(r-1)/k^{1/2}$, so that the same reasoning as in case~(i) above allows to show that for any~$\varepsilon>0$, there exists~$r=r(\varepsilon)$ such that~$V\notin B(I_k,r)$ implies $H\!D^{{\rm sh},S}_{P,T}(V)< s_{P,T}+\varepsilon$. This establishes that~$R^{{\rm sh},S}_{P,T}(\alpha)$ is $F$-bounded for~$\alpha>s_{P,T}$.

Finally, if~$s_{P,T}\geq 1/2$, then, with~$u_0\in\mathcal{S}^{k-1}$ such that~$P[|u_0' (X-T_P)|=0]=s_{P,T}$ (existence is guaranteed in Lemma~\ref{LemConsistency2}(ii), with~$K=\{T_P\}$), we have
$
H\!D^{\rm sc}_{P,T}(\sigma^2 V) 
\leq
P[ |u_0'(X-T_P)| \geq \sigma 
{\textstyle{\sqrt{u_0' V u_0}}}\, ]
\leq 
P[ |u_0'(X-T_P)| >0 ]
=
1-
s_{P,T}
$
for any~$\sigma^2>0$, so that~$H\!D^{{\rm sh},S}_{P,T}(V)\leq 1-s_{P,T}$. Therefore, $R^{{\rm sh},S}_{P,T}(\alpha)$ is empty for any~$\alpha>1-s_{P,T}=\alpha_{P,T}$.
\cqfd
\vspace{3mm}
 

{\sc Proof of Theorem~\ref{geodmaxdepthshape}.} The proof follows along the exact same lines as that of~Theorem~\ref{geodmaxdepth}, hence is not reported here. 
\cqfd
\vspace{3mm}


{\sc Proof of Theorem~\ref{Fishconstshape}.} 
(i) 
For any~$V\in\mathcal{P}_k^{S}$, Theorem~\ref{Fishconst}(i) readily implies that $H\!D^{{\rm sh},S}_{P,T}(V)=
\sup_{\sigma^2>0} H\!D^{\rm sc}_{P,T}(\sigma^2 V)\leq H\!D^{\rm sc}_{P,T}(\Sigma_0)$. Since~$H\!D^{\rm sc}_{P,T}(\Sigma_0)
\linebreak
=H\!D^{\rm sc}_{P,T}(S(\Sigma_0) V_0)\leq H\!D^{{\rm sh},S}_{P,T}(V_0)$, the result follows.
(ii) Before proceeding, note that since~$P$ is an elliptical probability measure with location~$\theta_0$, the affine-equivariance of~$T$ implies that~$T_P=\theta_0$. Ellipticity further entails that~$s_{P,T}=P[\{\theta_0\}]$. Moreover, we must have~$s_{P,T}<1/2$ (otherwise, $P[|Z_1|=0]\geq P[\{\theta_0\}]=s_{P,T}\geq 1/2$, so that~$0\in \mathcal{I}_{\rm MSD}[Z_1]$, a contradiction). Now, assume, ad absurdum, that there exists~$V\in\mathcal{P}_k^{S}\setminus\{V_0\}$ with~$H\!D^{{\rm sh},S}_{P,T}(V)=H\!D^{{\rm sh},S}_{P,T}(V_0)$. Since~$H\!D^{{\rm sh},S}_{P,T}(V_0)=H\!D^{\rm sc}_{P,T}(\Sigma_0)\geq 1/2$, we must have~$H\!D^{{\rm sh},S}_{P,T}(V)\geq 1/2$. Since~$\alpha_{P,T}=s_{P,T}<1/2$, there exists~$\sigma^2>0$ such that~$H\!D^{\rm sc}_{P,T}(\sigma^2 V)\geq s_P$. Therefore, Theorem~\ref{maxsigma2shape} ensures that~$H\!D^{{\rm sh},S}_{P,T}(V)=H\!D^{\rm sc}_{P,T}(\sigma^2_V V)$  for some~$\sigma^2_V>0$. We therefore have
$$
H\!D^{\rm sc}_{P,T}(\sigma^2_V V)=H\!D^{{\rm sh},S}_{P,T}(V)=H\!D^{{\rm sh},S}_{P,T}(V_0)=H\!D^{\rm sc}_{P,T}(\Sigma_0)=H\!D^{\rm sc}_{P,T}(S(\Sigma_0) V_0)
.
$$
 Theorem~\ref{Fishconst}(ii) then yields that~${\rm Sp}( (S(\Sigma_0))^{-1/2}\sigma_V V_0^{-1/2} V^{1/2})\subset \mathcal{I}_{\rm MSD}[Z_1]$. Since, by assumption, $\mathcal{I}_{\rm MSD}[Z_1]=\{1\}$, $V_0^{-1/2} V^{1/2}$ must be proportional to the identity matrix, which implies that~$V=V_0$. 
\cqfd
\vspace{3mm}
    

{\sc Proof of Theorem~\ref{Propquasiconcavityshape}.}
(i) Consider the scale functional~$S_{\rm tr}$ and \mbox{fix~$\varepsilon>0$.} By definition, there exist positive real numbers~$\sigma^2_{a}$ and~$\sigma^2_{b}$ such that
$$
H\!D^{\rm sc}_{P,T}(\sigma^2_{a} V_a)\geq H\!D^{{\rm sh},S_{\rm tr}}_{P,T}(V_a)-\varepsilon
\ \textrm{ and } \
 H\!D^{\rm sc}_{P,T}(\sigma^2_{b} V_b)\geq H\!D^{{\rm sh},S_{\rm tr}}_{P,T}(V_b)-\varepsilon.
 $$
  Consider then the linear path~$\Sigma_t=(1-t)\Sigma_a+t \Sigma_b$ from~$\Sigma_a=\sigma^2_{a} V_a$ to~$\Sigma_b=\sigma^{2}_{b} V_b$. Letting~$h(t)=t\sigma^2_b/((1-t)\sigma^{2}_a+t\sigma^{2}_b)$, the $S_{\rm tr}$-shape matrix associated with~$\Sigma_t$ is
\begin{eqnarray*}
\frac{k}{{\rm tr}[\Sigma_{t}]} 
\, 
\Sigma_{t}
&\!\!\!=\!\!\!&
\frac{k}{{\rm tr}[ (1-t)\sigma^2_a V_a+ t\sigma^2_b V_b]}
\,
( (1-t)\sigma^2_a V_a+ t\sigma^2_b V_b)
\\[2mm]
&\!\!\!=\!\!\! & 
\frac{k}{{\rm tr}[ (1-h(t))V_a+ h(t) V_b]}
\,
( (1-h(t)) V_a+ h(t) V_b)
=
V_{h(t)}
.
\end{eqnarray*}
Since~$h:[0,1]\to[0,1]$ is a one-to-one mapping, Theorem~\ref{Propquasiconcavity} yields that
\begin{eqnarray*}
\lefteqn{
H\!D^{{\rm sh},S_{\rm tr}}_{P,T}(V_t,P)
=
H\!D^{{\rm sh},S_{\rm tr}}_{P,T}\Big(
\frac{k}{{\rm tr}[\Sigma_{h^{-1}(t)}]} \, \Sigma_{h^{-1}(t)}
,P\Big)
\geq 
H\!D^{\rm sc}_{P,T}( \Sigma_{h^{-1}(t)})
}
\\[2mm]
& & 
\hspace{-5mm} 
\geq 
\min(H\!D^{\rm sc}_{P,T}(\Sigma_a),H\!D^{\rm sc}_{P,T}(\Sigma_b))
\geq
\min(H\!D^{{\rm sh},S_{\rm tr}}_{P,T}(V_a),H\!D^{{\rm sh},S_{\rm tr}}_{P,T}(V_b))	
-\varepsilon,
\end{eqnarray*}
for any~$t\in [0,1]$. Since this holds for any~$\varepsilon>0$, Part~(i) of the result is proved for~$S=S_{\rm tr}$. The proof for~$S=S_{11}$ is along the exact same lines, hence is omitted. As for Part~(ii), it strictly follows like Part~(ii) of Theorem~\ref{propgeodesicquasiconc}. 
\cqfd
\vspace{3mm}


{\sc Proof of Theorem~\ref{propgeodesicquasiconcshape}}. 
(i) Fix~$\varepsilon>0$. By definition, there exist~$\sigma^2_{a}>0$ and~$\sigma^2_{b}>0$ such that
$$
H\!D^{\rm sc}_{P,T}(\sigma^2_{a} V_a)\geq H\!D^{{\rm sh},S_{\det}}_{P,T}(V_a)-\varepsilon
\ \textrm{  and } \
H\!D^{\rm sc}_{P,T}(\sigma^2_{b} V_b)\geq H\!D^{{\rm sh},S_{\det}}_{P,T}(V_b)-\varepsilon
.
$$
 Consider then the geodesic path~$\tilde{\Sigma}_t=
\Sigma_a^{1/2} 
\big( \Sigma_a^{-1/2} \Sigma_b \Sigma_a^{-1/2} \big)^t 
\Sigma_a^{1/2}$ from~$\Sigma_a=\sigma^2_{a} V_a$ to~$\Sigma_b=\sigma^2_{b} V_b$. Then, since~$\det \tilde{\Sigma}_t=(\det \Sigma_a)^{1-t}(\det \Sigma_b)^t$, it is easy to check that the $S_{\det}$-shape matrix associated with~$\tilde{\Sigma}_t$ is
$
(\det \tilde{\Sigma}_t)^{-1/k} \tilde{\Sigma}_t
=
V_a^{1/2} 
\big( V_a^{-1/2} V_b V_a^{-1/2} \big)^t 
V_a^{1/2}
=:
\tilde{V}_t
.
$
Therefore, using Theorem~\ref{propgeodesicquasiconc}, we obtain
\begin{eqnarray*}
\lefteqn{
H\!D^{{\rm sh},S_{\det}}_{P,T}(\tilde{V}_t)
\geq 
H\!D^{\rm sc}_{P,T}((\det \tilde{\Sigma}_t)^{1/k} \tilde{V}_t)
=
H\!D^{\rm sc}_{P,T}(\tilde{\Sigma}_t)
}
\\[2mm]
& & 
\hspace{-3mm} 
\geq 
\min(H\!D^{\rm sc}_{P,T}(\Sigma_a),H\!D^{\rm sc}_{P,T}(\Sigma_b))
\geq
\min(H\!D^{{\rm sh},S_{\det}}_{P,T}(V_a),H\!D^{{\rm sh},S_{\det}}_{P,T}(V_b))	
-\varepsilon.
\end{eqnarray*}
Part~(i) of the result follows since~$\varepsilon>0$ is arbitrary above. As for Part~(ii), it is obtained again as in Part~(ii) of Theorem~\ref{propgeodesicquasiconc}. 
\cqfd
\vspace{3mm}


{\sc Proof of Theorem~\ref{Propquasiconcavityharmonicshape}.}
(i) Fix~$\varepsilon>0$. By definition, there exist positive real numbers~$\sigma^2_{a}$ and~$\sigma^2_{b}$ such that
$$
H\!D^{\rm sc}_{P,T}(\sigma^2_{a} V_a)\geq H\!D^{{\rm sh},S_{\rm tr}^*}_{P,T}(V_a)-\varepsilon
\ \textrm{ and } \
H\!D^{\rm sc}_{P,T}(\sigma^2_{b} V_b)\geq H\!D^{{\rm sh},S_{\rm tr}^*}_{P,T}(V_b)-\varepsilon
.
$$
 Consider then the harmonic path~$\Sigma_t^*=((1-t)\Sigma^{-1}_a+t \Sigma^{-1}_b)^{-1}$ from~$\Sigma_a=\sigma^2_{a} V_a$ to~$\Sigma_b=\sigma^{-2}_{b} V_b$. Then, letting~$h(t)=t\sigma^{-2}_b/((1-t)\sigma^{-2}_a+t\sigma^{-2}_b)$, the $S_{\rm tr}^*$-shape matrix associated with~$\Sigma_t^*$ is
\begin{eqnarray*}
\frac{{\rm tr}[(\Sigma_t^*)^{-1}]}{k} 
\, 
\Sigma_t^*
&=&
((1-t)\sigma^{-2}_a+ t\sigma^{-2}_b)
\,
((1-t)\sigma^{-2}_a V_a^{-1}+ t\sigma^{-2}_b V_b^{-1})^{-1}
\\[2mm]
& = &
\big(
(1-h(t))
 V_a^{-1}
 + 
h(t)
 V_b^{-1}
\big)^{-1}
=:
V^*_{h(t)}.
\end{eqnarray*}
Since~$h:[0,1]\to[0,1]$ is a one-to-one mapping, we obtain that, for any~$t\in [0,1]$,
\begin{eqnarray*}
\lefteqn{
H\!D^{{\rm sh},S_{\rm tr}^*}_{P,T}(V_t^*)
=
H\!D^{{\rm sh},S_{\rm tr}^*}_{P,T}
\bigg(
\,
\frac{{\rm tr}[(\Sigma^*_{h^{-1}(t)})^{-1}]}{k} 
\, 
\Sigma^*_{h^{-1}(t)}
\bigg)
\geq 
H\!D^{\rm sc}_{P,T}( \Sigma^*_{h^{-1}(t)} )
}
\\[2mm]
& &  
\hspace{-3mm} 
\geq 
\min(H\!D^{\rm sc}_{P,T}(\Sigma_a),H\!D^{\rm sc}_{P,T}(\Sigma_b))
\geq
\min(H\!D^{{\rm sh},S_{\rm tr}^*}_{P,T}(V_a),H\!D^{{\rm sh},S_{\rm tr}^*}_{P,T}(V_b))	
-\varepsilon
.
\end{eqnarray*}
Since this holds for any~$\varepsilon>0$, Part~(i) of the result is proved. Part~(ii) strictly follows like Part~(ii) of Theorem~\ref{propgeodesicquasiconc}. 
\cqfd
\vspace{3mm}





\bibliographystyle{imsart-nameyear.bst} 
\bibliography{ScatterDepth.bib}           
\vspace{3mm} 


\end{document}